\documentclass{article}
\usepackage{amssymb}
\usepackage{amsmath}
\usepackage{graphicx}
\usepackage{color}
\usepackage{multicol}

\usepackage{float}

\newtheorem{theorem}{Theorem}
\newtheorem{definition}{Definition}
\newtheorem{lemma}[theorem]{Lemma}
\newtheorem{rem}[theorem]{Remark}
\newtheorem{pro}[theorem]{Proposition}
\newtheorem{con}[theorem]{Conjecture}

\newtheorem{col}[theorem]{Collorary}

\newcommand{\inter}{\mathrm{int}}

\newcommand{\dom }{\,{\rm dom}\,}

\newcommand{\im }{\mbox{Im} \,}

\newcommand{\cover}[1]{\stackrel{#1}{\Longrightarrow}}

\newcommand{\invcover}[1]{\stackrel{#1}{\Longleftarrow}}

\newcommand{\W}{\mathcal{T}}
\newcommand{\inc}{\textrm{inc}}
\newcommand{\out}{\textrm{out}}

\newcommand{\outcomment}[1]{}

\def\qed{{\hfill{\vrule height5pt width3pt depth0pt}\medskip}}

\begin{document}
\begin{center} {\bf \LARGE  Shadowing of non-transversal heteroclinic chains }\\
 \vskip 0.5cm


\vskip 0.5cm

{\large Amadeu Delshams\footnote{Partially supported by the MINECO-FEDER Grant MTM2015-65715-P and the Russian
Scientific Foundation Grant 14-41-00044 at the Lobachevsky University of
Nizhny Novgorod.}, Adri\`a Simon\footnote{Partially supported by
the MINECO-FEDER Grant MTM2012-31714 and the Catalan
Grant 2014SGR504.}} \\
Departament de Matem\`atiques, Universitat Polit\`ecnica de Catalunya \\
Av. Diagonal 647, 08028 Barcelona \\
Amadeu.Delshams@upc.edu, adria.simon@gmail.com
 \vskip\baselineskip

   {\large Piotr Zgliczy\'nski}\footnote{Research has been supported by Polish National Science Centre grant 2011/03B/ST1/04780}     \\
 Jagiellonian University, Institute of Computer Science and Computational Mathematics, \\
{\L}ojasiewicza 6, 30--348  Krak\'ow, Poland \\ e-mail:
umzglicz@cyf-kr.edu.pl

\vskip 0.5cm
 \today

\end{center}

\begin{abstract}
We present a new result about the shadowing of nontransversal chain of heteroclinic connections based on the idea
of dropping dimensions. We illustrate this new mechanism with several examples.
As an application we discuss this mechanism in a simplification of a toy model system derived by Colliander \emph{et al.}
in the context of cubic defocusing nonlinear Schr\"odinger equation.
\end{abstract}

\section{Introduction}

In the present paper we deal with the problem of shadowing a nontransversal chain of heteroclinic connections between invariant sets
(fixed points, periodic orbits, etc). The motivation for us is the work \cite{CK} (see also \cite{GK}) on the transfer of energy to high frequencies in the nonlinear  Schrodinger equation (just NLS from now on).
From the dynamical systems viewpoint there is one remarkable feature of the construction in \cite{CK}, namely
that the authors were able to shadow a non-transversal highly degenerated chain of heteroclinic connections between
some periodic orbits. The length of the chain is arbitrary, but finite. Neither in \cite{CK} nor in \cite{GK} we were able to find a clear
geometric picture showing how this is achieved, so it could be easily applicable to other systems. In this work we present a mechanism, which we believe
gives a geometric explanation of what is happening. Moreover, we strive to establish an abstract framework, which will make it easier to apply this technique to other
systems, both PDEs and ODEs, in questions related to the existence of diffusing orbits.
The term \emph{diffusing orbit} relates to the Arnold's diffusion~\cite{Ar64}
for the perturbation of integrable Hamiltonian systems. Throughout the paper we
will often call diffusing orbit an orbit shadowing a  chain of heteroclinic connections,
and occasionally the existence of such an orbit will be referred to as the diffusion.

In our picture we think of evolving a disk of  dimension $k$  along a heteroclinic transition chain and when a given
transition is not transversal, then we `drop' one or more dimensions of our disk, i.e., we select a subdisk of lower dimension
``parallel to expanding directions in future transitions".  After at most $k$ transitions, our disk is a single point
and we cannot continue further. We will refer to this phenomenon as the \emph{dropping dimensions} mechanism.
Since this is a new mechanism, we have found it convenient to include several figures to illustrate the main
differences between transversal and non-transversal heteroclinic chains.
While thinking about disks has some geometric appeal,
we consider instead in our construction a thickened disk called h-set in the terminology of \cite{ZGi} and our approach is purely topological (just as the one presented in \cite{CK}).

The main technical tool used in our work is the notion of \emph{covering relations} as introduced in \cite{ZGi}, which differs from the notion used under the same name
in \cite{CK}. Similar ideas about the dropping exit dimensions appear implicitly also in the works \cite{BM+,WBS}.

In our work we present an abstract topological theorem about shadowing  chains of covering relations with dropping dimensions, and we show how such chains
of coverings can be obtained in the presence of chains of heteroclinic connections in two examples,
 a linear model 
and a  simplified Toy Model that the one in \cite{CK}, which however contains all the difficulties present in its prototype.
We intend to treat more complicated examples, in particular NLS from \cite{CK,GK} in subsequent papers.

The content of this paper can be described as follows. In Section~\ref{sec:NonTransverseDiffusion} we first describe the consequences of  the difference
between transversal and non-transversal intersection of invariant manifolds of fixed points and we present  the model problem with a
non-transversal heteroclinic chain and state our conjecture about the existence of shadowing orbits arbitrarily close to such chain.
We also introduce an example formed by a triangular system, where the existence of the diffusion is quite obvious.
In Section~\ref{sec:geomIdea} we explain the basic geometric idea of our dropping dimensions mechanism.
In Section~\ref{sec:covrel} we recall from \cite{ZGi} the notions of h-sets and the covering relation.
In Section~\ref{sec:topThm} we prove the main topological result on shadowing of chains of covering relations with dropping dimensions.
Using this new mechanism, in the next two sections we rigorously analyze two simple models,
a linear model in Section~\ref{sec:linModel-proof} 
and a simplified Toy Model in Section~\ref{sec:diffToymodel}.

\subsection{Notation}
\label{subsec:notation}

By $\mathbb{N}$, $\mathbb{Z}$, $\mathbb{Q}$, $\mathbb{R}$,
$\mathbb{C}$ we denote the set of natural, integer, rational, real
and complex numbers, respectively. We assume that $0 \in \mathbb{N}$. $\mathbb{Z}_-$ and
$\mathbb{Z}_+$ are nonpositive and nonnegative integers,  respectively.
By $S^1$ we will  denote the unit circle on the complex plane.

In $\mathbb{R}^n$ by $e_i$ for $i=1,\dots,n$ we will denote the $i$-th
vector from the canonical basis in $\mathbb{R}^{n}$, i.e. the $j$-th
coordinate of $e_i$ is equal to $1$, when $j=i$ and $0$ otherwise.

For $\mathbb{R}^n$  we will denote the norm of $x$ by $\|x\|$ and
when in some context the formula for the norm is not specified,
then it means that any norm can be used. For $x_0 \in
\mathbb{R}^s$, $B_s(x_0,r)=\{z \in \mathbb{R}^s : \|x_0
- z \| < r \}$ and $B_s=B_s(0,1)$.

Sometimes, if $V$ is a vector space with a norm,  then $B_V(a,r)$
will denote an open ball in $V$ centered at $a$ with radius $r$.

 For $z \in \mathbb{R}^u \times \mathbb{R}^s $ we will call
usually $x$ the first coordinate and $y$ the second one. Hence
$z=(x,y)$, where $x \in \mathbb{R}^u$ and $y \in \mathbb{R}^s$. We
will  use the projection maps $\pi_x(z)=x(z)=x$ and
$\pi_y(z)=y(z)=y$.  For functions $f: \mathbb{R}^u \times
\mathbb{R}^s \to \mathbb{R}^u \times \mathbb{R}^s$ we will use the
shortcuts $f_x=\pi_x f$ and $f_y=\pi_y f$.

Let $z \in \mathbb{R}^n$ and $U \subset \mathbb{R}^n$ be a compact
set and $f:U \to \mathbb{R}^n$ be continuous map, such that $z
\notin f(\partial U)$. Then the local Brouwer degree \cite{S} of
$f$ on $U$ at $z$ is defined and will be denoted by $\deg(f,U,z)$.
See for example the Appendix in \cite{ZGi} and references given there for
the properties of $\deg(f,U,z)$.

If $V,W$ are two vector spaces, then by $\mbox{Lin}(V,W)$ we will
denote the set of all linear maps from $V$ to $W$. When
$V=\mathbb{R}^k$ and $W=\mathbb{R}^m$, we will identify
$\mbox{Lin}(\mathbb{R}^k,\mathbb{R}^m)$ with the set of matrices
with $m$ columns and $k$ rows,  denoted by
$\mathbb{R}^{k \times m}$.


\section{Non-transverse diffusion, the statement of the problem, some examples}
\label{sec:NonTransverseDiffusion}

In this section we introduce the geometric assumptions under which we expect to  construct the orbits shadowing
a non-transversal heteroclinic chain. Our approach is motivated by the work \cite{CK} on the NLS.

One of the key ingredients of the constructions of the energy transfer in \cite{CK}  consists on finding an orbit which visits the neighborhoods of $N$ invariant 1-dimensional
objects in a $N$-dimensional complex system. Each object is connected with the previous and the following one with heteroclinic connections,
so the authors look for a solution that 'concatenates' these connections. This kind of scheme seems similar
to the Arnold diffusion\cite{Ar64}, but we plan to explain that it is a different phenomenon
since we do not have a transverse intersection between the invariant manifolds. In addition,
the proposed mechanism  could be applied to integrable systems in contrast to the Arnold diffusion, which is
a phenomenon that only takes place in non integrable systems.

\subsection{Transverse versus Non-Transverse}
\label{subsec:TransvsNoTrans}

In this subsection we will  explain the difference between the transverse and the non-transverse situation. The hint about the idea of dropping dimensions will be given.

To do so, we are going to consider a two dimensional map with four fixed points:
\[p_0=(0,0)\quad p_1=(1,0) \quad p_2=(1,1) \quad p_3=(2,1).\]

We are going to assume also that each point $p_i$ has a one dimensional stable manifold, $\W^{s}(p_i)$, and a one dimensional unstable manifold, $\W^{u}(p_i)$, both tangent to some linear subspaces. That is
\begin{itemize}
\item $\W^s(p_0)$ is tangent to the subspace generated by $e_2$ at $p_0$ and $\W^u(p_0)$ is tangent to the subspace generated by $e_1$ at $p_0$.
\item $\W^s(p_1)$ is tangent to the subspace generated by $e_1$ at $p_1$ and $\W^u(p_1)$ is tangent to the subspace generated by $e_2$ at $p_1$.
\item $\W^s(p_2)$ is tangent to the subspace generated by $e_2$ at $p_2$ and $\W^u(p_2)$ is tangent to the subspace generated by $e_1$ at $p_2$.
\item $\W^s(p_3)$ is tangent to the subspace generated by $e_1$ at $p_3$ and $\W^u(p_3)$ is tangent to the subspace generated by $e_2$ at $p_3$.
\end{itemize}

We are now going to consider two different scenarios. The first one consists on assuming that the unstable manifold of a point $p_i$ intersects transversally the stable manifold of the following point $p_{i+1}$.

The second will be given by a non-transverse intersection of the manifolds, with one of the  branches of $\W^u(p_i)$ coinciding with one of the branches of $\W^s(p_{i+1})$

The schematic situation is the following (unstable manifolds are in red and stable manifolds are in blue):
\begin{multicols}{2}
\begin{figure}[H]
\centering
	\includegraphics[width=\linewidth]%
	{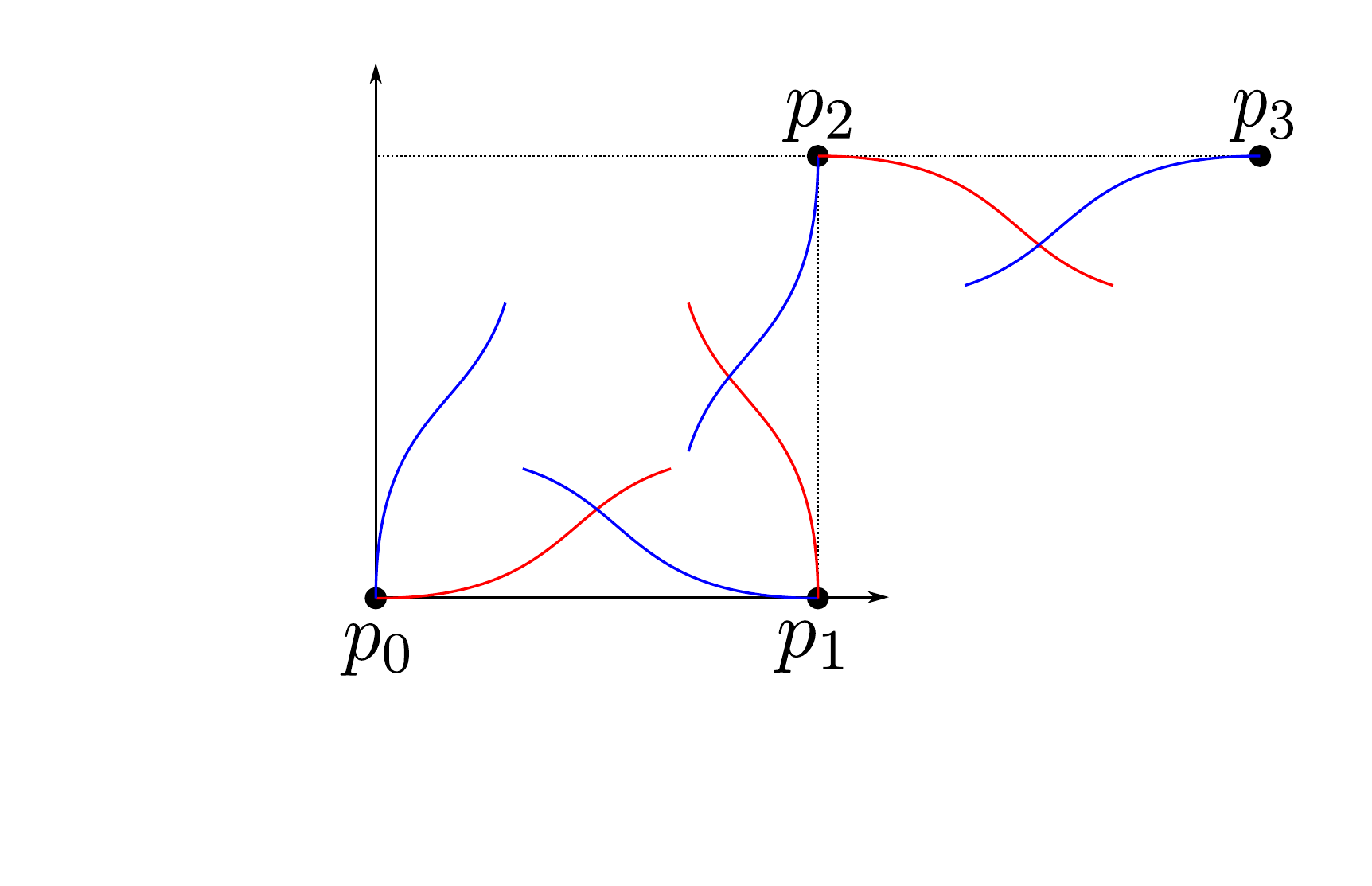}%
\end{figure}
\begin{figure}[H]
\centering
	\includegraphics[width=\linewidth]%
	{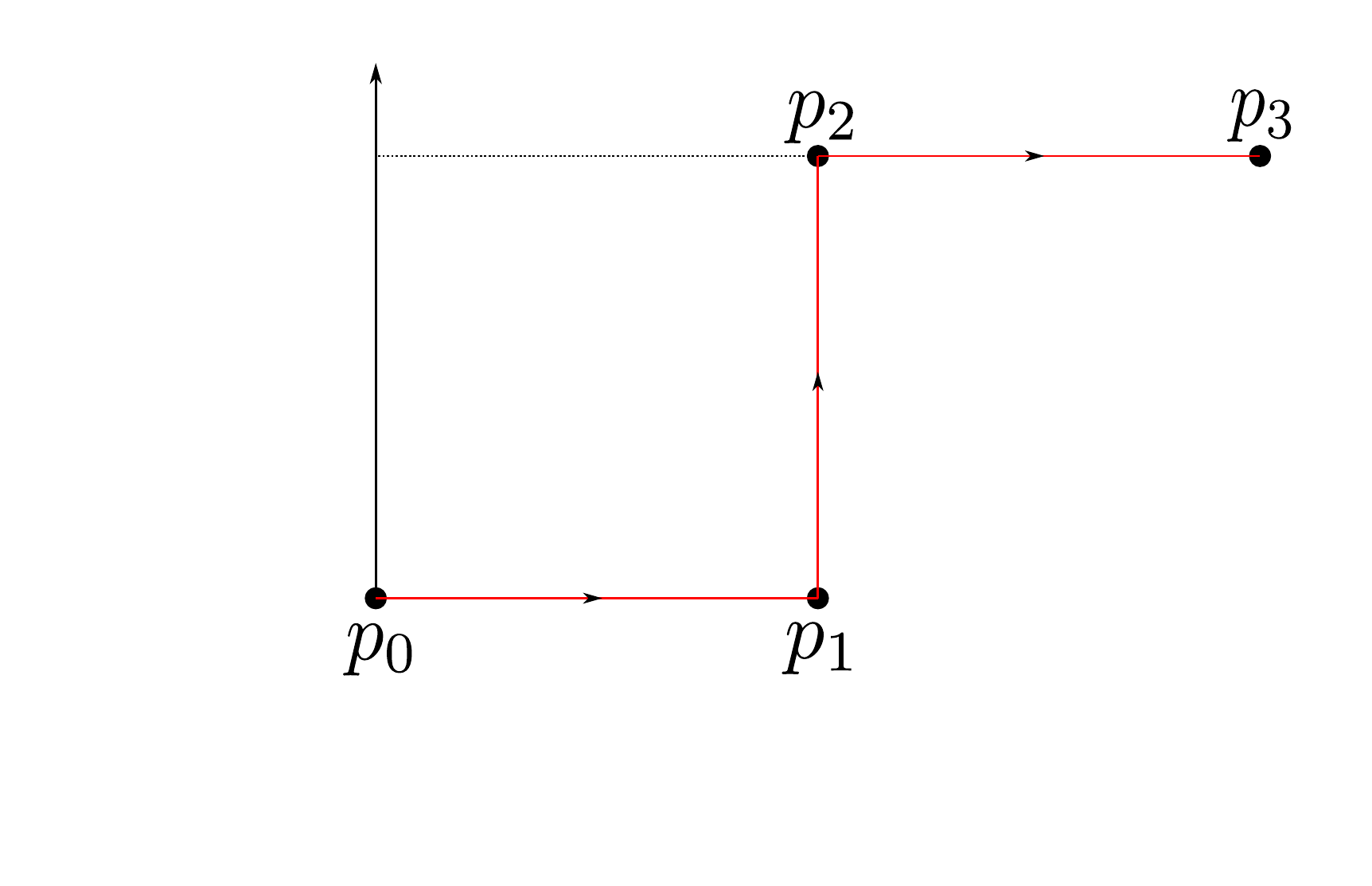}%
\end{figure}
\end{multicols}
We wonder if it is possible to connect $p_0$ with $p_3$ through the map, in both situations. To do so, we consider a ball containing the first fixed point $p_0$:
\begin{multicols}{2}
\begin{figure}[H]
\centering
	\includegraphics[width=\linewidth]%
	{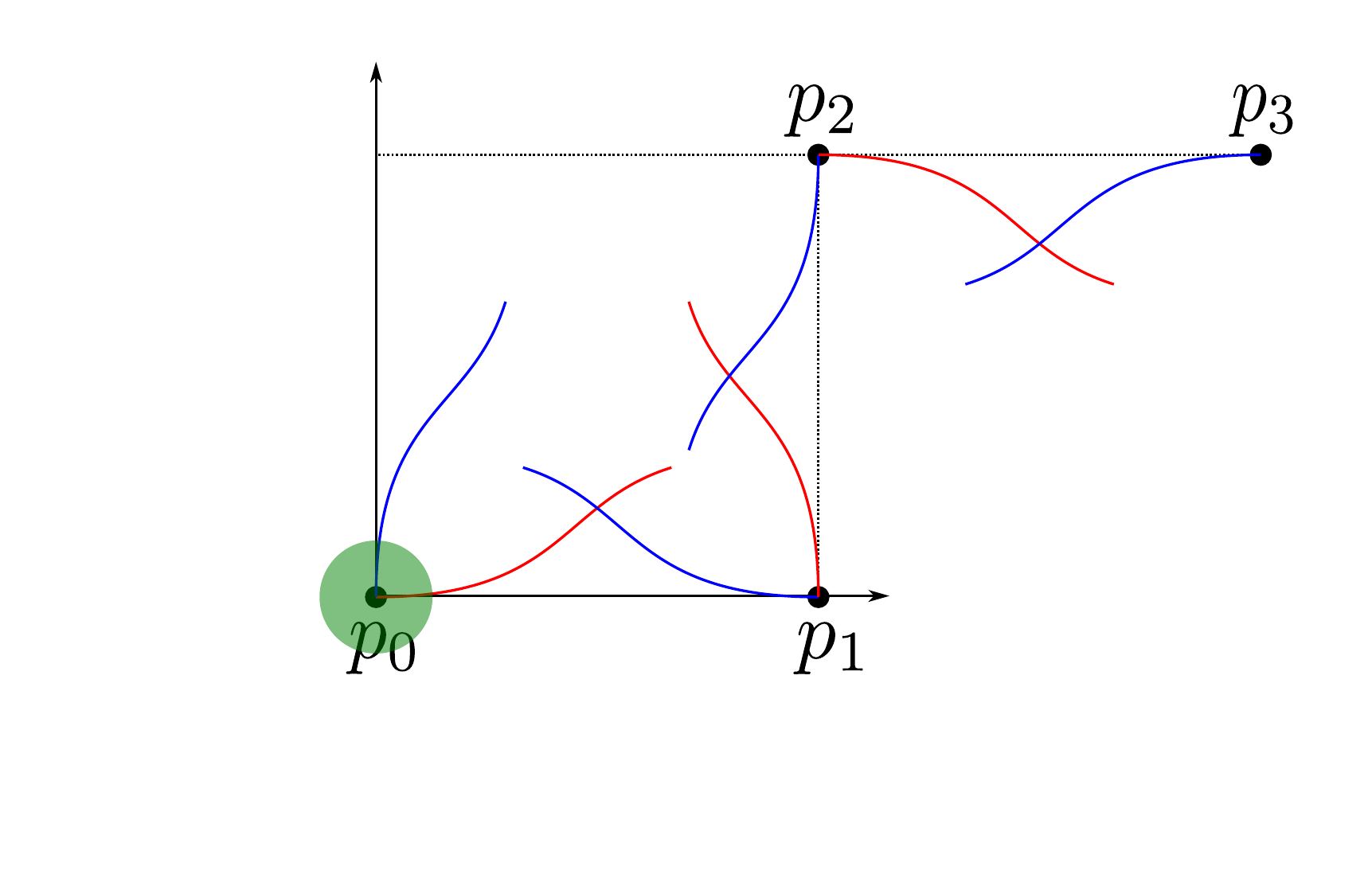}%
\end{figure}
\begin{figure}[H]
\centering
	\includegraphics[width=\linewidth]%
	{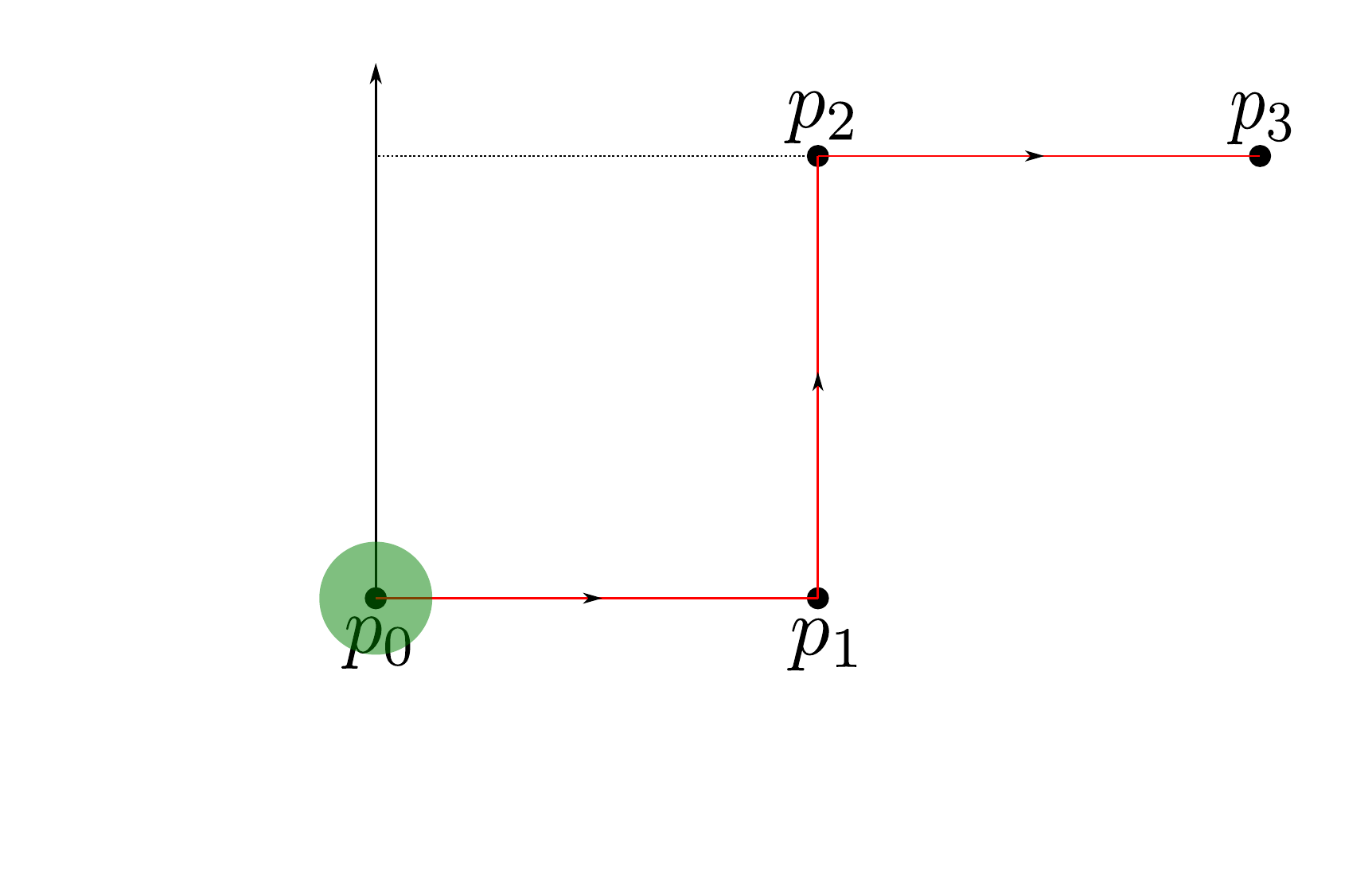}%
\end{figure}
\end{multicols}
If we compute iterates of the ball through the maps we can expect that it is expanded in the unstable direction and contracted in the stable direction:
\begin{multicols}{2}
\begin{figure}[H]
\centering
	\includegraphics[width=\linewidth]%
	{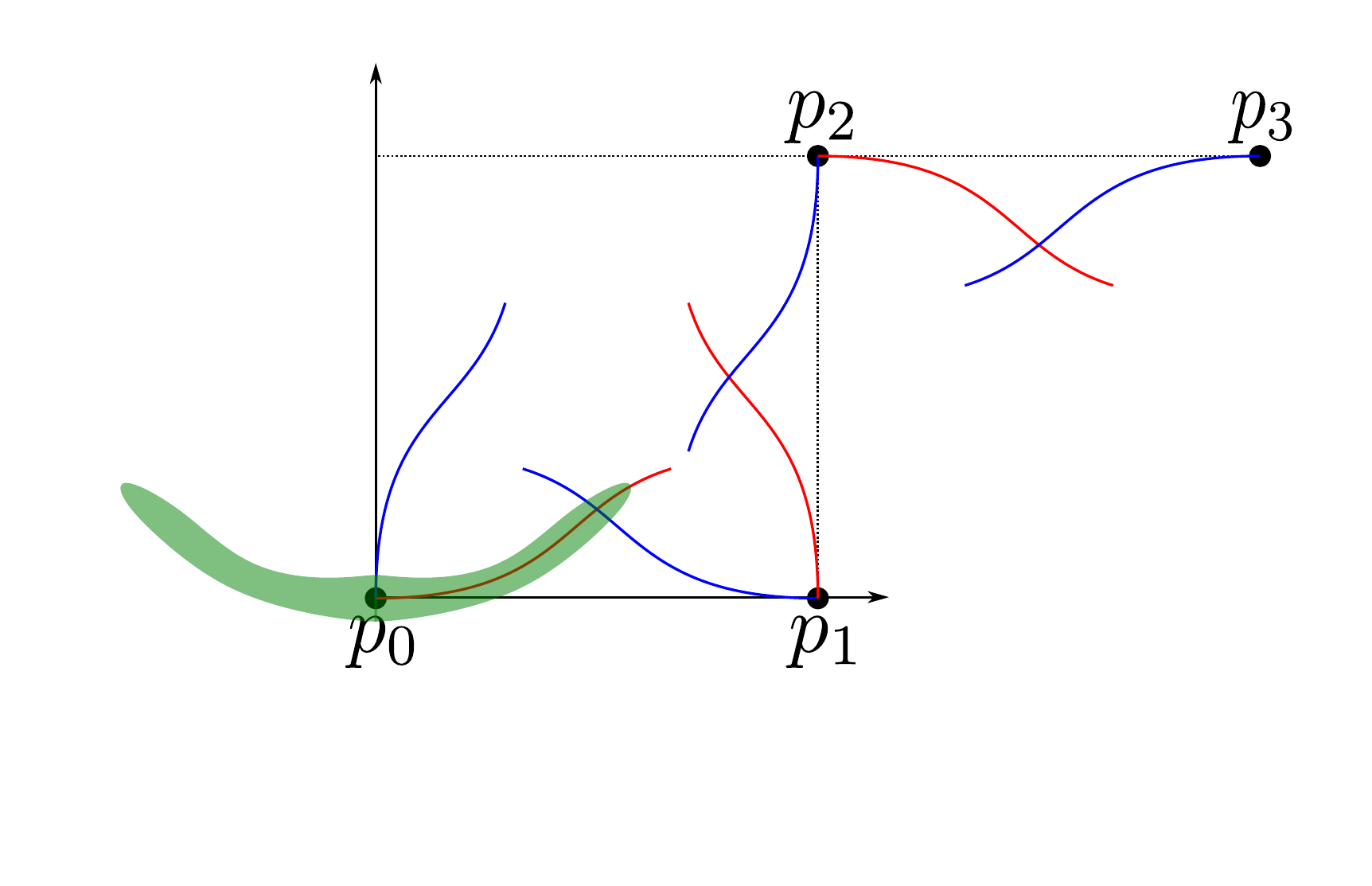}%
\end{figure}
\begin{figure}[H]
\centering
	\includegraphics[width=\linewidth]%
	{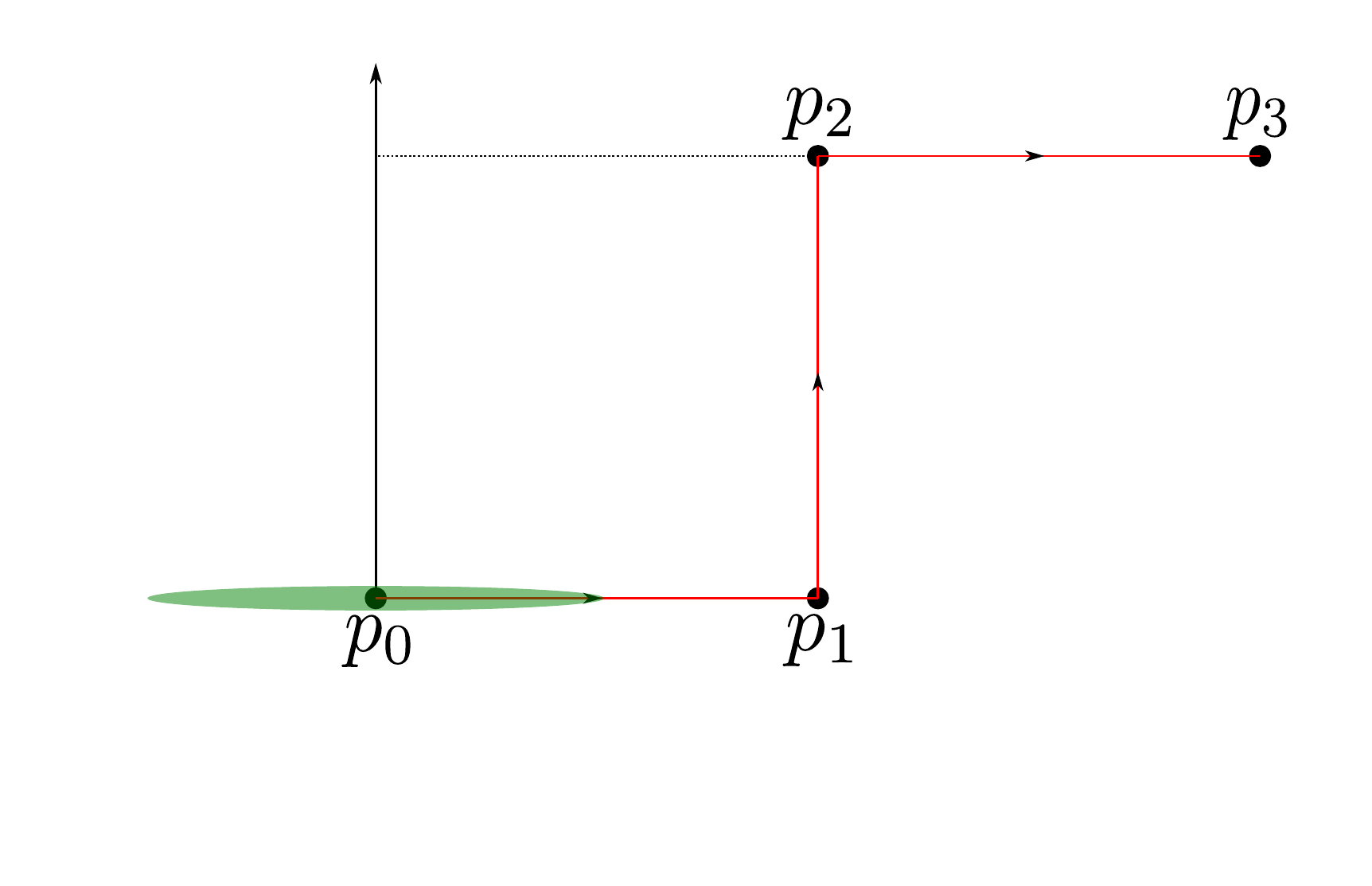}%
\end{figure}
\end{multicols}
Notice that in both cases the domain intersects the stable manifold of the following fixed point, $p_1$. In the transverse case the domain contains a heteroclinic point. In the non-transverse situation this is obvious because the manifolds are coincident. We can now restrict our domain precisely around that intersection point for the transverse case and at some place in the right-hand side of the fixed point $p_0$ for the non-transverse case.
\begin{multicols}{2}
\begin{figure}[H]
\centering
	\includegraphics[width=\linewidth]%
	{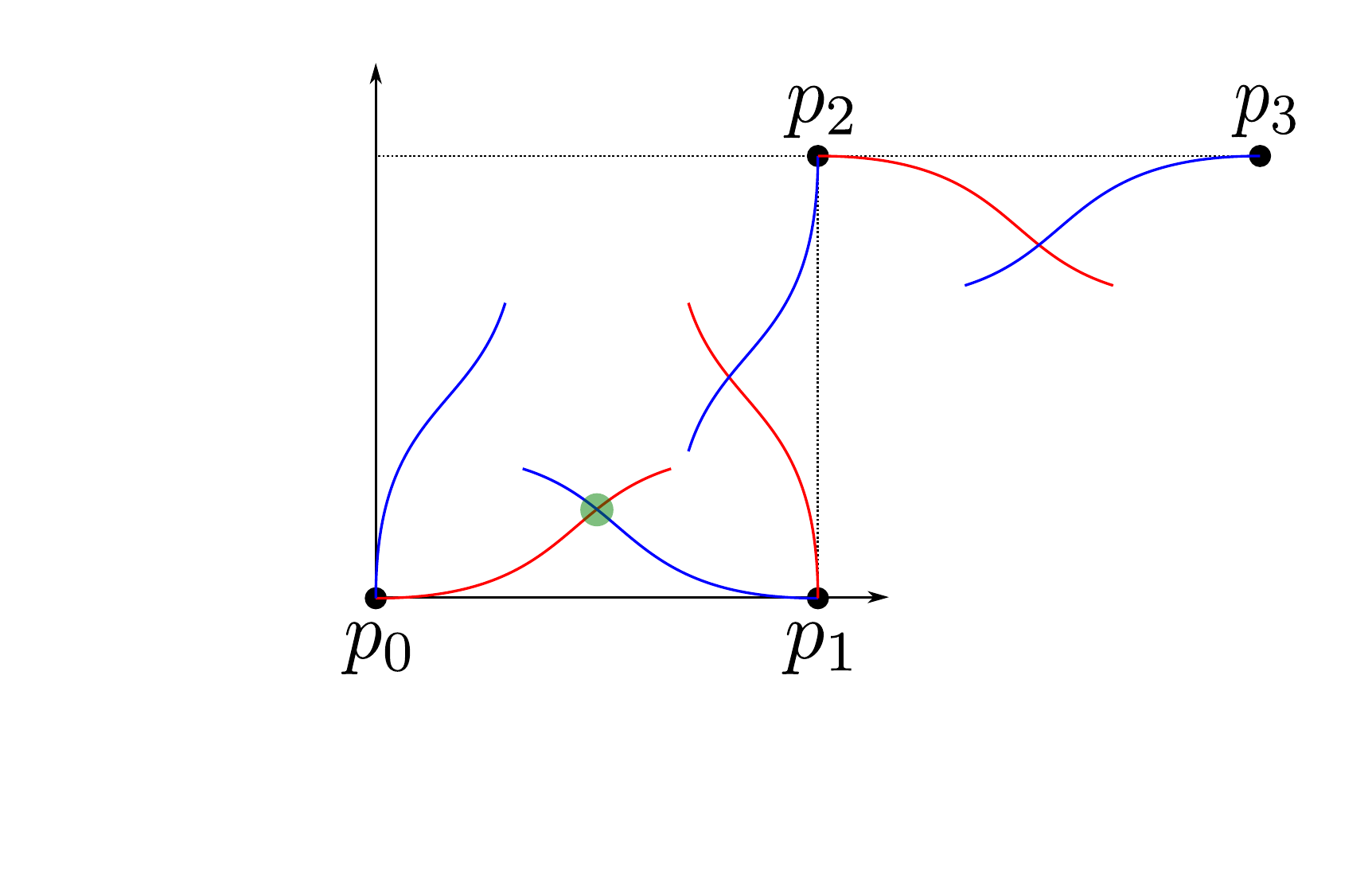}%
\end{figure}
\begin{figure}[H]
\centering
	\includegraphics[width=\linewidth]%
	{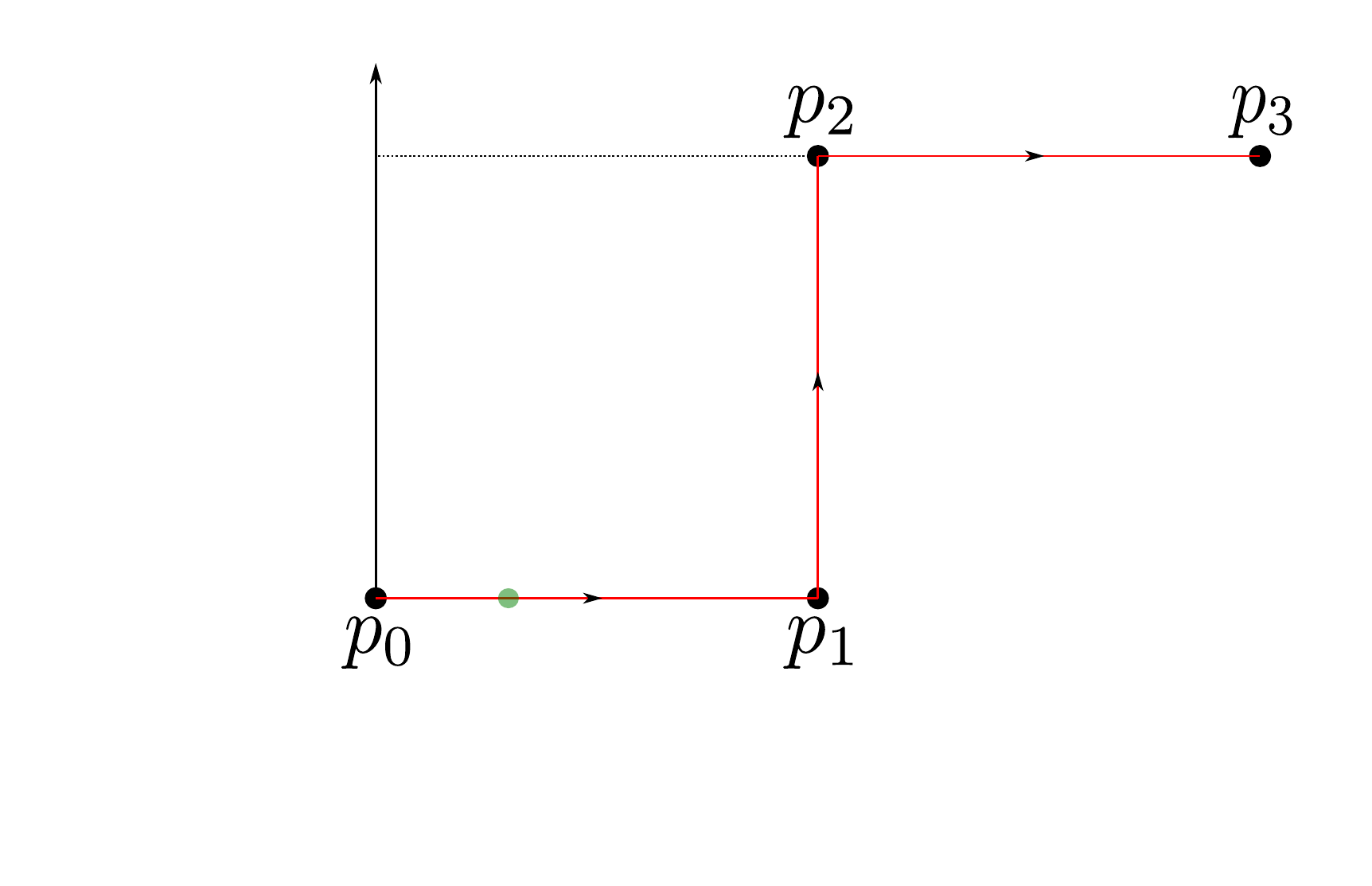}%
\end{figure}
\end{multicols}
The forward iterates of the restricted domain approach the fixed point $p_1$ since, in both cases, our domains contain heteroclinic points. The domains will not only approach $p_1$ but also, after some iterates, will spread along the unstable manifold of $p_1$:

\begin{multicols}{2}
\begin{figure}[H]
\centering
	\includegraphics[width=\linewidth]%
	{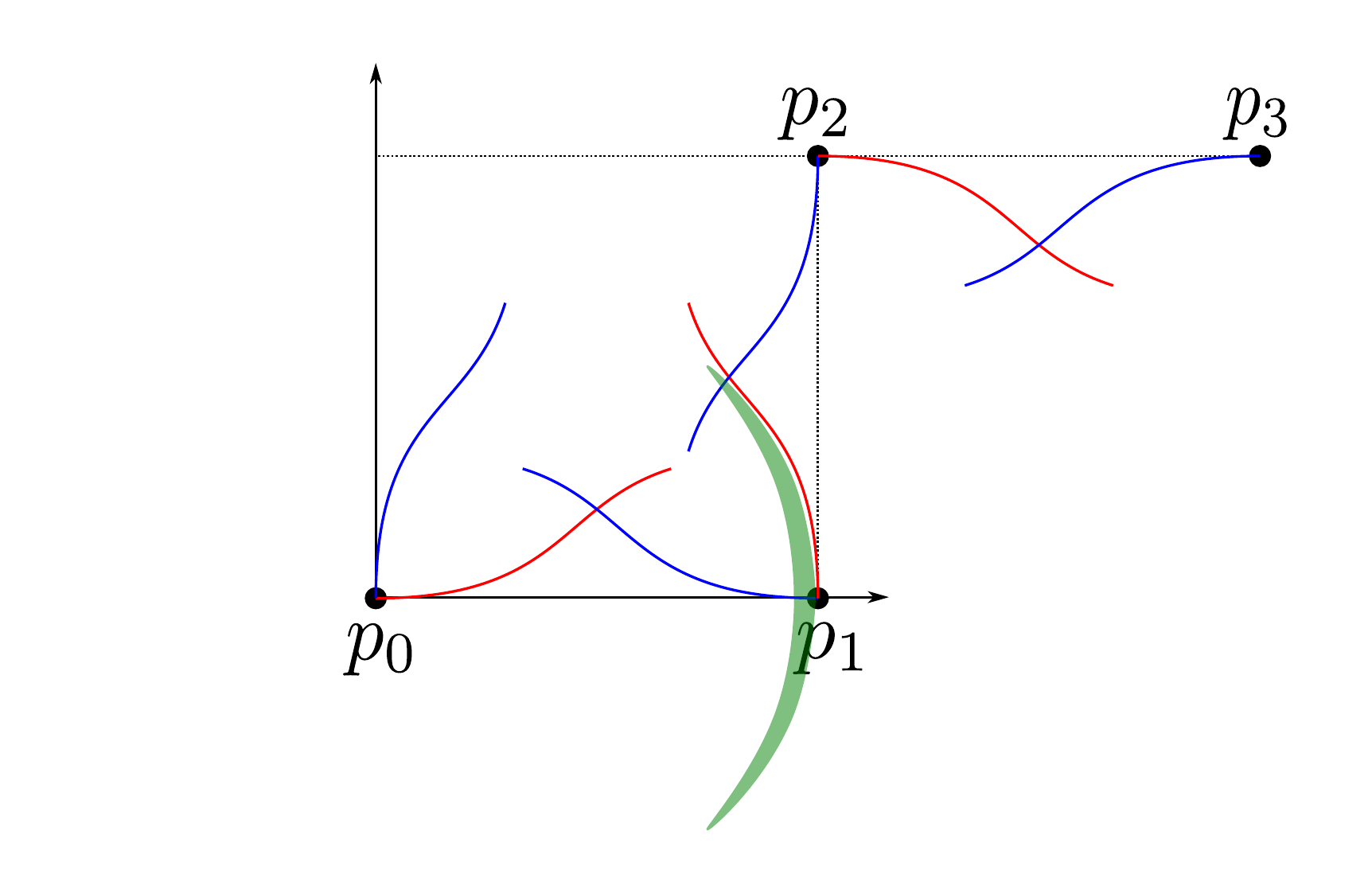}%
\end{figure}
\begin{figure}[H]
\centering
	\includegraphics[width=\linewidth]%
	{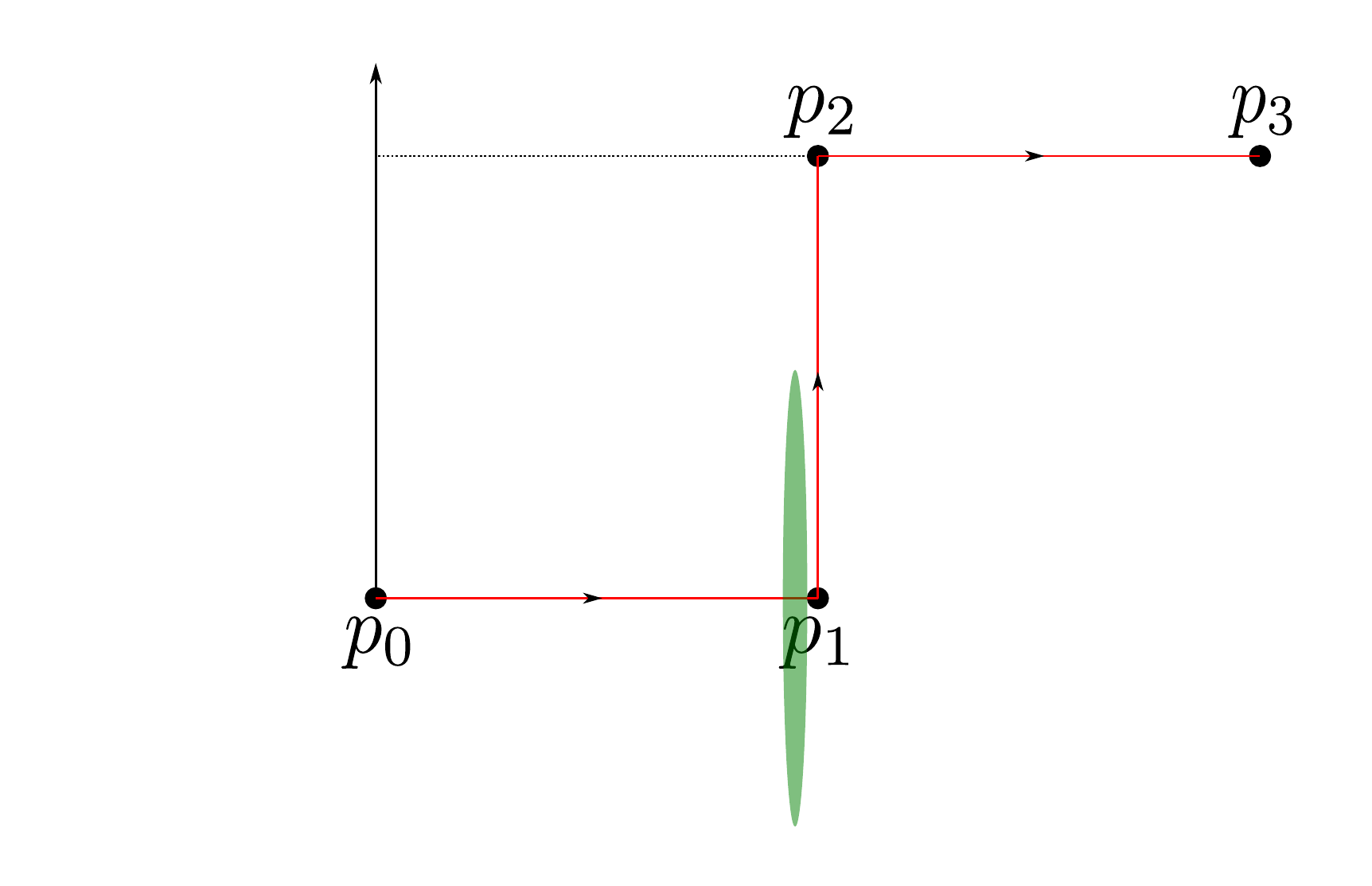}%
\end{figure}
\end{multicols}
While our domain spreads along the unstable manifold of $p_1$, here we find the first big difference between the transverse and the non-transverse case.  In the transverse case,  it is clear that our domain will intersect the stable manifold of $p_2$. In the non-transverse situation, our domain will never cross the stable manifold of $p_2$. Then, we restrict our domain in the intersection for the transverse case and in the upper part of $p_1$ since we want to reach $p_3$:

\begin{multicols}{2}
\begin{figure}[H]
\centering
	\includegraphics[width=\linewidth]%
	{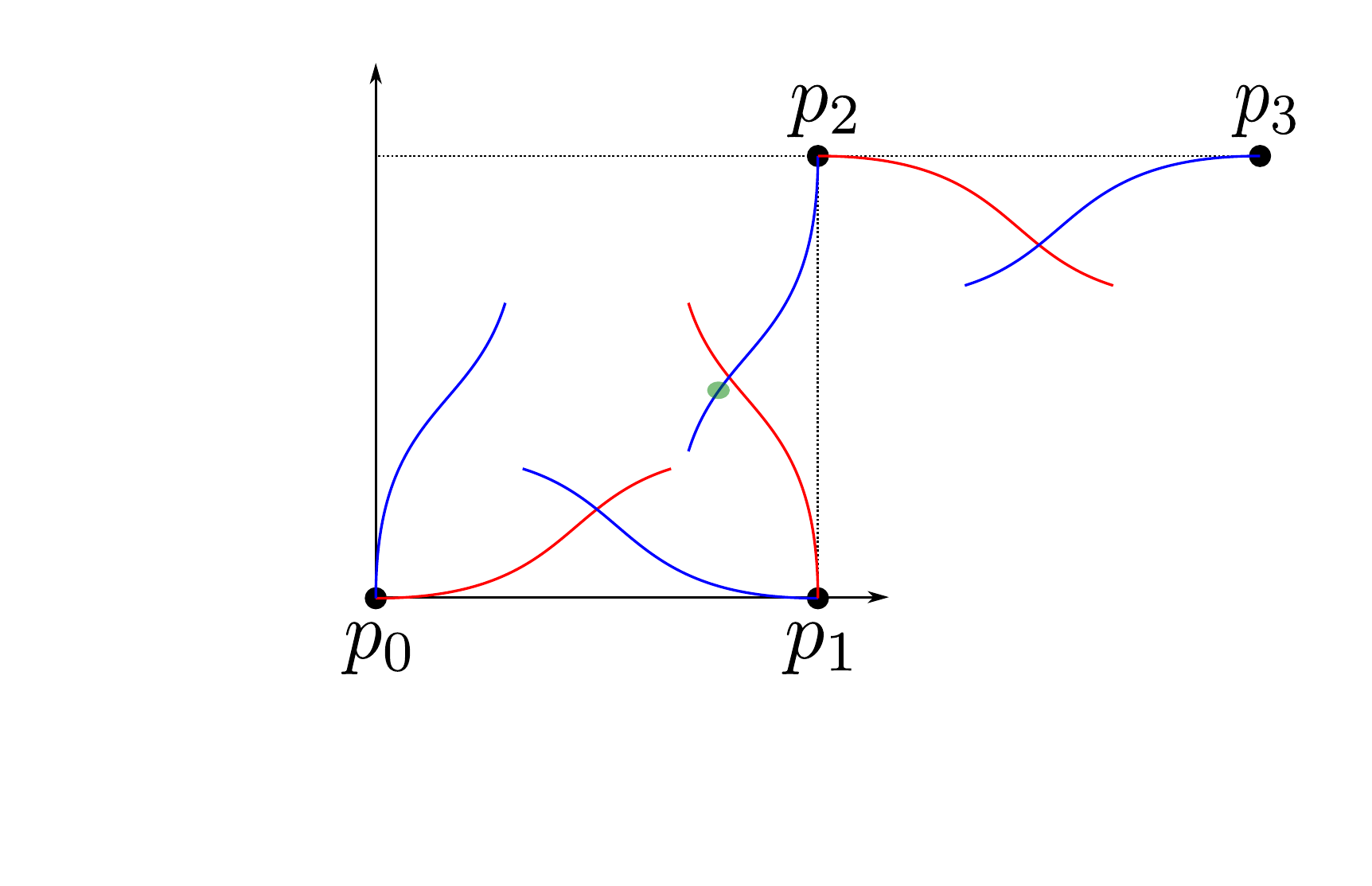}%
\end{figure}
\begin{figure}[H]
\centering
	\includegraphics[width=\linewidth]%
	{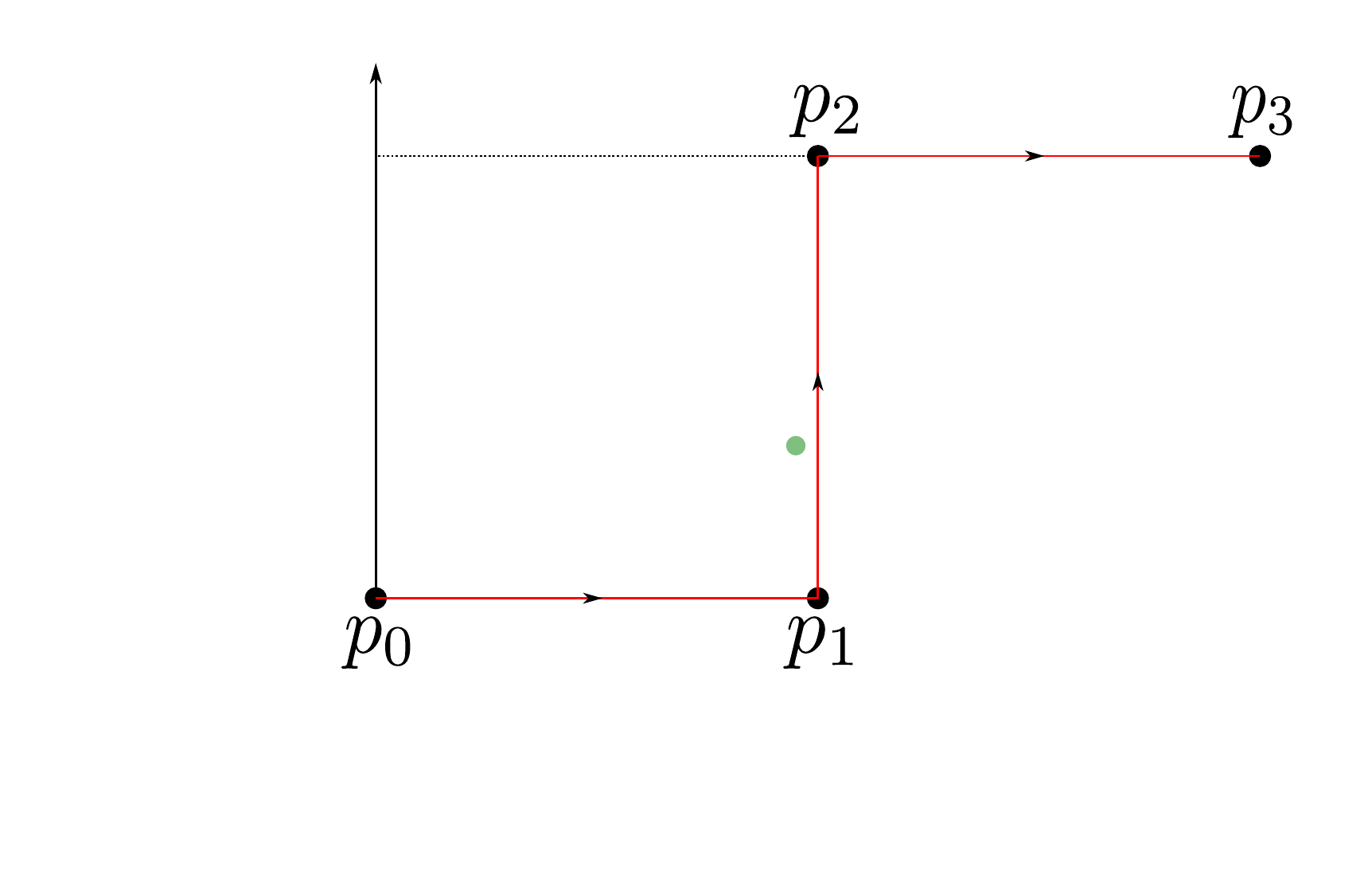}%
\end{figure}
\end{multicols}
Using, in the transverse case, the same argument as before, since our domain contains a heteroclinic point in the transverse situation, forward iterates will spread our domain along the unstable manifold of $p_2$. For the non transverse case we will reach the proximity of $p_2$ after some iterates, but our domain will be trapped and will not visit the following fixed point, $p_3$:
\begin{multicols}{2}
\begin{figure}[H]
\centering
	\includegraphics[width=\linewidth]%
	{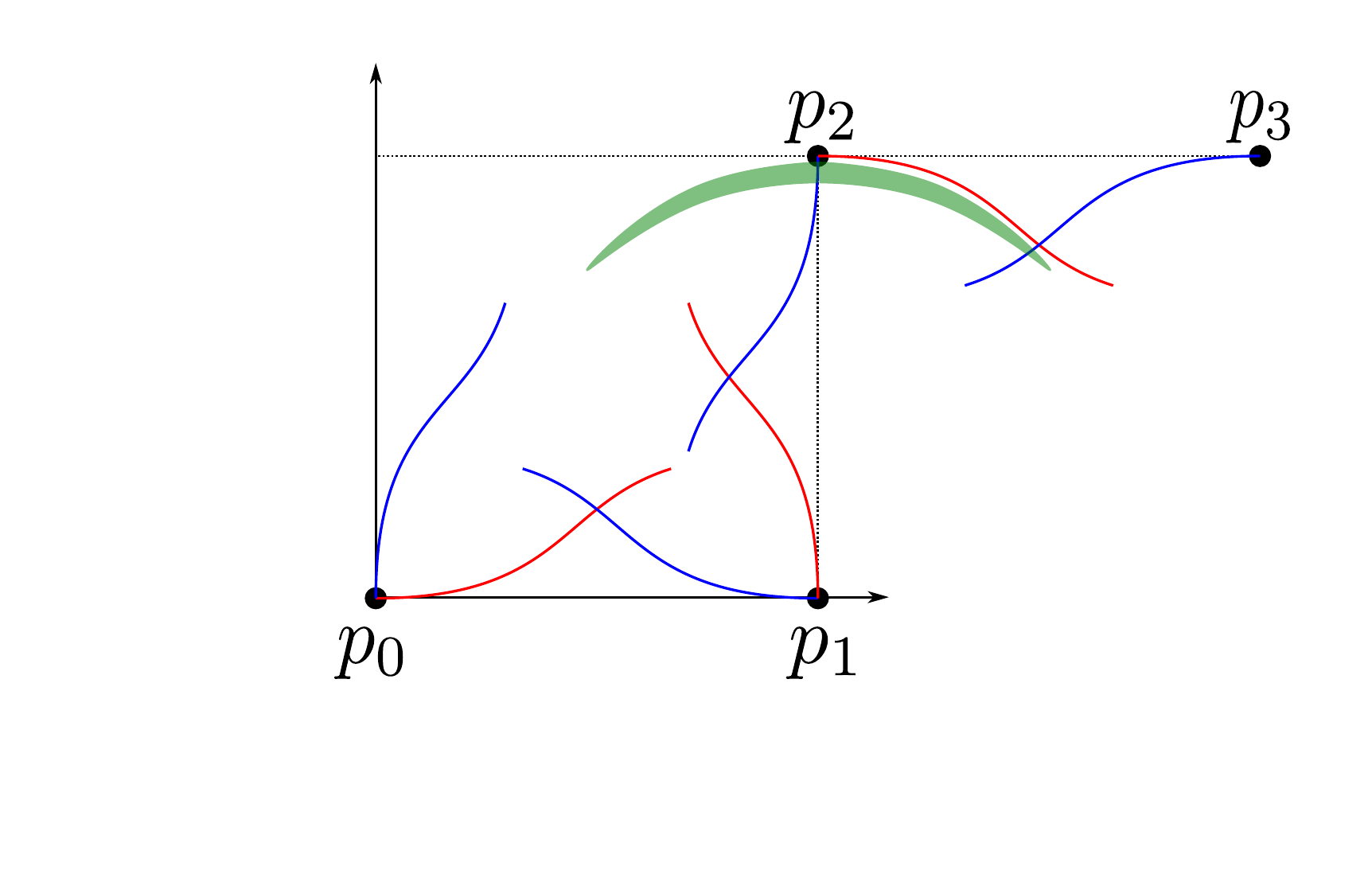}%
\end{figure}
\begin{figure}[H]
\centering
	\includegraphics[width=\linewidth]%
	{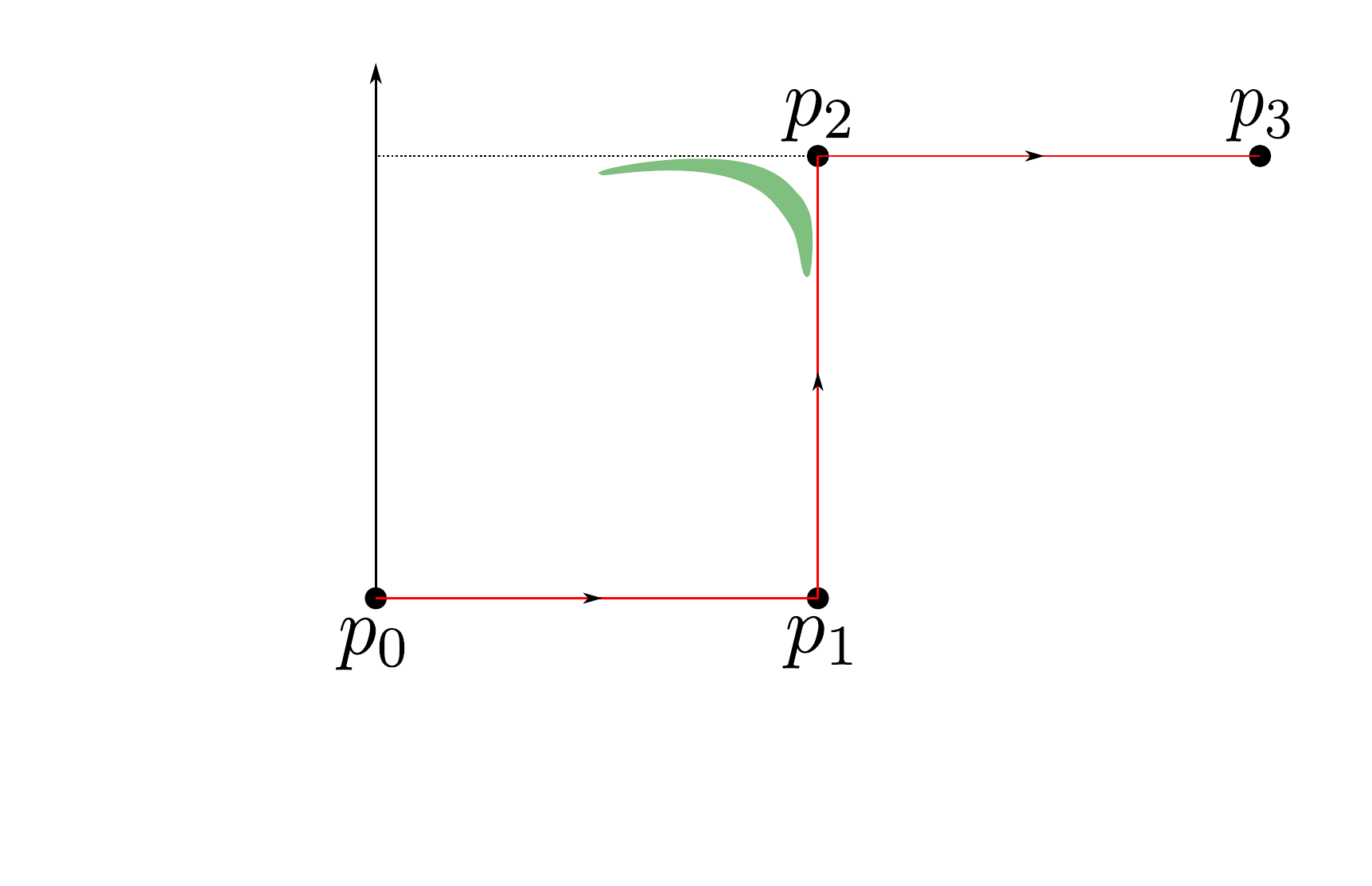}%
\end{figure}
\end{multicols}
In the transverse case, we could continue and see that the domain will visit $p_3$.

%

The above analysis, shows that, on the one hand, in the transverse situation there are no geometric obstructions in shadowing the heteroclinic chain. On the other hand, in the non-transverse case, we can see that, in general, we cannot visit as many invariant objects as we want. So, now, we wonder why the authors of \cite{CK} can connect $N$ periodic orbits in the Toy Model System. The main reason is the large dimension of the system and the fact that each connection takes place in a direction that has not been used in the past.

Indeed notice that, if in the non-transverse example the last fixed point $p_3$ is located in a new dimension (that means that the system is three dimensional) it might be possible continue with the argument and visit $p_3$.

In the next subsection we are going to generate an example for which it is clear that one can shadow a non-transverse heteroclinic chain.

\subsection{Examples with diffusion in a non-transverse situation - the triangular system}
\label{subsec:Examples}

The triangular system is
\begin{equation}
\dot{x}=F(x)  \label{eq:triang-ode}
\end{equation} with
\begin{equation}
\left\{\begin{array}{rcl}
F_1(x)&=&\lambda_1 x_1-\lambda_1x_1^2\\
F_i(x)&=&(\lambda_i-\mu_i)x_i-\lambda_ix_i^2+\mu_ix_ix_{i-1}\qquad \text{for }1<i\leq n
\end{array}\right.
\label{eq:F-triang-ode}
\end{equation} with $\lambda_i>0$ for $1\leq i\leq n$ and $\mu_i\in\mathbb{R}$ for $1<i<n$.

Note that we can integrate (\ref{eq:triang-ode}) recursively, since
$$\dot{x}_i=f_i(t)x_i-\lambda_i x_i^2\Rightarrow
x_i(t)=\frac{x_i(0)\exp\left(\int_0^t{f_i(s)ds}\right)}{1+\lambda_i \int_0^t{\exp\left(\int_0^s{f_i(r)dr}\right)ds}},$$
with $f_i(t)=\lambda_i-\mu_i+\mu_ix_{i-1}(t)$ and $x_{-1}(t)=1$.

For $i=0,1,\dots,n$ we denote by $p_i$ the following fixed point of  (\ref{eq:triang-ode})
\begin{equation}
  p_i=(p_{i,1},\dots,p_{i,n}), \quad p_{i,j}=1, j\leq i-1, \quad p_{i,j}=0, \quad j\geq i.
\end{equation}
Therefore we have
\begin{eqnarray*}
  p_0=(0,\dots,0), \ p_1=(1,0,0,\dots,0), \ p_2=(1,1,0,\dots,0),\dots, p_n=(1,1,\dots,1).
\end{eqnarray*}

Observe that the interval $C_i=(1,\dots,1,x_i,0,\dots,0)$, where the length of the initial sequence of $1$'s is equal $i-1$, for $x_i \in (0,1)$ is the heteroclinic orbit connecting $p_{i}$ and $p_{i+1}$.  The question is whether we can follow this heteroclinic chain, i.e., whether there exists an orbit
which starts in the close vicinity of $p_0$ then it visits consequtively the neighborhoods of $p_1$, $p_2$, \dots, $p_n$.

\begin{figure}[H]
\centering
	\includegraphics[width=0.6\textwidth]%
	{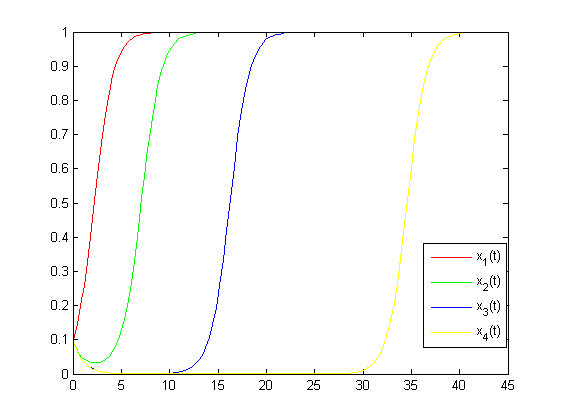}%
\caption{Solution of system \eqref{eq:triang-ode} for $\lambda=1$ and $\mu=2$. The initial condition is: $x_1(0)=x_2(0)=x_3(0)=x_4(0)=\frac{1}{10}$.}
\label{fig:mu-large}
\end{figure}

\begin{figure}[H]
\centering
	\includegraphics[width=0.6\textwidth]%
	{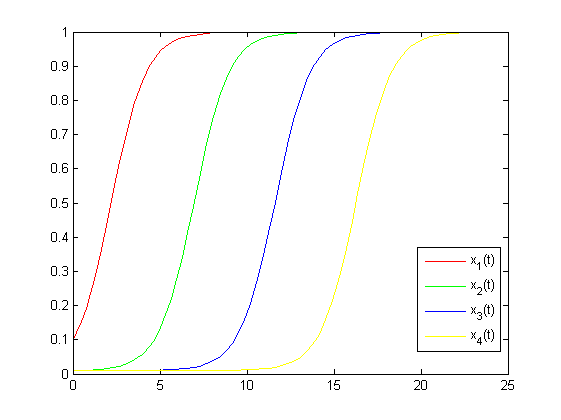}%
\caption{Solution of system \eqref{eq:triang-ode} for $\lambda=1$ and $\mu=1$. The initial condition is: $x_1(0)=\frac{1}{10},\quad x_2(0)=x_3(0)=x_4(0)=\frac{1}{100}$}
\label{fig:mu=lambda}
\end{figure}

\begin{figure}[H]
\centering
	\includegraphics[width=0.6\textwidth]%
	{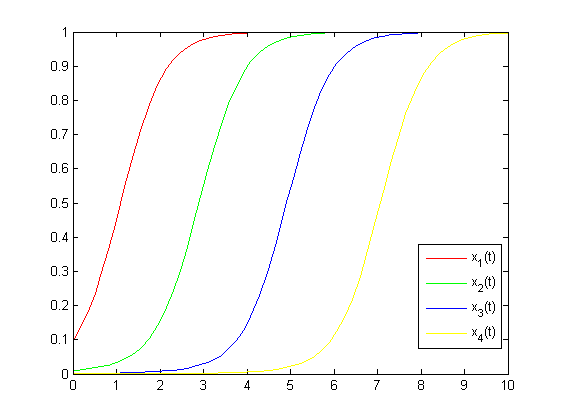}%
\caption{Solution of system \eqref{eq:triang-ode} for $\lambda=2$ and $\mu=1$. The initial condition is: $x_1(0)=10^{-1},\quad x_2(0)=10^{-2}\quad x_3(0)=10^{-3}\quad x_4(0)=10^{-4}.$}
\label{fig:lambda-large}
\end{figure}

We can check now the linear behavior around the equilibrium points computing the derivative of the vector field:

$$DF(p_0)=\left(\begin{array}{cccccc}
\lambda_1 &  &  &  & &\\
 & \lambda_2-\mu_2 & & & &\\
& &\ddots& & &\\
& & & \lambda_i-\mu_i & &\\
& & & & \ddots &\\
& & & & & \lambda_n-\mu_n
\end{array}\right)$$

$$DF(p_i)=\left(\begin{array}{cccccccc}
-\lambda_1 &  &  &  & &&&\\
\mu_2 & -\lambda_i & & & &&&\\
& \ddots &\ddots& & &&&\\
& & \mu_i & -\lambda_i & &&&\\
& & & &\lambda_{i+1} & &&\\
& & & & & \lambda_{i+2}-\mu_{i+2} &&\\
& & & & & & \ddots& \\
& & & & & & & \lambda_n-\mu_n
\end{array}\right)$$
so we have different possible choices for the parameters.
\begin{itemize}
\item If $\mu_i> \lambda_i$ for all $i=1,\dots, n$, each point $p_j$ has only one unstable direction defined by $e_{j+1}$, while the rest of the directions are stable.
\item If $\mu_i= \lambda_i$ for all $i=1,\dots, n$, each point $p_j$ has only one unstable direction, defined by $e_{j+1}$. All the ``past'' directions, defined by $\{e_1,\dots,e_j\}$ are stable while all the ``future'' directions, defined by $\{e_{j+2},\dots,e_n\}$, are linearly neutral.
\item If $\mu_i< \lambda_i$ for all $i=1,\dots, n$,  at each point, $p_j$, all the ``past'' directions, defined by $\{e_1,\dots,e_j\}$, are stable while all the ``future'' directions, defined by $\{e_{j+2},\dots,e_n)\}$, are linearly unstable.
\end{itemize}
Let us now display some numerical integration of the system for these three possibilities for $n=4$. For the sake of concreteness we are going to assume $\lambda=\lambda_i$ and $\mu=\mu_i$ for all $i=1,\dots,n$. In Figures~\ref{fig:mu-large}, \ref{fig:mu=lambda} and \ref{fig:lambda-large} we display the numerical evidence for the existence of diffusive orbits.

Of course, we could expect these numerical evidences just by looking at the equations. It is clear that \[\lim_{t\to\infty}{x_1(t)}=1.\]
By induction, assuming \[l_{i-1}=\lim_{t\to\infty}{x_{i-1}(t)}\]
and looking for equilibria of $x_i(t)$ for $t\to\infty$, it has to satisfy
\[\lim_{t\to\infty}{\dot{x}_{i}(t)}=0,\]
i.e. \[0=\lambda_il_i(1-l_i)-\mu_il_i(1-l_{i-1}).\]
So, we get $l_i=0$ or $l_i=1$ but since $\dot{x}_i(t)>0$ for $t$ large enough we can conclude $l_i=1$. This reasoning with some small effort can be turned into the rigorous proof. 

The main conclusion of this part is that we have obtained an easy example for which we can ensure the transition chain even if the intersection between the invariant manifolds is not transverse, regarding the high dimension and the disposition of the equilibrium points and the heteroclinic connections: each one in a new direction not used before. In addition, the system is integrable by quadratures which goes against the notion of the Arnold's diffusion.
However, we have detected the reason why the connection could be possible. The geometric mechanism relies on the fact that we are dealing with a high dimensional system and that each new connection is defined by a direction that \emph{has not been used before}.

\subsection{The conjecture}
\label{subsec:model-example}
We are now in condition to present a conjecture, which can be proved by our method discussed in the following sections under some additional
assumptions.

Let $n_i >0$ for $i=1,\dots,L$  and let $n_1+n_2 + \dots + n_L=n$, $d_i=n_1+\dots+n_i$.

For $i=1,\dots,L$ the subspaces $V_i$ which are spanned by
$\{e_{k_{i-1}+j}\}_{j=1,\dots,n_i}$. In this notation
$\mathbb{R}^n = \bigoplus_{i=1}^L V_i$.

We will use the following notation: for $l \in \mathbb{N}$,
$0^l,1^l$ will denote the sequences of length $l$ consisting of
$l$ $0$'s or $l$ $1$'s, respectively.

Assume that we have a diffeomorphism  $f:\mathbb{R}^n \times
\mathbb{R}^{w_u} \times \mathbb{R}^{w_s} \to \mathbb{R}^n\times \mathbb{R}^{w_u} \times  \mathbb{R}^{w_s}$, $w=w_u + w_s$ with the
following properties:
\begin{itemize}
\item there exists a sequence of fixed points
   \begin{eqnarray*}
   p_0&=&(0^n,0^w), \\
   p_1&=&(1,0^{n_1-1},0,\dots,0^w), \\
    p_2&=&(1,0^{n_1-1},1,0^{n_2-1},0\dots,0^w),\\
     &\vdots& \\
    p_k&=&(1,0^{n_1-1},1,0^{n_2-1},\dots,1,0^{n_k-1},0^{n_{k+1}+\dots + n_L},0^w), \\
    &\vdots& \\
    p_{L}&=&(1,0^{n_1-1},1,0^{n_2-1},\dots,1,0^{n_L-1},0^w)
  \end{eqnarray*}
\item for any $i=1,\dots,L$,  the interval connecting $p_{i-1}$ and $p_{i}$ denoted by $C_i$
\begin{equation*}
C_i=\{ z = (1,0^{n_1-1},1,0^{n_2-1},\dots,1,0^{n_{i-1}-1},t,0^{-1+n_i+\dots + n_L},0^w ), \quad t \in [0,1]   \}
\end{equation*}
  is invariant under $f$ and for any $z \in C_i$
  \begin{equation}
    \lim_{k \to \infty} f^k(z)=p_i, \quad   \lim_{k \to -\infty} f^k(z)=p_{i-1}
  \end{equation}
\item at $p_i$ the subspace  $V_i$ defined the  exit directions  and these directions are 'dominating' for the scattering
   when passing by $p_i$.    This statement  is vague, because the scenario we are going to present most likely can be realized
   under various sets of assumptions.
\end{itemize}

Only the first $n$-directions really count, the others will be treated as the entry directions
(in the sense of covering relations, see Section~\ref{sec:covrel}).

Our  result is is
\begin{con}
\label{con:long-transition} Under the above assumptions for any
$\epsilon >0$ there exists a point $z $ and  a
sequence of integers $k_1 < k_2 < \dots < k_L$, such that
\begin{eqnarray*}
  \|z - p_0\| &<& \epsilon , \\
  \|f^{k_i}(z) - p_{i}\| &<& \epsilon, \quad i=1,\dots,L
\end{eqnarray*}
\end{con}

Our idea of the proof of this result requires a construction of covering relations (see
(\ref{eq:cov-M-tildeM}), and  then
Conjecture~\ref{con:long-transition} follows directly from
Theorem~\ref{thm:top-chn-long-transition}.

We will show how a construction of suitable h-sets and coverings can be done for a linear model in Section~\ref{sec:linModel-proof} and for the simplified version of the toy model from \cite{CK} in Section~\ref{sec:diffToymodel}.

\section{The geometric idea of dropping dimensions }
\label{sec:geomIdea}

In this section we will explain  our idea of dropping dimensions along a nontransversal heteroclinic chain.

Let us start with a simplified version of the example from Section~\ref{subsec:model-example}.
 Let $f:\mathbb{R}^n\rightarrow\mathbb{R}^n$ be a diffeomorphism with the following properties:
\begin{enumerate}
 \item The points $p_i=(1,\stackrel{i}{\dots},1,0,\stackrel{n-i}{\dots},0)$ are fixed under $f$ for $i=0\dots n$.
 \item The segments $C_i$ that connect the points $p_{i-1}$ and $p_{i}$, $$C_i=\{(1,\stackrel{i-1}{\dots},1,t,0,\stackrel{n-i}{\dots},0),0\leq t\leq1\}$$ for $1\leq i\leq n$ are invariant under $f$ and, for all $x\in C_i$: $$\lim_{k\rightarrow\infty}{f^k(x)}=p_i \qquad\lim_{k\rightarrow-\infty}{f^k(x)}=p_{i-1}.$$
\item  At each point $p_i$ the $i$-th direction is stable and the $(i+1)$-th is unstable. This means:
$$Df(p_i)e_i=\mu_ie_i,\quad |\mu_i|<1$$ $$Df(p_i)e_{i+1}=\lambda_ie_{i+1},\quad |\lambda_i|>1$$
\item The past directions, defined by $e_1,\dots,e_{i-1}$, are contracting directions around the fixed point $p_i$ but with a rate greater than $\mu_i$. The future directions, defined by $e_{i+2},\dots,e_n$, are expanding directions around the fixed point $p_i$ but with a rate lower than $\lambda_i$.
\end{enumerate}

\begin{pro}
\label{conj:very-simple-model}
 Under the previous assumptions, for all $\epsilon>0$ there exists a point $x$ and a sequence of integers $0=k_0<k_1<\dots<k_n$ such that: $$||f^{k_i}(x)-p_i||<\epsilon \quad i=0,\dots,n.$$
\end{pro}
\begin{rem}
 Notice that we connect $n+1$ points in a $n$ dimensional space. We cannot guarantee that the result is valid for more points.
\end{rem}

\begin{rem}
 Observe that we are assuming that there are only two dominant coordinates around each fixed point. That means that this could not be applied to the Toy Model System in NLS \cite{CK,GK}, where we have four dominant directions. In Section~\ref{sec:diffToymodel} we treat such system.
\end{rem}

Proposition~\ref{conj:very-simple-model} is proved in Section~\ref{sec:linModel-proof} for the linear model.

\subsection{The idea of the proof: dropping dimensions}
\label{subsubsec:sketch}

Here we  sketch the basic idea of the proof  of Proposition \ref{conj:very-simple-model} with some pictures. We are going to consider only a two dimensional map. So, consider $f:\mathbb{R}^2\rightarrow\mathbb{R}^2$ with three fixed points:
\[p_0=(0,0),\quad p_1=(1,0), \quad p_2=(1,1),\]
with invariant segments $C_1$ and $C_2$ defined as
\[C_1=\{(x_1,x_2)\,:\,0\leq x_1\leq 1,\,x_2=0\}\quad
C_2=\{(x_1,x_2)\,:\, x_1=0,\,0\leq x_2\leq 1\}\]
\begin{figure}[H]
\centering
	\includegraphics[width=0.7\linewidth]%
	{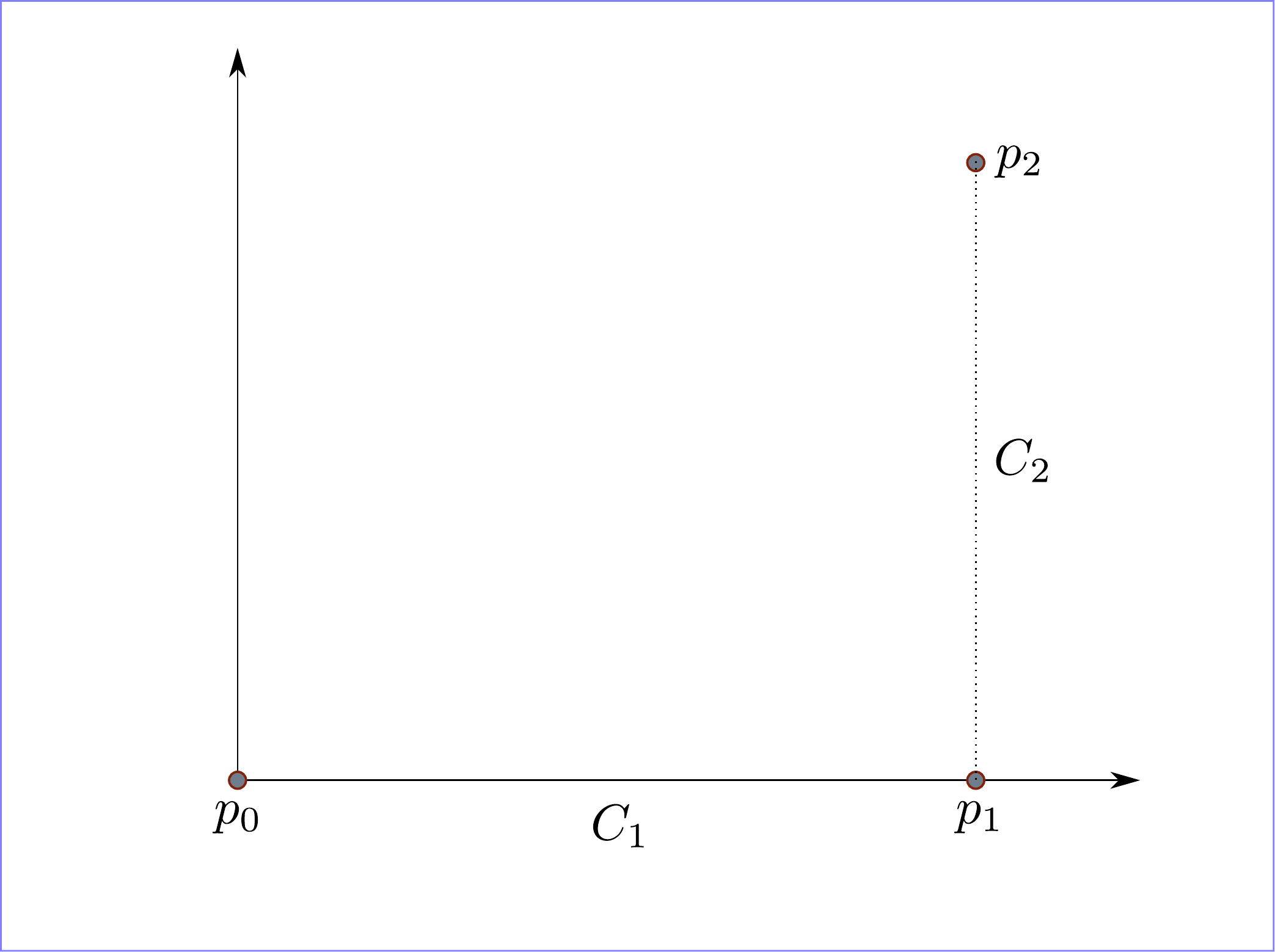}%
\end{figure}
Assume also that the derivatives of the map around the fixed points have the following structure:
\begin{equation}
Df(p_0)=\left(\begin{array}{cc}\lambda_{0,1} & 0 \\ 0 & \lambda_{0,2}\end{array}\right),\
Df(p_1)=\left(\begin{array}{cc}\mu_{1,1} & 0 \\ 0 & \lambda_{1,2}\end{array}\right),\
Df(p_2)=\left(\begin{array}{cc}\mu_{2,1} & 0 \\ 0 & \mu_{2,2}\end{array}\right),
\label{linstab}\end{equation}
where $\lambda_{0,1},\lambda_{0,2},\lambda_{12}>1$ and $0<\mu_{1,1},\mu_{2,1},\mu_{2,2}<1$.

We start by considering a domain (ball) $D_0$ of full dimension centered around $p_0$.
\begin{figure}[H]
\centering
	\includegraphics[width=0.7\linewidth]%
	{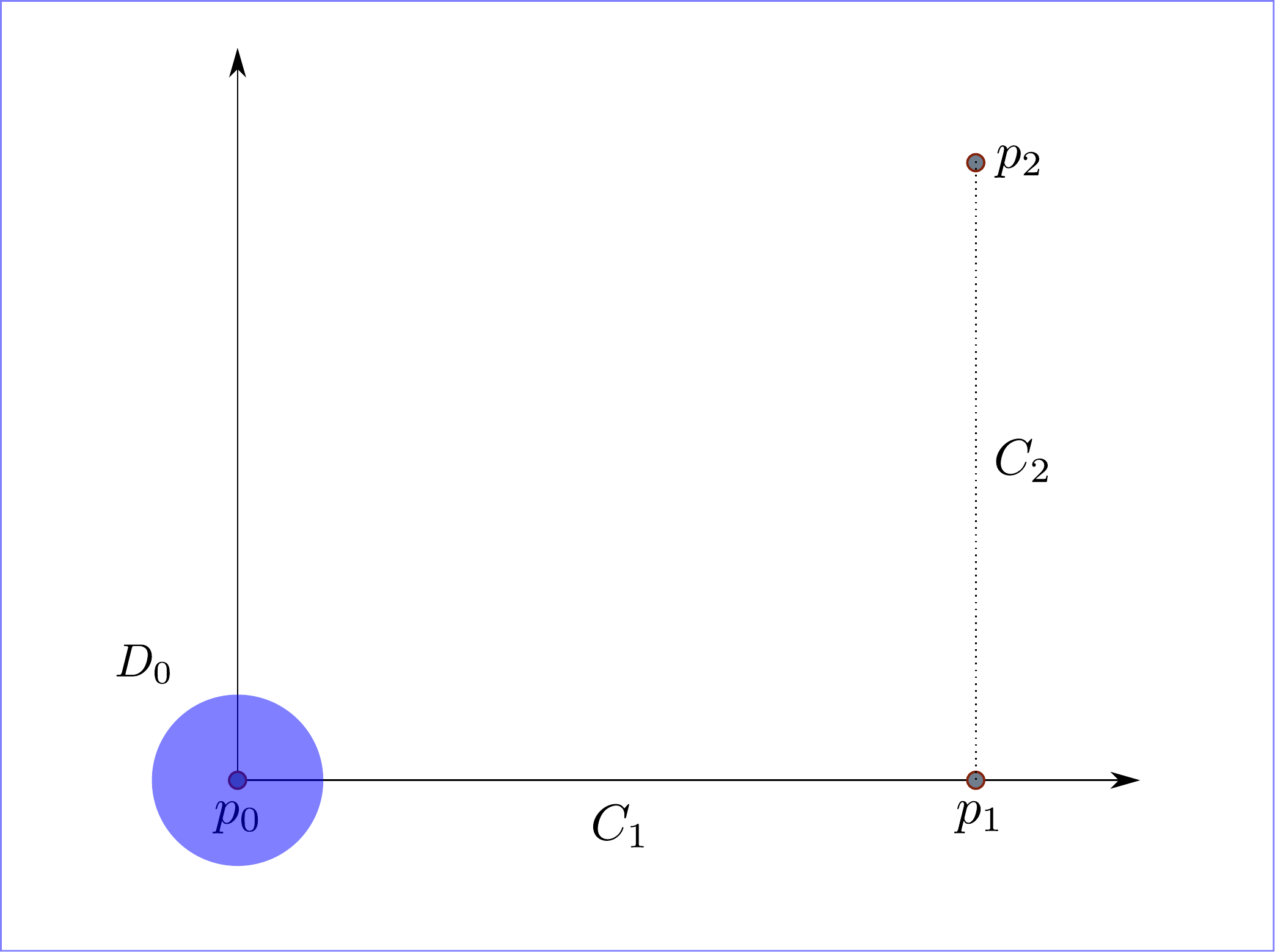}%
\end{figure}
Given the linear stability from \eqref{linstab}, we can assume that after one iteration of the map, our initial ball $D_0$ will be expanded in both directions:
\begin{figure}[H]
\centering
	\includegraphics[width=0.7\linewidth]%
	{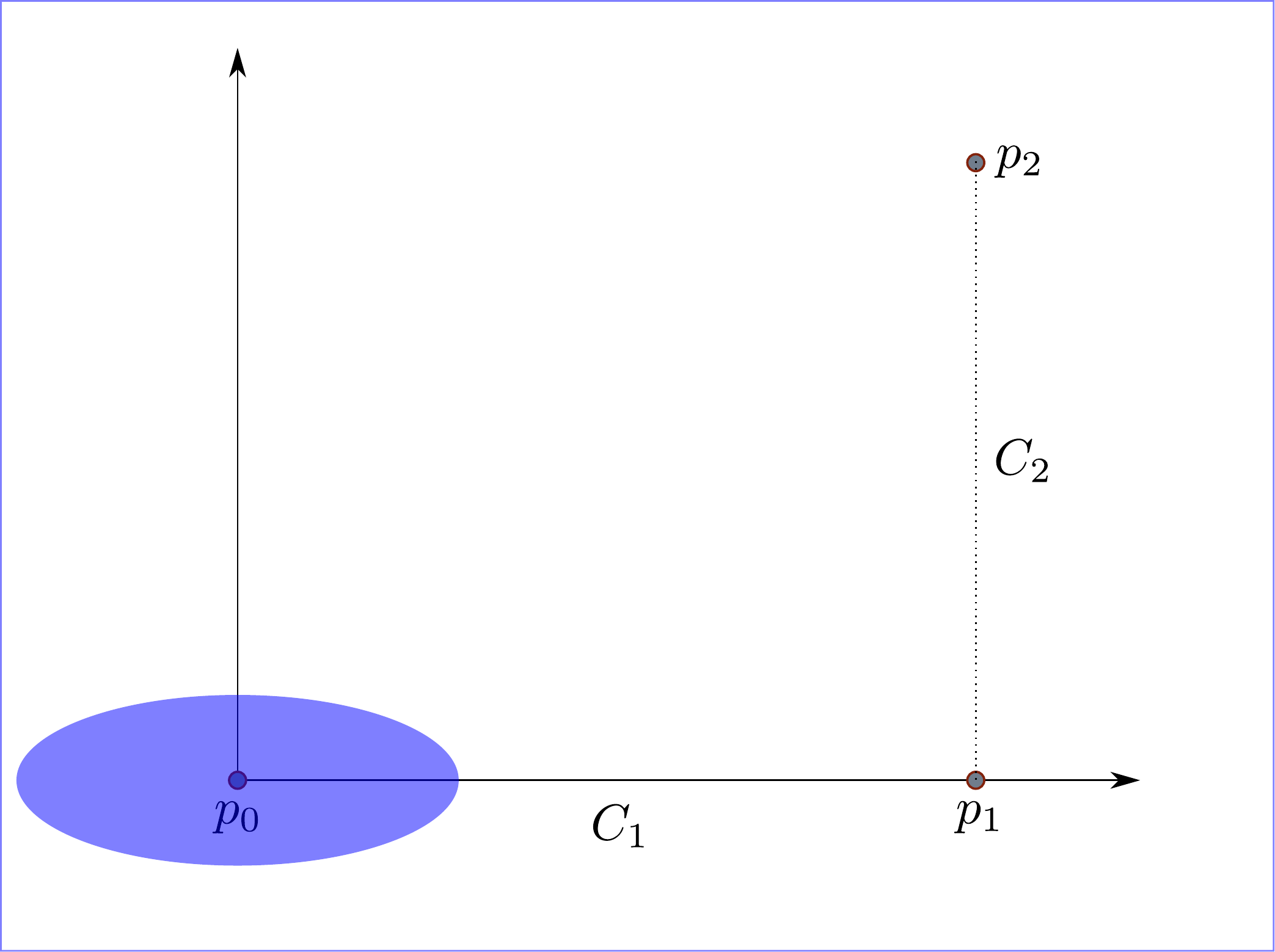}%
\end{figure}
It is now time to make a decision: from all the possible directions, we are only interested in the one defined by the outgoing heteroclinic connection,
that is, the segment $C_1$. Then we consider a section $S_0=\{x_1=\sigma\}$ where $\sigma$ is some small parameter:
\begin{figure}[H]
\centering
	\includegraphics[width=0.7\linewidth]%
	{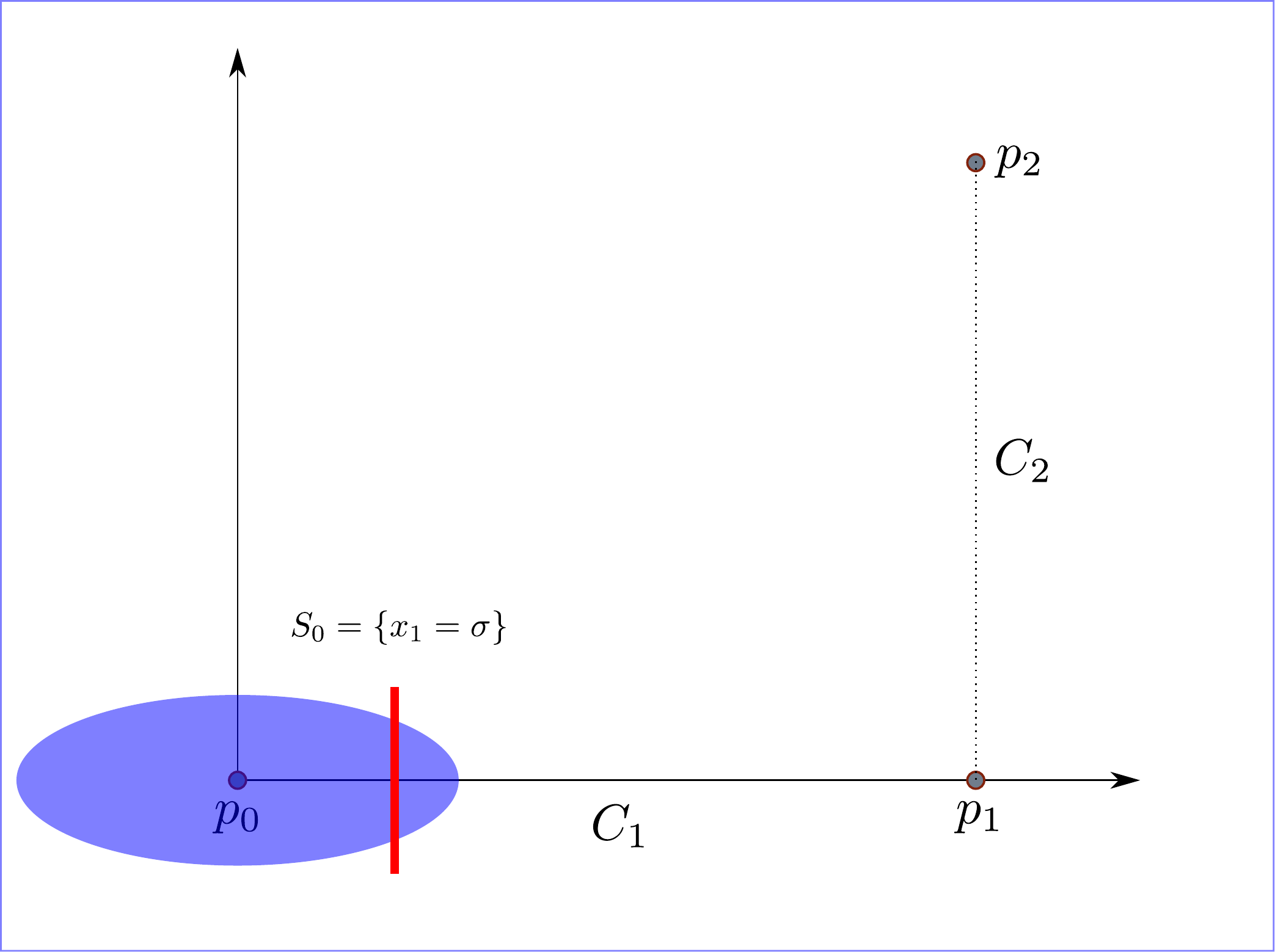}%
\end{figure}
Since we are only interested in the points of our ball close to the heteroclinic connection, we intersect the domain with the section $S_0$. We say that we have dropped the $x_1$ direction. We will not use this direction in future steps. Our domain, $\bar{D}_0$, has one dimension less than $D_0$, that is, it has dimension one.
\begin{figure}[H]
\centering
	\includegraphics[width=0.7\linewidth]%
	{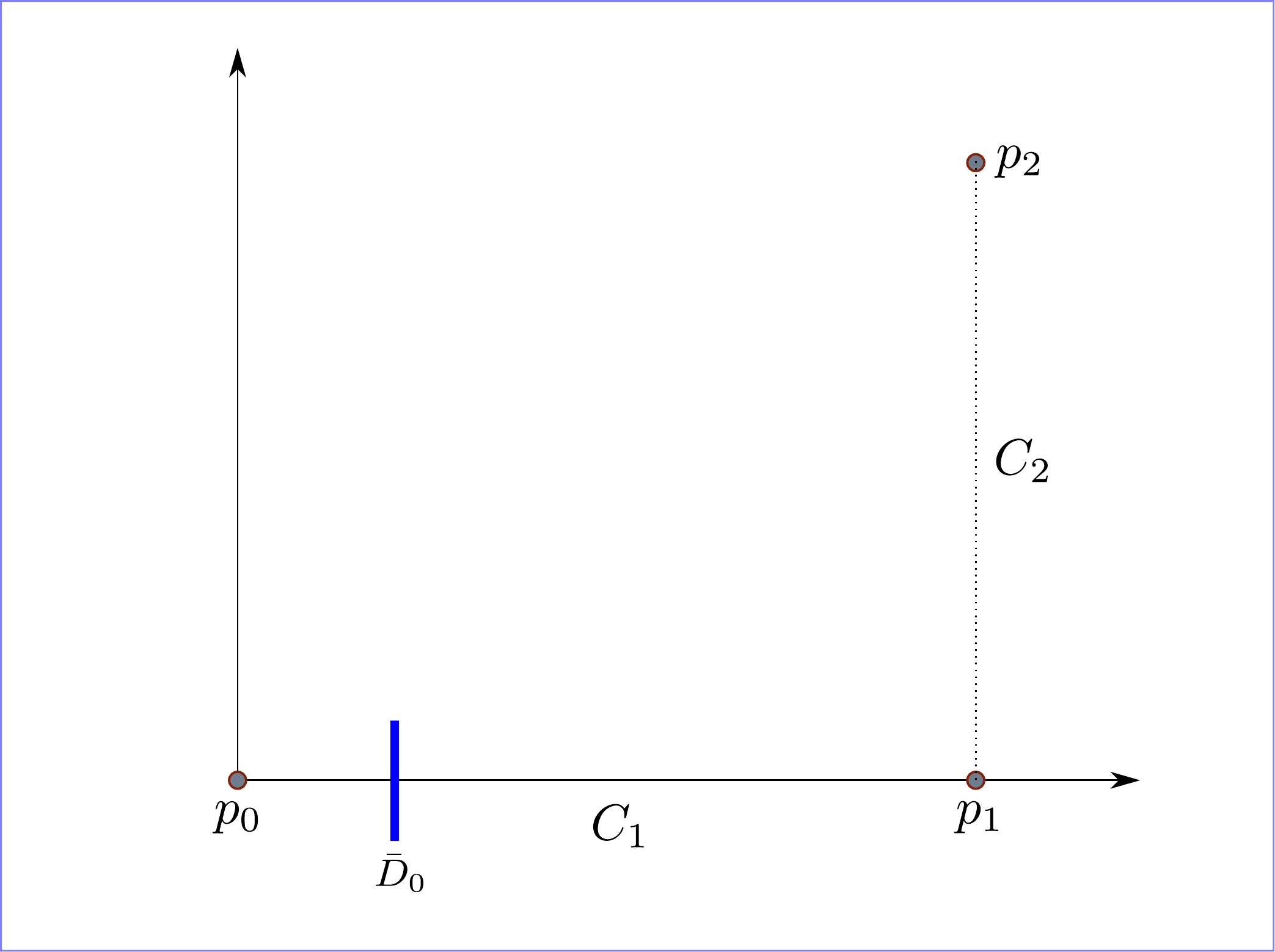}%
\end{figure}
Then we continue with this domain. After several iterations of the map since the domain is close to the heteroclinic connection, we can ensure that $\bar{D}_0$ will approach $p_1$ and we obtain a domain $D_1$:
\begin{figure}[H]
\centering
	\includegraphics[width=0.7\linewidth]%
	{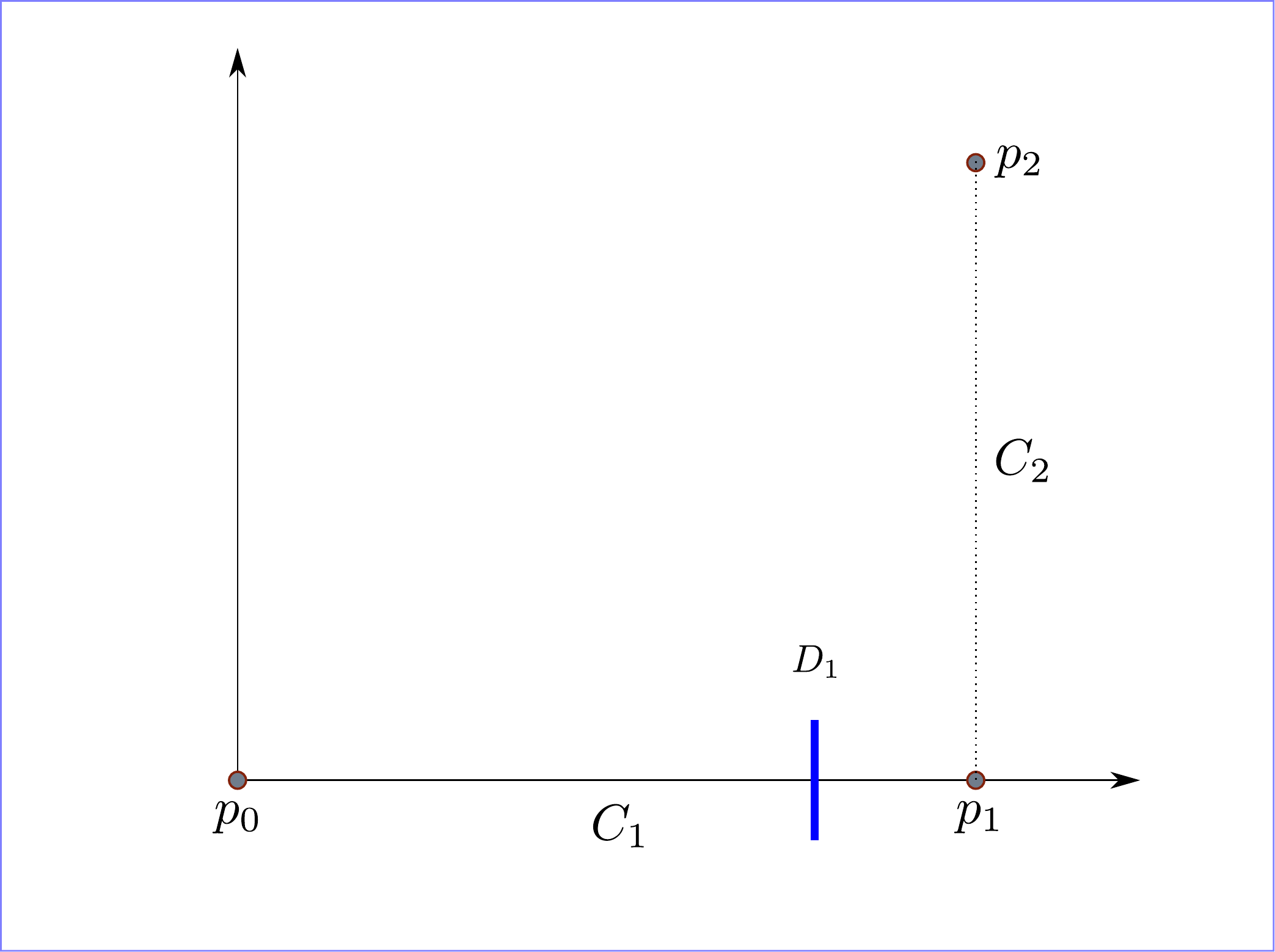}%
\end{figure}
We can use, again, the linear prediction of the map, \eqref{linstab}, to be sure that our domain is expanded in the $x_2$ direction.
\begin{figure}[H]
\centering
	\includegraphics[width=0.7\linewidth]%
	{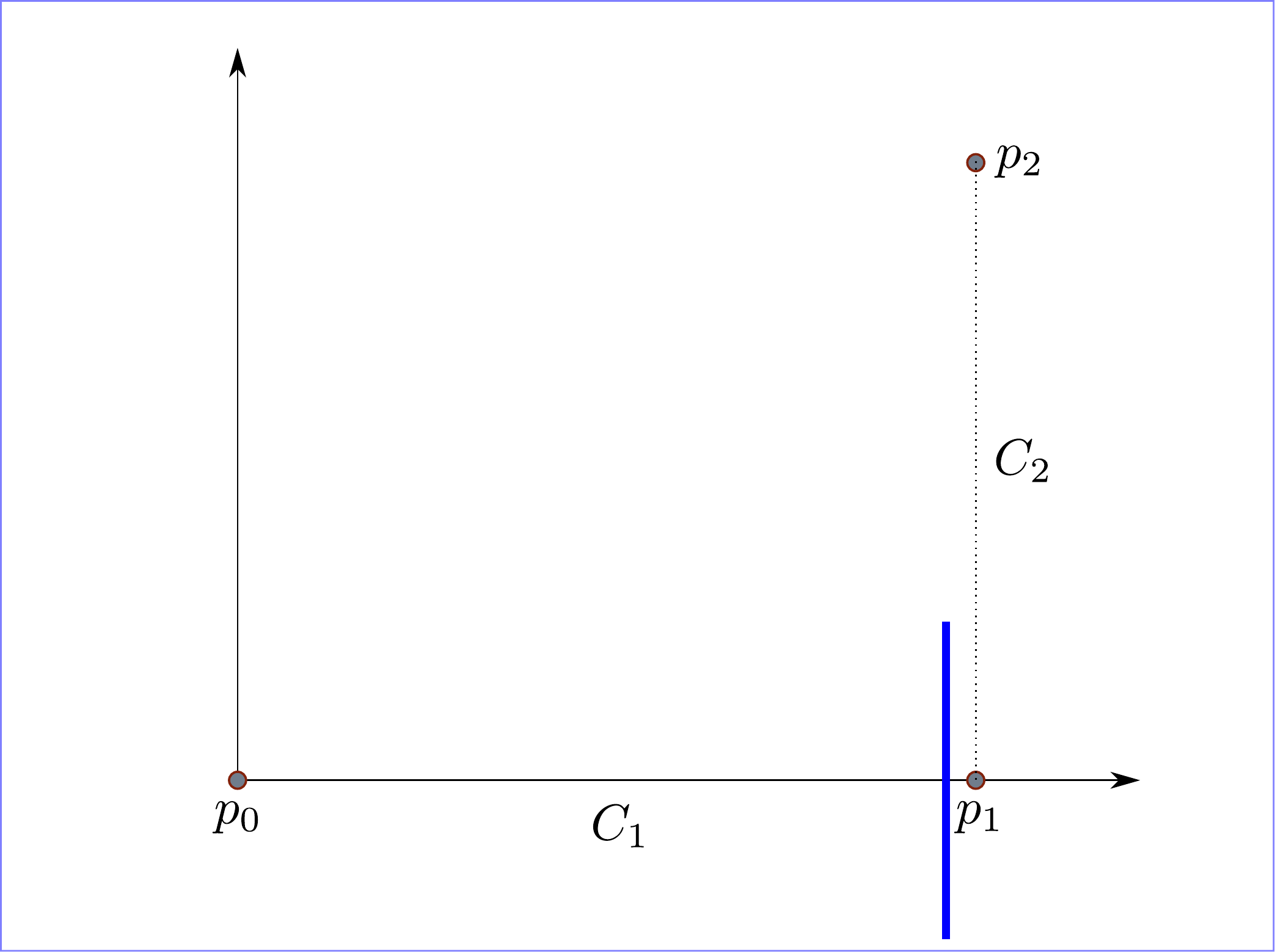}%
\end{figure}
We use now the same argument. From all the possible directions that $f(D_1)$ covers, we want to escape through the one defined by the heteroclinic connection
to $p_2$. So we put a section defined in the same spirit as before: $S_1=\{x_2=\sigma\}$.
\begin{figure}[H]
\centering
	\includegraphics[width=0.7\linewidth]%
	{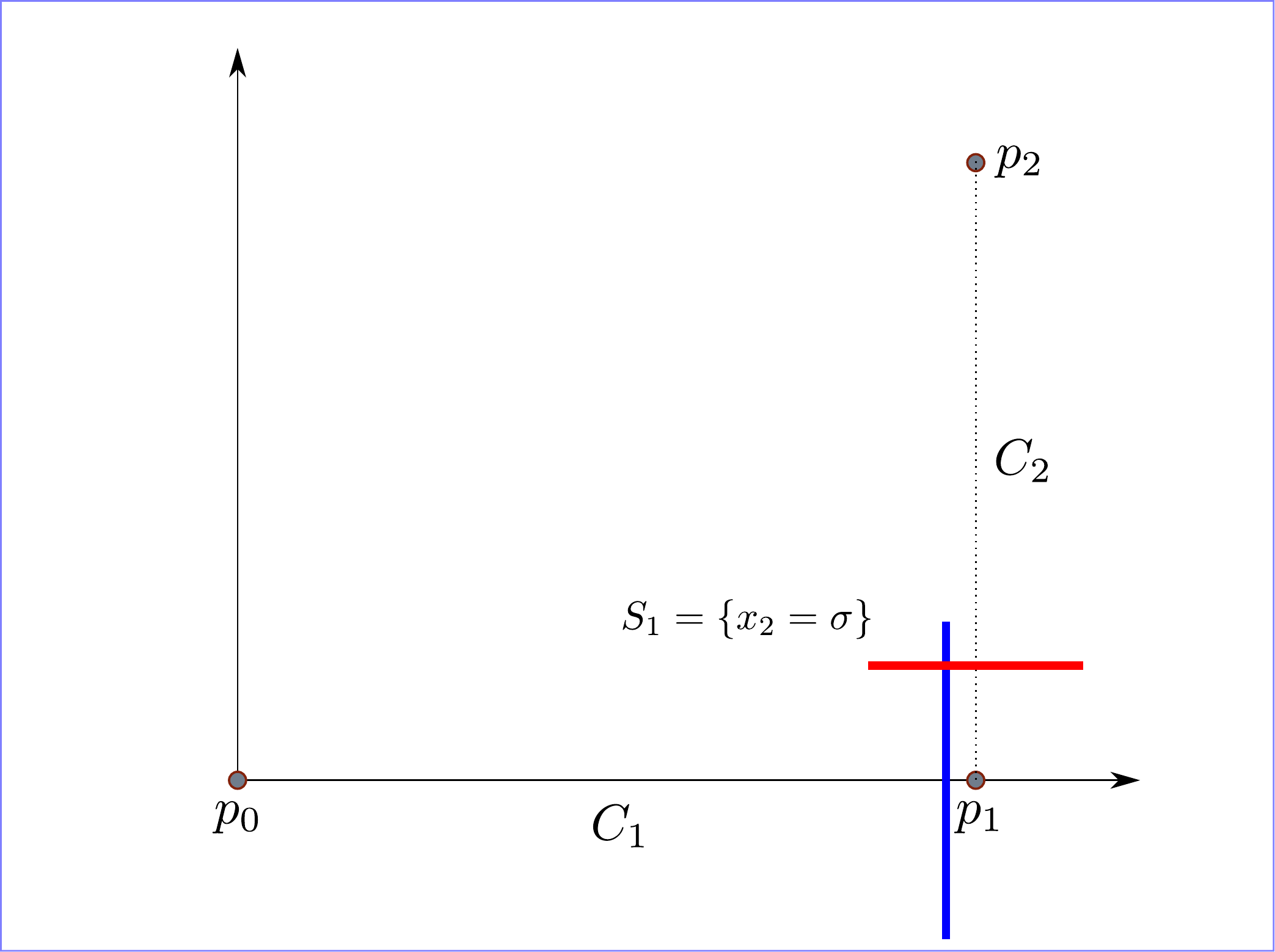}%
\end{figure}
We restrict now our domain in its intersection with the section $S_1$. The resulting domain $\bar{D}_1$ will have, then, one dimension less than $D_1$, which means that it will have dimension zero.
\begin{figure}[H]
\centering
	\includegraphics[width=0.7\linewidth]%
	{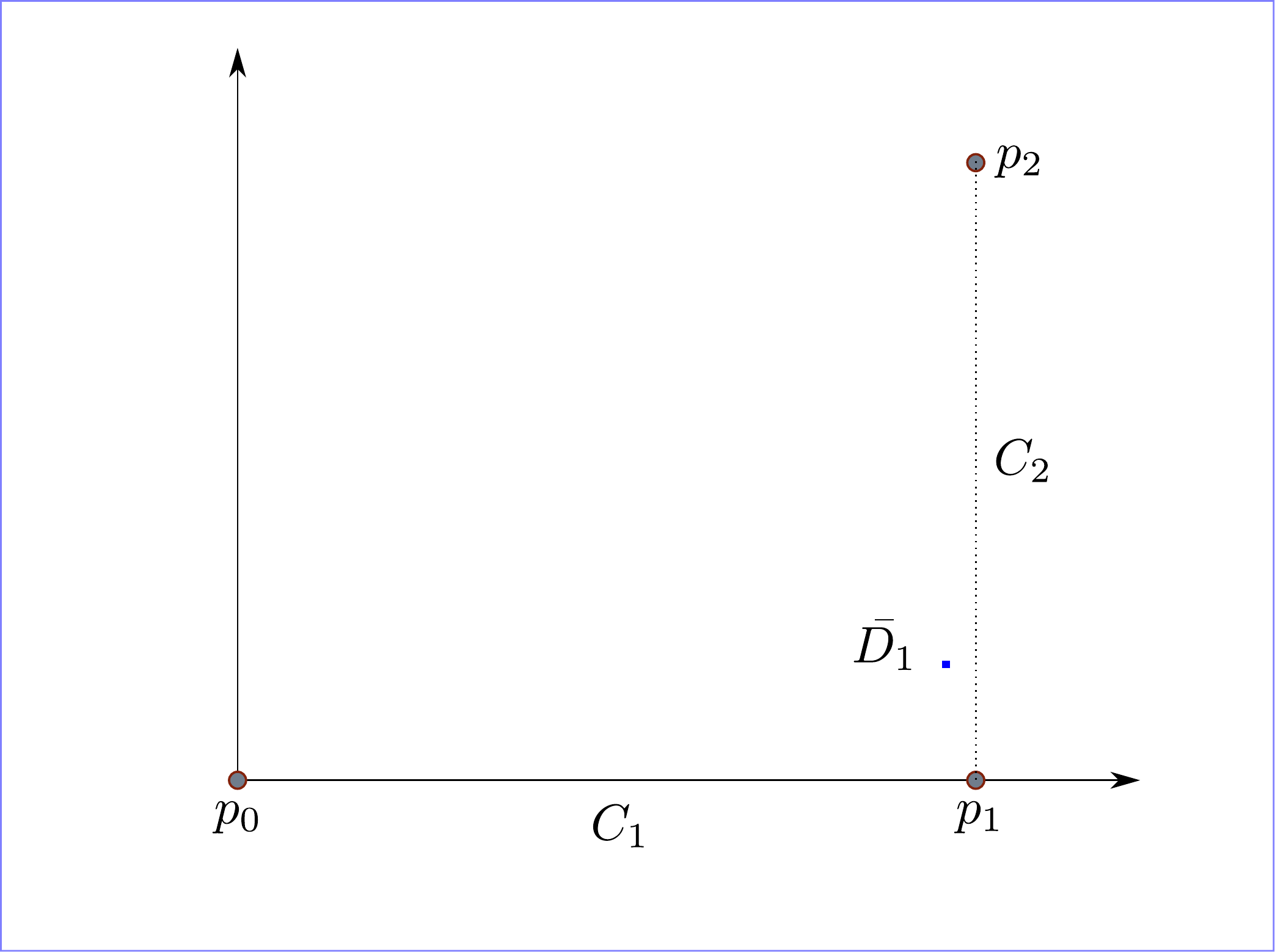}%
\end{figure}
We have no more dimensions to drop, since our initial domain becomes a single point. This point is close to the heteroclinic defined in $C_2$, so that we are sure that after some iterates, it will approach the final fixed point $p_2$:
\begin{figure}[H]
\centering
	\includegraphics[width=0.7\linewidth]%
	{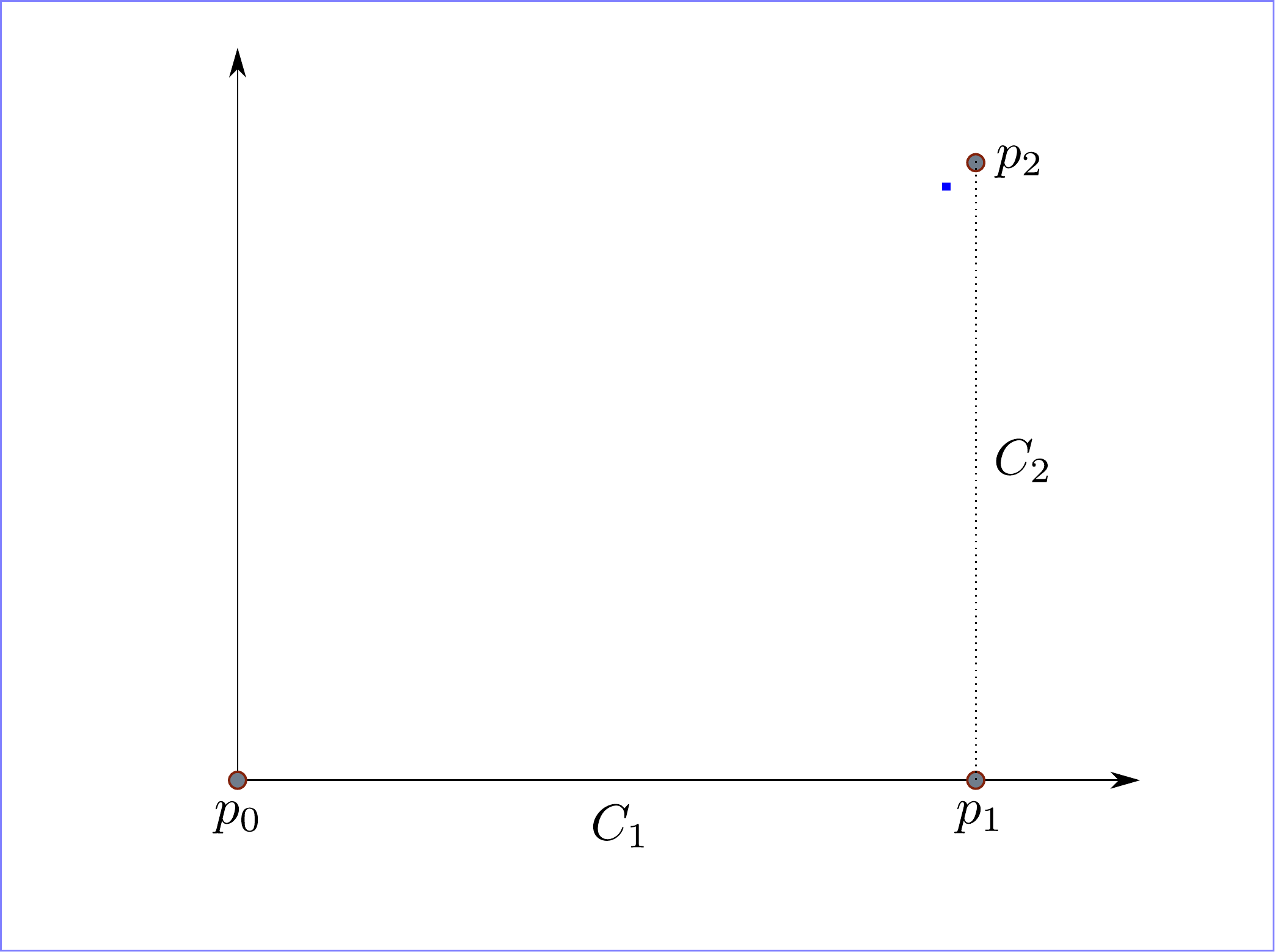}%
\end{figure}


\section{h-sets, covering relations}
\label{sec:covrel}

The goal of this section is to recall from \cite{ZGi} the notions of h-sets and covering relations, and to state the theorem
about the existence of point realizing the chain of covering relations. This will be the main technical tool in proving the existence of the orbits
shadowing the  heteroclinic chain in the next sections.

\subsection{h-sets and covering relations}
\label{subsec:covrel}

\begin{definition} \cite[Definition 1]{ZGi}
\label{def:covrel} An $h$-set $N$ is a quadruple \\
$(|N|,u(N),s(N),c_N)$ such that
\begin{itemize}
 \item $|N|$ is a compact subset of ${\mathbb R}^n$
 \item $u(N),s(N) \in \{0,1,2,\dots,n\}$ are such that $u(N)+s(N)=n$
 \item $c_N:{\mathbb R}^n \to
   {\mathbb R}^n={\mathbb R}^{u(N)} \times {\mathbb R}^{s(N)}$ is a
   homeomorphism such that
      \begin{displaymath}
        c_N(|N|)=\overline{B_{u(N)}} \times
        \overline{B_{s(N)}}.
      \end{displaymath}
\end{itemize}
We set
\begin{eqnarray*}
   \dim(N) &:=& n,\\
   N_c&:=&\overline{B_{u(N)}} \times \overline{B_{s(N)}}, \\
   N_c^-&:=&\partial B_{u(N)} \times \overline{B_{s(N)}}, \\
   N_c^+&:=&\overline{B_{u(N)}} \times \partial B_{s(N)}, \\
   N^-&:=&c_N^{-1}(N_c^-) , \quad N^+=c_N^{-1}(N_c^+).
\end{eqnarray*}
\end{definition}

Hence a $h$-set $N$ is a product of two closed balls in some
coordinate system. The numbers $u(N)$ and $s(N)$ are called the
exit and entry dimensions, respectively.
The subscript $c$ refers to the new coordinates given by
the homeomorphism $c_N$. Observe that if $u(N)=0$, then
$N^-=\emptyset$ and if $s(N)=0$, then $N^+=\emptyset$. In the
sequel to make the notation less cumbersome we will often drop the
bars in the symbol $|N|$ and we will use $N$ to denote both the
h-sets and its support.

We will call $N^-$ \emph{the exit set of N} and $N^+$
\emph{the entry set of $N$}. These names are motivated by the
Conley index theory \cite{C,MM} and the role that these sets will play in the
context of covering relations.

\begin{definition}\cite[Definition 6]{ZGi}
\label{def:covw} Assume that $N,M$ are $h$-sets, such that
$u(N)=u(M)=u$ and $s(N)=s(M)=s$. Let $f:N \to {\mathbb R}^n$ be a
continuous map. Let $f_c= c_M \circ f \circ c_N^{-1}: N_c \to
{\mathbb R}^u \times {\mathbb R}^s$. Let $w$ be a nonzero integer.
We say that
\begin{displaymath}
  N\cover{f,w} M
\end{displaymath}
($N$ $f$-covers $M$ with degree $w$) iff the following conditions
are satisfied
\begin{description}
\item[1.] There exists a continuous homotopy $h:[0,1]\times N_c \to {\mathbb R}^u \times {\mathbb R}^s$,
   such that the following conditions hold true
   \begin{eqnarray}
      h_0&=&f_c,  \label{eq:hc1} \\
      h([0,1],N_c^-) \cap M_c &=& \emptyset ,  \label{eq:hc2} \\
      h([0,1],N_c) \cap M_c^+ &=& \emptyset .\label{eq:hc3}
   \end{eqnarray}
\item[2.] If $u >0$, then there exists a  map $A:{\mathbb R}^u \to {\mathbb
R}^u$ such that
   \begin{eqnarray}
    h_1(p,q)&=&(A(p),0), \mbox{ for $p \in \overline{B_u}(0,1)$ and $q \in
    \overline{B_s}(0,1)$,}\label{eq:hc4}\\
      A(\partial B_u(0,1)) &\subset & {\mathbb R}^u \setminus
      \overline{B_u}(0,1).  \label{eq:mapaway}
   \end{eqnarray}
  Moreover, we require that
\begin{equation}
  \deg(A,\overline {B_u}(0,1),0)=w, \label{eq:deg-A}
\end{equation}
\end{description}

We will call condition \eqref{eq:hc2} \emph{the exit condition}
and condition \eqref{eq:hc3} \emph{the entry
condition}.

\end{definition}
Note that in the case $u=0$, if $N \cover{f,w} M$, then $f(N)
\subset \inter M$ and  $w=1$.

In fact in the above definition $s(N)$ and $s(M)$ can be
different, see \cite[Def. 2.2]{W2}.


\begin{rem}If the map $A$  in condition 2 of Def.~\ref{def:covw} is a linear
map, then  condition (\ref{eq:mapaway}) implies that
\begin{displaymath}
  \deg(A,\overline {B_u}(0,1),0)=\pm 1.
\end{displaymath}
Hence condition (\ref{eq:deg-A}) is fulfilled with $w =\pm 1$.

  In fact, this is the most common situation in the applications of
  covering relations.
\end{rem}

Most of the time we will not be interested in the value of $w$ in the
symbol $N \cover{f,w} M$ and we will often drop it and write  $N
\cover{f} M$, instead. Sometimes we may even drop the symbol $f$, if known from the context,
and write $N \cover{} M$.

To handle the inverse maps in the context of covering relations the following two definitions are useful.

\begin{definition} \cite[Definition 3]{ZGi}
\label{def:hsetT}
 Let $N$ be a $h$-set.
We define a $h$-set $N^T$ as follows
\begin{itemize}
 \item $|N^T|=|N|$
 \item $u(N^T)=s(N)$,  $s(N^T)=u(N)$
 \item We
 define a homeomorphism $c_{N^T}:{\mathbb R}^n \to   {\mathbb R}^n={\mathbb R}^{u(N^T)} \times {\mathbb
R}^{s(N^T)}$,
 by
      \begin{displaymath}
        c_{N^T}(x)= j(c_{N}(x)) ,
      \end{displaymath}
      where $j: {\mathbb R}^{u(N)} \times {\mathbb R}^{s(N)} \to {\mathbb R}^{s(N)} \times {\mathbb R}^{u(N)}$
      is given by $j(p,q)=(q,p)$.
\end{itemize}
\qed
\end{definition}
Observe that $N^{T,+}=N^-$ and $N^{T,-}=N^+$. This operation is
useful in the context of inverse maps.

\begin{definition}\cite[Definition 7]{ZGi}
\label{def:backcov} Assume $N,M$ are $h$-sets, such that
$u(N)=u(M)=u$ and $s(N)=s(M)=s$.  Let $g: {\mathbb R}^n \supset
\Omega   \to {\mathbb R}^n$. Assume that $g^{-1}:|M| \to {\mathbb
R}^n$ is well defined and continuous. We say that $N \invcover{g}
M$ ($N$ $g$-backcovers $M$ ) iff $M^T \cover{g^{-1}} N^T$.
\end{definition}

\subsection{Main theorem about chains of covering relations}

\begin{theorem}[Thm. 9]\cite{ZGi}
\label{thm:cov}
 Assume $N_i$, $i=0,\dots,k$, $N_k=N_0$ are
$h$-sets and for each $i=1,\dots,k$ we have either
\begin{equation}
  N_{i-1} \cover{f_i,w_i} N_{i} \label{eq:dirgcov}
\end{equation}
or
\begin{equation}
  N_i \subset \dom(f_i^{-1}) \quad \mbox{and} \quad  N_{i-1} \invcover{f_i,w_i} N_{i}.
  \label{eq:invgcov}
\end{equation}

 Then there exists a point $x \in \inter N_0$, such that
\begin{eqnarray}
   f_i \circ f_{i-1}\circ \cdots \circ f_1(x) &\in& \inter N_i, \quad i=1,\dots,k \\
  f_k \circ f_{k-1}\circ \cdots \circ f_1(x) &=& x
\end{eqnarray}
\end{theorem}
The reader is referred to \cite{ZGi} for a proof.  The basic idea of the proof of this theorem is the homotopy and the local Brouwer degree.

The following corollary is an immediate consequence of Theorem
\ref{thm:cov}.
\begin{col}
\label{col:corinvgencov} Let $N_i$, $i\in\mathbb{Z}_+$ be h-sets.
Assume that  for each $i\in\mathbb{Z}_+$ we have either
\begin{equation}
  N_{i-1} \cover{f_i,w_i} N_{i}
\end{equation}
or
\begin{equation}
  N_{i} \subset \dom(f_i^{-1})   \quad  \mbox{and}
    \quad N_{i-1} \invcover{f_i,w_i} N_{i}.
\end{equation}

Then there exists a point $x \in \inter N_0$, such that
\begin{eqnarray}
   f_i \circ f_{i-1}\circ \cdots \circ f_1(x) &\in& \inter N_i, \quad
   i\in \mathbb{Z}_+.
\end{eqnarray}
Moreover, if $N_{i+k}=N_i$ for some $k>0$ and all $i$, then the
point $x$ can be chosen so that
\begin{eqnarray}
   f_{k}\circ f_{k-1}\circ \cdots \circ f_1(x)=x.
\end{eqnarray}
\end{col}

\subsection{Natural structure of a h-set}
Observe that all the conditions appearing in the definition of the
covering relation are expressed in `internal' coordinates $c_{N}$
and $c_{M}$. Also the homotopy is defined in terms of these
coordinates.  Sometimes this makes statements and notation seem
a bit cumbersome. With this in mind we introduce the notion of a
natural structure on a h-set.

\begin{definition}
\label{def:hset-nat-stru} We will say that $N = \{(x_0,y_0)\} +
\overline{B}_u(0,r_1) \times \overline{B}_s(0,r_1) \subset
\mathbb{R}^u \times \mathbb{R}^s$ is an \emph{$h$-set} with a
natural structure if:
\newline
 $u(N)=u$,
$s(N)=s$, $c_{N}(x,y)= \left(\frac{x-x_0}{r_1},\frac{y-y_0}{r_2}
\right)$.

\end{definition}

\subsection{The operation of dropping exit dimensions}

\begin{definition}
\label{def:dropping-exit-dim}
Assume that we have a decomposition $\mathbb{R}^n=\mathbb{R}^{u_1}  \oplus \mathbb{R}^t \oplus   \mathbb{R}^{s_1}$ and the norm for $(x_1,x_2,x_3) \in \mathbb{R}^{u_1}  \oplus \mathbb{R}^t \oplus  \mathbb{R}^{s_1}$ is $\|(x_1,x_2,x_3)\|=\max (\|x_1\|,\|x_2\|,\|x_3\|)$.

Assume that $N$ is an h-set, with $u(N)=u_1+t$ and $s(N)=s_1$.  In view of the norm on $\mathbb{R}^n$ we have
\begin{equation}
  c_{|N|}=\left(\overline{B}_{u_1} \oplus \overline{B}_{t} \right) \oplus \overline{B}_{s_1}
\end{equation}
where the parentheses enclose the exit directions.

Let us denote by $V$ the subspace $\{0\} \times \mathbb{R}^t \times \{0\}$. We define a new h-set $R_V (N)$  by setting
\begin{itemize}
  \item $|R_V (N)|=|N|$
  \item $u(R_V (N))=u_1$, $s(R_V (N))=s_1+t$
  \item $c_{R_V (N)}=c_N$
\end{itemize}
\end{definition}
Roughly speaking,  $R_V(N)$ is obtained from $N$ by relabeling some exit coordinates in $N$ as the entry directions.

\section{The mechanism of dropping dimensions---the main topological theorem }
\label{sec:topThm}

\subsection{h-sets $M_i$ and $\widetilde{M}_i$}
\label{subsec:NM-def}
\outcomment{The basic ingredient of the method is a
careful treatment of the effect of 'scattering' of orbits coming
along one of the entry directions and then being scattered along
one selected exit direction. We will see, that in the process
we drop some directions, but under suitable assumptions we may
continue as long as we are still left with some exit
directions. This is pretty vague statement, but hopefully it will
become more clear as we progress.
}

Our setting is motivated by Proposition~\ref{con:long-transition}.
For the sake of completeness we recall here from
Section~\ref{sec:NonTransverseDiffusion} some assumptions.

Let $n_i >0$ for $i=1,\dots,L$  and let $n_1+n_2 + \dots + n_L=n$ and $w=w_u + w_s$, where $w_u,w_s \in \mathbb{N}$.

For $i=1,\dots,L$ let $V_i$ be a subspace with $\dim V_i=n_i$.  Let $W_u$ and $W_s$ be two subspaces of dimensions $w_u$ and $w_s$, respectively.

In $\mathbb{R}^n\times \mathbb{R}^{w_u}\times \mathbb{R}^{w_s}$ we will represent points as
$(z_1,\dots,z_{L},z_{L+1},z_{L+2})$, where $z_i \in V_i$ for $i=1,\dots,
L$, $z_{L+1} \in W_u$ and $z_{L+2} \in W_s$. In each of the spaces $V_i$, $W_j$ we have some fixed basis and an isomorphism with some $\mathbb{R}^d$ equipped with the metric, so we can define  balls
in this subspace.

First, we  put sections (hyperplanes of codimension $n_i$) in the vicinity of each
point $p_i$ (not to be confused with the Poincar\'e sections for ODEs): 
 the exit section  $S_i$ is given by conditions
 \begin{eqnarray}
 z_{i}=\delta_i=(\Delta_{i},0^{n_i-1}) \in V_i, \quad i=1,\dots,L,
 \end{eqnarray}
  where $\Delta_i >0$.

We also define
\begin{equation}
  \delta_{L+1}=0.
\end{equation}

For $1 \leq i\leq L+1$ we define the set $M_i$ (centered on the section $S_i$ for $i \leq L$ and on $p_{L+1}$ for $M_{L+1}$)
\begin{eqnarray}
  M_i=p_i + \delta_i +  \left(\Pi_{j=1}^L \overline{B}_{n_j}(0,t_{i,j})\right) \times \overline{B}_{w_u} (0,t_{i,L+1}) \times \overline{B}_{w_s} (0,t_{i,L+2}),   \label{eq:M_i}
\end{eqnarray}
where $t_{i,j}$ are positive real numbers.

We equip the set $M_i$ with two  different h-set structures. To
define the first one, denoted by $M_i$ for $i=1,\dots,L,L+1$,  we
declare $V_i \oplus V_{i+1} \oplus \dots \oplus V_L \oplus W_u$ as the exit directions
and  $V_1 \oplus V_2 \oplus \cdots V_{i-1} \oplus W_s$ as the entry directions. For the second one, we set
\begin{equation}
\widetilde{M}_i=R_{V_i}(M_i), \quad i=1,\dots,L.
\end{equation}
This means that (see Definition~\ref{def:dropping-exit-dim}) we declare $V_{i+1}
\oplus \dots \oplus V_L \oplus W_u$  as the exit directions  (i.e. when compared to $M_i$ we
drop the subspace $V_i$ of exit directions). We will not need $\widetilde{M}_{L+1}$.

 Observe that $\widetilde{M}_{L}$ and $M_{L+1}$ have $W_u$ as the exit directions  and $M_{L}$ has
 $V_L \oplus W_u$ as the exit directions.

Now we assume that
\begin{eqnarray}
  \widetilde{M}_i \cover{f^{l_i}} M_{i+1} \quad i=1,\dots,L. \label{eq:cov-M-tildeM}
\end{eqnarray}
The above covering  relations are expected by combining the transition where we drop the connection
direction (plus some others we decide to treat from that point on as
the entry ones) with the local hyperbolic behavior near
$p_{i+1}$, whereas in other directions for both
covering relations where the dynamics might not help us we just
adjust the sizes to obtain the correct inequalities. For this purpose we need to increase the sizes
 during the transition if these are treated as the entry directions or to decrease
 the sizes if they are treated as the exit ones.


\subsection{ The  main topological shadowing theorem}
\label{subsec:topApproach}

About the same time as this work was under development a theorem about shadowing a chain covering relations with decreasing number of exit directions appeared
in works \cite{BM+,WBS}, where a slightly different technique of proof was used, but  it still is based on the same covering relations we are using.

\begin{theorem}
\label{thm:top-chn-long-transition}
Assume that the following covering relations are satisfied
\begin{eqnarray}
  \widetilde{M}_i=R_{V_i}(M_i) &\cover{f^{l_i}} & M_{i+1} \quad i=1,\dots,L. \label{eq:cov-M-M-top}
\end{eqnarray}
Let $k_i=\sum_{j=1}^i l_j$.

Then there exists $q$ such that
\begin{eqnarray*}
  q &\in& M_1  , \\
  f^{k_i}(q) &\in& M_{i+1}  \quad i=1,\dots,L.
\end{eqnarray*}
\end{theorem}

\noindent
\textbf{Proof:}
Equation (\ref{eq:M_i}) allows us to introduce the coordinates on $M_i$ through the map
\begin{eqnarray}
  \mathcal{C}_i : M_i \to  \Pi_{i=1}^{L} \overline{B_{V_i}}(0,1) \times \overline{B_{W_u}}(0,1) \times \overline{B_{W_s}}(0,1), \label{eq:Ci-coord} \\
  \mathcal{C}_i(z_1,\dots,z_L,z_{L+1},z_{L+2})=\left( \frac{z_1 - p_{i,1} - \delta_{i,1}}{t_{i,1}}, \dots , \frac{z_{L+2}- p_{i,L+2} - \delta_{i,L+2}}{ t_{i,L+2}}\right). \nonumber
\end{eqnarray}
Observe that the above coordinates $\mathcal{C}_i$, up to a permutation required to put the exit direction first, are the ones from the natural structure of h-set.

From now on we will use these coordinates. Without any loss of generality  we will assume that $M_i=M_{i,c}=\Pi_{i=1}^{L} \overline{B_{V_i}}(0,1) \times \overline{B_{W_u}}(0,1) \times \overline{B_{W_s}}(0,1)$.

We will prove the following statement, which implies the assertion of our theorem.

 For any
$\bar{y} \in B_{W_s}(0,1)$, $\bar{x} \in B_{W_u}(0,1)$, $\eta_i \in B_{V_i}(0,1)$, $i=1,\dots,L$,
there exists $q$ such that
\begin{eqnarray}
  q &\in& M_1, \ \pi_{W_s}(q)=\bar{y}, \ \pi_{V_1}(q)=\eta_1, \label{eq:topthm1}\\
  f^{k_i}(q) &\in&  M_{i+1}, \ \pi_{V_{i+1}} f^{k_i}(q) = \eta_{i+1}  \quad i=1,\dots,L, \label{eq:topthm2} \\
  \pi_{W_u}(f^{k_{L+1}}(q))&=&\bar{x}. \label{eq:topthm3}
\end{eqnarray}

In the sequel we will denote $f^{l_i}$ by $f_i$.

To obtain $q_1 \in  M_1$ satisfying (\ref{eq:topthm1}--\ref{eq:topthm3}) it is enough to find a sequence $\{q_i\}_{i=1}^{L+1}$ satisfying the following conditions
\begin{eqnarray}
  y(q_1) - \bar{y} &=& 0,  \label{eq:y10} \\
  z_{i}(q_i) - \eta_i &=& 0, \quad i=1,\dots,L  \label{eq:q-on-sec} \\
  f_{i}(q_i) - q_{i+1}&=&0, \quad i=1,\dots,L, \label{eq:q-traj} \\
  z_{L+1}(q_{L+1}) - \bar{x} &=& 0.  \label{eq:q-last-p}
\end{eqnarray}
which we will consider in the set
\begin{eqnarray*}
D=\Pi_{i=1}^{L+1} M_i.
\end{eqnarray*}
Let us remind the reader that the supports of $M_i$ and
$\widetilde{M}_i$ coincide, but $M_i^\pm$ and
$\widetilde{M}_i^\pm$ differ.

 Observe that the number of equations in system
(\ref{eq:y10}--\ref{eq:q-last-p}) coincides with the number of
variables in $D$. Indeed the equation count goes as follows:
\begin{itemize}
\item (\ref{eq:y10}) gives $w_s$ equations
\item (\ref{eq:q-on-sec}) consists of $n_1+n_2 + \dots + n_L=n$
equations
\item (\ref{eq:q-traj}) consists of $L \cdot (n+w_u + w_s)$ equations
\item (\ref{eq:q-last-p}) gives $w_u$ equations,
\end{itemize}
which gives $(L+1)(n+w_u + w_s)$ equations in the system, which coincides with the dimension of the set $D$.

If $w_u=0$, then $\bar{x}=0$ and equation (\ref{eq:q-last-p}) is dropped from further considerations when defining
maps $F$, $H_t$. Analogously, when $w_s=0$ then $\bar{y}=0$ and we drop equation (\ref{eq:y10}).

Let us denote by $F$ the map given by the left hand side of system
(\ref{eq:y10}--\ref{eq:q-last-p}). We have for $q=(q_1,\dots,q_{L+1})\in D$
\begin{equation}
  F(q)= \begin{pmatrix}
  y(q_1) - \bar{y} & \\
  z_{i}(q_i) - \eta_i & i=1,\dots,L   \\
  f_{i}(q_i) - q_{i+1} & i=1,\dots,L, \\
   z_{L+1}(q_{L+1}) - \bar{x} &
  \end{pmatrix}.
    \label{eq:defFn4}
\end{equation}

We will prove that system (\ref{eq:y10}--\ref{eq:q-last-p})  has a
solution in $D$, by using the homotopy argument to show that the
local Brouwer degree $\deg (F,\inter D,0)$ is nonzero.

Let $h_i$ for $i=1,2,\dots,L$ be the homotopies  from the covering relations
\eqref{eq:cov-M-M-top}.

We imbed $F$ into a one-parameter family of maps (a homotopy)
$H_t$ as follows
\begin{equation}
  H_t(q)= \begin{pmatrix}
 y(q_1) - (1 - t) \bar{y} &  \\
  z_{i}(q_i) - (1 - t) \eta_i   & i=1,\dots,L   \\
  h_{t,i}(q_i) - q_{i+1} &  i=1,\dots,L, \\
  z_{L+1}(q_{L+1}) - (1-t) \bar{x} &
    \end{pmatrix}
    \label{eq:defHn4}
\end{equation}
It is easy to see that  $H_0(q)=F(q)$.

We show that if $q \in \partial D$ then $H_t(q) \neq 0$ holds for all $t \in [0,1]$.
This will imply that $\deg(H_t,D,0)$ is
defined for all $t \in [0,1]$ and does not depend on $t$.

Let $q \in \partial D$. Then $q_i \in \partial M_i$ for some $i=1,\dots,L+1$.
We will use the following decomposition of $\partial M_i$ for $i=1,\dots,L$:
$\partial M_i=M_i^+ \cup (\widetilde{M}_i^+ \cap M_i^-) \cup
\widetilde{M}_i^-$, whereas for $i=L+1$, since $\widetilde{M}_{L+1}$ is not defined,
$\partial M_i=M_i^+ \cup M_i^-$. It may happen that $M_{L+1}^-=\emptyset$.

\begin{itemize}
\item{the case  $q_i \in M_i^+$}.

If  $i>1$, then we consider  $\widetilde{M}_{i-1}
\cover{f_{i-1}} M_i$ and we see from (\ref{eq:hc3}) that $q_i \notin h_{t,i}(q_{i-1})(M_{i-1})$. Therefore in this case
$ h_{t,i}(q_{i-1}) - q_i \neq 0$.

If  $i=1$, then $y(q_1) \neq (1-t)\bar{y}$, because in this case $y(q_1) \in \partial B_{w_s}(0,1)$, hence
$\|y(q_1)\|=1 > \|\bar{y}\|$.

\item{the case $i \leq L$ and $ q_i \in \widetilde{M}_i^+ \cap M_i^-$.} Then $q_i \in \partial B_{V_i}(0,1)$,
hence $\|z_i(q_i)\|=1 > \|\eta_i\|$.

\item{the case $i \leq L$ and $q_i \in \widetilde{M}_i^- $}.

From the exit condition (\ref{eq:hc2}), in the covering relation $\widetilde{M}_{i} \cover{f_{i}} M_{i+1}$ it follows $h_{t,i}(q_i) \notin M_{i+1}$, therefore
$h_{t,i}(q_i) - q_{i+1} \neq 0$.

\item{the case $i = L+1$ and $q_{L+1} \in M_{L+1}^- $}.

We have $\|z_{L+1}(q_{L+1})\|=1 > \|\bar{x}\|$.

\end{itemize}

We have proved that $\deg(H_t,\inter D,0)$ is defined. By the
homotopy invariance we have
\begin{equation}
  \deg(F,\inter D,0)=\deg(H_1,\inter D,0). \label{eq:deg-cont}
\end{equation}

In the sequel the points in $M_i$ (and $\tilde{M}_i$) will be denoted by $(z_{i,1},\dots,z_{i,L},z_{i,L+1},y_i)$, where $z_{i,k} \in V_k$ for
$k=1,\dots,L$, $z_{i,L+1} \in W_u$ and $y_i \in W_s$.

Observe that $H_1(q)=0$ is the following system of \emph{linear
equations}
\begin{eqnarray*}
  y_1=0, \\
  z_{1,1}=0,  \\
  (0,A_1(z_{1,2},\dots,z_{1,L+1}),0) -
  (z_{2,1},z_{2,2},\dots,z_{2,L+1},y_2)= 0, \\
  z_{2,2}= 0, \\
  (0,0,A_2(z_{2,3},\dots,z_{2,L+1}),0) -
  (z_{3,1},z_{3,2},\dots,z_{3,L+1},y_3)= 0, \\
  \dots \\
  (0,\dots,A_{L}(z_{L,L+1}),0) -
  (z_{L+1,1},z_{L+1,2},\dots,z_{L+1,L+1},y_{L+1})= 0, \\
   z_{L+1,L+1} = 0,
\end{eqnarray*}
where $A_i$ is a linear map which appears at
the end of the homotopy $h_i$.

It is not hard to see that $q=0$ is the only solution of this
system. For the proof observe that $y_i=0$ for $i=1,2,\dots$
because the first term in each equation involving $A_i$ has zero
on the last ($y$) coordinate. To prove that $z_{i,j}=0$ for
$i,j=1,\dots,L,L+1$, we should start from the two bottom equations to infer
that $z_{L+1,i}=0$ for $i=1,\dots,L+1$, and since $A_{L}$ is an
isomorphism then also $z_{L,L+1}=0$. Now we consider $z_{L,i}$
from the next two equations from the bottom and so on.

Therefore $\deg(H_1,\inter D,0) = \pm 1$.
This and \eqref{eq:deg-cont} implies that
\begin{equation}
  \deg(F,\inter D,0)=\pm 1,
\end{equation}
hence there exists a solution of equation $F(q)=0$ in $D$. This
finishes the proof.

 \qed

\subsection{Generalization}

In the theorem below we allow chains of coverings relations combined with dropping some directions.
\begin{theorem}
\label{thm:top-gen}

Assume that we have h-sets $N_i$ and $M_j$ (and $\tilde{M}_j$ when some exit dimensions have been dropped) and
the following covering relations are satisfied
\begin{eqnarray*}
  N_{0,0} \cover{f_{0,0}} N_{0,1} \cover{f_{0,1}} &\cdots& \cover{f_{0,i_0}} N_{0,i_0+1}=M_0, \\
  \tilde{M}_{0}=  N_{1,0} \cover{f_{1,0}} N_{1,1} \cover{f_{1,1}} &\cdots& \cover{f_{1,i_1}} N_{1,i_1+1}=M_1, \\
   \tilde{M}_{1}=  N_{2,0} \cover{f_{2,0}} N_{2,1} \cover{f_{2,1}} &\cdots& \cover{f_{2,i_1}} N_{2,i_2+1}=M_2, \\
   &\dots& \\
   \tilde{M}_{L}=  N_{L,0} \cover{f_{L,0}} N_{L,1} \cover{f_{L,1}} &\cdots& \cover{f_{L,i_1}} N_{L,i_L+1}=M_L.
\end{eqnarray*}

Then there exists $q_0,\dots,q_l$, such that
\begin{eqnarray*}
  q_k &\in& N_{k,0},  \quad f_{k,j}\circ \cdots \circ f_{k,1} \circ f_{k,0}(q_k) \in N_{k,i_j+1}, \quad  j=0,\dots,i_k, \quad k=0,\dots,L \\
  q_{k+1}&=& f_{k,i_k}\circ \cdots \circ f_{k,1}\circ f_{k,0}(q_k), \quad k=0,\dots,L-1.
\end{eqnarray*}

\end{theorem}
\textbf{Proof:} Conceptually the same as the proof of Theorem~\ref{thm:top-chn-long-transition}.
\qed

\section{Diffusion in the linear model}
\label{sec:linModel-proof}
We prove now Proposition~\ref{conj:very-simple-model}  for $f$ being a linear model.
 To formulate the precise assumptions about our linear model we need first to introduce some notations.

Let $z=(x_1,\dots,x_n)$. For $i=0,\dots,n$, define $z_i=(z_{i,p},z_{i,\inc},z_{i,\out},z_{i,f})$ where
\begin{itemize}
\item $z_{i,p}=(x_1,\dots,x_{i-1})$ are the past coordinates
\item $z_{i,\inc}=x_i$ is the incoming coordinate
\item $z_{i,\out}=x_{i+1}$ is the outgoing coordinate
\item $z_{i,f}=(x_{i+2},\dots,x_n)$ are the future coordinates.
\end{itemize}
These are the local coordinates around each fixed point $p_i$.

We assume that we have a sequence of linear maps: $f_i$ for $i=0,\dots,n$ and affine maps
$f_{i-1,i}$ for $i=1,\dots,n$.

We assume that the map $f$ is equal to the linear map $f_i$ around
the fixed point $p_i$ and defined as $f_{i-1,i}$  close to the
heteroclinic connection. Therefore will refer to $f_i$' as
\emph{the local maps} and $f_{i-1,i}$'s will be called
\emph{the transition maps}.

\subsection{Local maps  }

Let
$f_i(z_i)=\left(f_{i,p}(z_i),f_{i,\inc}(z_i),f_{i,\out}(z_i),f_{i,f}(z_i)\right)$
the decomposition of the map $f_i$ in terms of the previous
splitting of the coordinates $z_i$. We  assume that:
\begin{align}
\label{LinearCondition}
f_{i,p}(z_i)&=A_{i,p} z_{i,p}\notag\\
f_{i,\inc}(z_i)&=\mu_iz_{i,\inc}\\
f_{i,\out}(z_i)&=\lambda_iz_{i,\out}\notag\\
f_{i,f}(z_i)&=A_{i,f} z_{i,f},\notag
\end{align}
where $A_{i,p}$ and $A_{i,f}$ are matrices  that satisfy
\begin{eqnarray*}
\left|A_{i,p}z_{i,p}\right|\leq\mu_{i,p}|z_{i,p}|\qquad\left|A_{i,f}z_{i,f}\right|\geq\lambda_{i,f}|z_{i,f}|,
 \end{eqnarray*}
 with $|\mu_i|, \mu_{i,p}<1$ and $1<\lambda_{i,f},|\lambda_i|$. The norm that we are using here and
 for the rest of the proof is the maximum norm, $|.|=||.||_\infty$.


Let us fix $\epsilon >0$, $0 <\sigma < \epsilon$ and $0<\eta<1$
for all $i=0,\dots,n$.

For each $i=0,\dots,n$, we want to define $h$-sets that will be centered in the following points $q_{i,\inc}$ and $q_{i,\out}$:
\begin{itemize}
\item $q_{i,\inc}=(0,\sigma,0,0)$
\item $q_{i,\out}=(0,0,\sigma,0)$.
\end{itemize}
 Notice that $q_{i,\inc}$ is located close to the fixed point
$p_i$ in the direction of the incoming heteroclinic connection, defined by
the segment $C_i$. The point $q_{i,\out}$ is also located close to
the fixed point $p_i$, but in the direction of the outgoing
heteroclinic, defined by the segment $C_{i+1}$. Define the sets:
\begin{eqnarray}
N_i^\inc&=&\{z_i\in\mathbb{R}^{n}\,:\,\left|z_i-q_{i,\inc}\right|\leq \epsilon\}\label{Niinc}\\
N_i^\out&=&\left\{z_i\in\mathbb{R}^{n}\,:\,\left| z_{i,p}\right| \leq (1-\eta)\epsilon, \,
\left|z_{i,inc}\right|  \leq (1-\eta)\epsilon, \right. \nonumber \\
 & &\left.        \left|z_{i,out}- \sigma \right|\leq (1+\eta)\epsilon, \,
 \left|z_{i,f}\right|\leq (1+\eta)\epsilon \right\}.\label{Niout}
\end{eqnarray}
 We equip these sets with the natural $h$-set structure (see Definition~\ref{def:hset-nat-stru}).
 We declare the directions $(z_{i,p},z_{i,\inc})$ as the entry directions and
$(z_{i,\out},z_{i,f})$ as the exit directions in both cases.

\begin{lemma}
\label{lem:lmod-cov-fp} For every $i=0,\dots,n-1$, there exists an
integer $k_i$ such that the following covering relation holds:
$$N_i^\inc \cover{f_i^{k_i}}  N_i^\out.$$
\end{lemma}
\textbf{Proof:}
Since the map is linear we only have to prove that the entry (stable) components
of $N_i^\inc$ are mapped inside $N_i^\out$ and that the exit  (unstable) directions
of $N_i^\inc$ cover the exit components of $N_i^\out$. This is, the boundary of
the exit directions of $N_i^\inc$ is mapped outside $N_i^\out$.

Let us start with the past components. We have to show that
$\left|f^{k_i}_{i,p}(z_i)\right|<(1-\eta)\epsilon$ for $z_{i}\in
N_i^\inc$. But
\begin{equation*}
\left|f^{k_i}_{i,p}(z_i)\right|=\left|A^{k_i}_{i,p} z_{i,p}\right|\leq
 \mu^{k_i}_{i,p}\left|z_{i,p}\right|\leq\mu^{k_i}_{i,p}\epsilon,
 \end{equation*}
 and the requested inequality holds for
 \begin{equation}
\label{kifirst}
   k_i > \frac{\ln (1-\eta)}{\ln \mu_{i,p}}.
 \end{equation}

Consider the incoming component, $z_{i,\inc}$. We want $k_i$ such
that $f^{k_i}_{i,\inc}(z_i) < (1-\eta)\epsilon$ for $z_{i}\in
N_i^\inc$. But
\begin{equation*}
\left|f^{k_i}_{i,\inc}(z_i)\right|=\mu_i^{k_i}\left|z_{i,\inc}\right|\leq
   \mu_i^{k_i}\left(\sigma+\epsilon\right).
\end{equation*}
If we take
\begin{equation}
\label{kisecond}
  k_i>\frac{\ln{\frac{\sigma+\epsilon}{(1-\eta)\epsilon}}}{\ln{\mu_i^{-1}}},
\end{equation}
 we obtain the desired inequality.

Now we study the exit components. Take $z_{i}\in N_i^\inc$ such
that its outgoing component $z_{i,\out}$ satisfies
$\left|z_{i,\out}\right|=\epsilon$. We want to see that
$$\left|f_{i,\out}^{k_i}(z_i)\right| > \sigma+(1+\eta)\epsilon.$$
 Notice that we have:
\begin{equation*}
\left|f_{i,\out}^{k_i}(z_i)\right|=\lambda_i^{k_i}\left|z_{i,\out}\right|=\lambda_i^{k_i}\epsilon.
\end{equation*}
If we take $k_i$ such that
\begin{equation}
\label{kithird}
k_i>\frac{\ln{\frac{\sigma+(1+\eta)\epsilon}{\epsilon}}}{\ln{\lambda_i}},
\end{equation}
we obtain the desired inequality.

Finally, for the future components we proceed in the same way.
Take $z_{i}\in N_i^\inc$ such that its future component $z_{i,f}$
satisfies $\left|z_{i,f}\right|=\epsilon$. We want to see that
\[\left|f_{i,f}^{k_i}(z_i)\right| > (1+\eta)\epsilon.\] But
\[\left|f_{i,f}^{k_i}(z_i)\right|=\left|A_{i,f}^{k_i}z_{i,f}\right|\geq
\lambda_{i,f}^{k_i}\left|z_{i,f}\right|=\lambda_{i,f}^{k_i}\epsilon,\]
and the requested inequality holds for when
\begin{equation}
\label{kifourth}
  k_i > \frac{\ln (1+\eta)}{\ln \lambda_{i,f}}.
\end{equation}

To finish the proof of Lemma~\ref{lem:lmod-cov-fp}, we take $k_i$ large enough to satisfy the
derived above lower bounds \eqref{kifirst}-\eqref{kifourth} for $k_i$.
\qed

\subsection{Dropping of one direction}

Now we are going to equip $N_i^{\out}$ with another $h$-set structure,
$\widetilde{N_i^{\out}}$. We are going to put the outgoing coordinate
$z_{i,\out}$ in the set of entry directions.
Notice that $\widetilde{N_i^{\out}}$ is the same as $N_i^{\out}$ as sets.
We are only changing the declaration of entry and exit coordinates, that is,
the $h$-set structure.

Notice that it is precisely at this moment where we drop
the outgoing direction. This argument is equivalent to the one in
Section \ref{subsubsec:sketch} where we intersected some domain with a section
of co-dimension one located in the desired outgoing direction.
Notice that $\widetilde{N_i^{\out}}$ have the same number of entry
(and exit) components than $N_{i+1}^{\inc}$.

\subsection{Transition along the heteroclinic connection}
 We define the map close to the heteroclinic segment just as a translation,
 $f_{i,i+1}$. For points in $\widetilde{N_i^{\out}}$, that is
 for points of the form $q_{i,\out}+z_i$
 the map $f_{i,i+1}$ is defined as:
\begin{equation}
\label{TranslationCondition}
f_{i,i+1}(q_{i,\out}+z_i)=q_{i+1,\inc}+z_i.
\end{equation}
Notice that, with the transition written in this way we do not have to perform
a change of variables that would locate the fixed point $p_{i+1}$ at the origin.
The change is included in the transition.

Our goal is to prove that $\widetilde{N_i^{\out}}$ covers $N_{i+1}^{\inc}$.
If we write the transition map in terms of $z_i$ and $z_{i+1}$ we have:
\begin{align*}
(z_{i,p},z_{i,\inc})&=z_{i+1,p}\\
z_{i,\out}&=z_{i+1,\inc}\\
z_{i,f}&=(z_{i+1,\out},z_{i+1,f}).
\end{align*}

With this relation, taking into account the relative sizes of the
entry and exit directions in $\widetilde{N_i^{\out}}$ and
$N_{i+1}^\inc$ we can conclude that:
\begin{lemma}
\label{lem:cv-tran-map}
The following covering relation hold:
$$\widetilde{N_i^{\out}} \cover{f_{i,i+1}}  N_{i+1}^\inc,$$ for all $i=0,\dots,n-1$.
\end{lemma}

\subsection{The conclusion}

By combining Lemmas~\ref{lem:lmod-cov-fp},\ref{lem:cv-tran-map} with Theorem~\ref{thm:top-gen} we obtain the following
\begin{theorem}
\label{thm:lin-model}
 Under the previous assumptions \eqref{LinearCondition}, \eqref{TranslationCondition},
 for all $\epsilon>0$ there exists a point $x$ and
 a sequence of integers $0=k_0<k_1<\dots<k_n$ such that:
 $$||f^{k_i}(x)-p_i||<\epsilon \quad i=0,\dots,n.$$
\end{theorem}

\section{Diffusion in a simplified Toy Model from \cite{CK}}
\label{sec:diffToymodel}

\subsection{The toy model from \cite{CK}}
The toy model system from \cite[page 59, eq. (31)]{CK}
is given by
\begin{equation}
  \frac{d}{dt} b_j = -i |b_j|^2 b_j + 2i \overline{b}_j(b_{j-1}^2 + b_{j+1}^2), \quad j \in \mathbb{Z} \label{eq:toy-model}
\end{equation}
where $b_j \in \mathbb{C}$, $j \in \mathbb{Z}$, $i=\sqrt{-1}$, and
$\overline{z}$ denote the complex conjugation of $z$.

System (\ref{eq:toy-model}) is Hamiltonian with
\begin{equation}
  \partial_t b_j = -2i \frac{\partial H}{\partial \overline{b}_j},
  \quad  \partial_t \overline{b}_j = 2i \frac{\partial H}{\partial
  b_j},
\end{equation}
where
\begin{equation}
  H(b)= \sum_j \left( \frac{1}{4} |b_j|^4 - \mbox{Re} (\overline{b}_j^2 b^2_{j-1})\right)=  \sum_j \left( \frac{1}{4}|b_j|^4 - \frac{1}{2} \overline{b}_j^2 b^2_{j-1}-
   \frac{1}{2} b_j^2 \overline{b}^2_{j-1}
  \right)
\end{equation}

Another conserved quantity for \eqref{eq:toy-model} called
\emph{the mass} is given by
\begin{equation}
  M(b)= \sum_j |b_j|^2.
\end{equation}

Following \cite{CK,GK} we are interested in the dynamics
(\ref{eq:toy-model}) on the hypersurface $M(b)=1$, and more specifically
in the diffusing solutions, which transport the mass
from modes with $j$ small to modes with large $j$. This phenomenon
is one of the main ingredients in the transfer of energy to high
frequencies in the cubic defocusing NLS on a 2D torus established in \cite{CK,GK}.

Observe from (\ref{eq:toy-model}) that on the hypersurface $M(b)=1$ there exists  the following family of periodic
solutions indexed by $j \in \mathbb{Z}$: $b_j(t)= \textrm{e}^{it}$ and $b_k(t)=0$ for $k \neq j$.
Following \cite{CK,GK} we denote the $j$-th
orbit in this family by $\mathbb{T}_j$. It turns out there exist heteroclinic
connections from $\mathbb{T}_j$ to $T_{j\pm 1}$. The diffusing
orbits constructed in \cite{CK,GK} follow the chain of
heteroclinic connections $\mathbb{T}_0 \to \mathbb{T}_1 \to \dots
\to \mathbb{T}_N$ for arbitrary $N \in \mathbb{N}$.

\subsubsection{Coordinates from \cite{CK} and the local form of the toy model}

In the analysis of (\ref{eq:toy-model}) in \cite{CK} (see also \cite{GK})  new coordinates were introduced,
which turn each  periodic orbit $\mathbb{T}_j$ into a fixed
point.

These coordinates are defined as follows.

We fix $j$ and we introduce new coordinates $r,\theta,c_k$ for $k \neq j$
\begin{equation}
  b_j=re^{i\theta}, \qquad b_k=c_k e^{i \theta}.
\end{equation}
Let
\begin{equation}
  \omega=e^{i 2\pi/3}=\left(\frac{-1}{2},\frac{\sqrt{3}}{2}\right).
\end{equation}
For $c \in \mathbb{C}$ we set
\begin{equation}
  c=\omega c^- + \omega^2 c^+,  \label{eq:defc+c-}
\end{equation}
where $c^-,c^+ \in \mathbb{R}$.

The above decomposition means that we represent a complex number $c$
in the basis $\{\omega,\omega^2\}$ over the field $\mathbb{R}$.

Then we introduce  $c_{j\pm 1}^\pm \in \mathbb{R}$ by
\begin{equation}
  c_{j\pm 1}=\omega c^-_{j \pm 1} + \omega^2 c^+_{j \pm 1}. \label{eq:cjpm1}
\end{equation}

The form of the toy model system in these coordinate is given by the following lemma.

\begin{lemma}\cite[Prop. 3.1, page 69]{CK}
\label{lem:eq-CK-Tj}
Consider the toy model with the constraint $M(b)=1$. Let
$c=(c_k)_{k \neq j}$, $c_*=(c_k)_{k \neq j-1,j,j+1}$.

Then the equations for the toy model system have the following form
\begin{eqnarray}
   \dot{c}_{j-1}^-&=& -\sqrt{3} c_{j-1}^- + O_1(c^2 c_{j-1}^-) + O_2(c_{k \neq j-1}^2 c_{j-1}^+), \label{eq:CK-der-cjm1-}  \\
   \dot{c}_{j-1}^+&=& \sqrt{3} c_{j-1}^+ + O_3(c^2 c_{j-1}^+) + O_4(c_{k \neq j-1}^2 c_{j-1}^-), \label{eq:CK-der-cjm1+}  \\
   \dot{c}_{j+1}^-&=& -\sqrt{3} c_{j+1}^- + O_1(c^2 c_{j+1}^-) + O_2(c_{k \neq j+1}^2 c_{j+1}^+), \label{eq:CK-der-cjp1-}  \\
   \dot{c}_{j+1}^+&=& \sqrt{3} c_{j+1}^+ + O_3(c^2 c_{j+1}^+) + O_4(c_{k \neq j+1}^2 c_{j+1}^-), \label{eq:CK-der-cjp1+}  \\
   \dot{c_k}&=& i c_k + O(c^2 c_k), \quad k \neq j\pm 1, j,  \label{eq:CK-der-ck}
\end{eqnarray}
where $c_{k \neq j-1}=(c_k)_{k \neq j, k \neq j-1}$
Moreover, all these $O_?$ are uniform with respect to $k$ and $j$
\end{lemma}

\subsubsection{The transition between consecutive charts}
\label{subsubsec:tran-chart}
In the coordinate systems `centered' at the $j$-th torus we have
\begin{eqnarray}
  r_j^2&=&1 - \sum_{k \neq j} |b_k|^2= 1 - \sum_{k \neq j} |c_k|^2=1 - \sum_{k \neq j} r_k^2, \\
  b_j&=&r_j e^{i \theta_j}, \quad b_k=c_k e^{i \theta_j}, k \neq j  \\
  c_{j+1}&=&\omega c_{j+1}^- + \omega^2 c_{j+1}^+, \quad  c_{j-1}=\omega c_{j-1}^- + \omega^2 c_{j-1}^+.
\end{eqnarray}

We look for the relation between coordinates in the $j$-th and $(j+1)$-th charts.
First we will prove the following simple lemma.
\begin{lemma}
If $c=\omega c_- + \omega^2 c_+$, then
\begin{equation*}
  \frac{1}{c}= \frac{c_+}{|c|^2} \omega + \frac{c_-}{|c|^2} \omega^2.
\end{equation*}
\end{lemma}
\textbf{Proof:}
It is easy to see that
\begin{eqnarray*}
 \overline{c_- \omega + c_+ \omega^2}=c_- \omega^2 + c_+ \omega.
\end{eqnarray*}
Since $c^{-1}=\frac{\overline{c}}{|c|^2}$, the assertion follows.
\qed

We will denote by the variable with tilde the coordinates expressed in the $(j+1)$-th chart.

Since
\begin{eqnarray}
  r_{j+1}&=&|c_{j+1}|=\left( (c_{j+1}^-)^2 + (c_{j+1}^+)^2 - c_{j+1}^- c_{j+1}^+ \right)^{1/2}, \\
   e^{i \theta_{j+1}}&=& \frac{c_{j+1}}{|c_{j+1}|} e^{i \theta_j},  \\
   e^{i (\theta_j - \theta_{j+1})}&=&\frac{|c_{j+1}|}{c_{j+1}}= \frac{c_{j+1}^+}{|c_{j+1}|}\omega+ \frac{c_{j+1}^-}{|c_{j+1}|}\omega^2
\end{eqnarray}
we therefore obtain
\begin{eqnarray}
 \tilde{c}_k &=& c_k e^{i(\theta_j - \theta_{j+1})}=c_k \cdot  \frac{|c_{j+1}|}{c_{j+1}}, \qquad  k < j-1 \  \mbox{or} \ k > j+2, \\
 \tilde{c}_{j-1} &=& c_{j-1} \cdot  \frac{|c_{j+1}|}{c_{j+1}} = \left(c_{j-1}^- \omega  + c_{j-1}^+ \omega^2\right) \frac{|c_{j+1}|}{c_{j+1}}, \\
 \tilde{c}_j&=&r_j e^{i(\theta_j - \theta_{j+1})} = r_j \frac{|c_{j+1}|}{c_{j+1}} =
  r_j \left(\frac{c_{j+1}^+}{|c_{j+1}|}\omega+ \frac{c_{j+1}^-}{|c_{j+1}|}\omega^2\right) \\
  \tilde{c}^-_j&=& \frac{r_j c_{j+1}^+}{|c_{j+1}|}   \quad   \tilde{c}^+_j= \frac{r_j c_{j+1}^-}{|c_{j+1}|} \\
  \tilde{c}_{j+2}&=& c_{j+2} \frac{|c_{j+1}|}{c_{j+1}} = \tilde{c}_{j+2}^- \omega +  \tilde{c}_{j+2}^+ \omega^2.
\end{eqnarray}

\subsection{Simplified toy model}
\label{subsec:c-r-O-full-model}

Our model is a simplification of the toy model system from
\cite{CK,GK}. However, while it is simpler analytically than the original, the problems and obstacles to be overcome to prove
the existence of diffusive orbits are very similar.

Our phase space is defined by a sequence of coordinate charts indexed by $j \in \mathbb{Z}$.
 Each of these maps  has at its center a fixed point $\mathbb{T}_j$.

The $j$-th node centered chart (we will sometimes refer to it as $j$-th chart) uses the following coordinates
\begin{itemize}
\item $c_k \in \mathbb{C}$, $k\leq j-2$ or $k\geq j+2$
\item $y_-,x_-,y_+,x_+ \in \mathbb{R}$.
\end{itemize}

The coordinates $y_\pm$ and $x_\pm$ are related to the ones from the previous subsection as follows
\begin{equation*}
  y_-=c_{j-1}^-, \quad x_-=c_{j-1}^+, \quad  y_+=c_{j+1}^-, \quad x_+=c_{j+1}^+.
\end{equation*}

Let us fix $\sigma>0$. For example we can take $\sigma=0.5$ or $\sigma=0.1$.
We assume that our system preserves the mass $\sum_j |c_j|^2$.

\subsubsection{Local evolution close to a fixed point}
 The evolution in the $j$-th  chart is given by the following ODE
\begin{eqnarray}
  \dot{y}_-&=& -y_- + O(x_- (y_+)^2),  \label{eq:sys-res-0-full-xy-1}  \\
  \dot{x}_-&=&  x_- + O(y_- (x_+)^2),   \\
  \dot{y}_+&=& -y_+ + O(x_+ (y_-)^2),  \\
  \dot{x}_+&=&  x_+ + O(y_+ (x_-)^2), \label{eq:sys-res-0-full-xy-4} \\
  \dot{c_k}&=& i c_k(1+O(c^2)), \qquad k\leq j-2 \quad \mbox{or} \quad k\geq j+2, \label{eq:sys-res-0-full-center}
\end{eqnarray}
where we assume that all $O()$ terms satisfy
\begin{equation}
  |O(z)| \leq K |z|.  \label{eq:O-full-K-bound}
\end{equation}
Let us stress that each $O(\dots)$ function may depend on all the variables.

Comparing with (\ref{eq:CK-der-cjm1-}--\ref{eq:CK-der-ck}), we make the following simplifications:
 \begin{itemize}
 \item we have left only the resonant terms in the first four equations.   Removing the non-resonant terms makes  perfect sense, as this can be achieved by a suitable coordinate change (see \cite{GK}).
 \item $\sqrt{3}$ was replaced by $1$ in the first four equations (only for simplicity of the notation).
\end{itemize}

For a given $T>0$ we will denote by $\varphi_T$ a shift by time $T$ along a trajectory of (\ref{eq:sys-res-0-full-xy-1}--\ref{eq:sys-res-0-full-center})
in the domain of the $j$-th chart.

\subsubsection{The dynamics near a heteroclinic connection}

In the domain of the $j$-th chart we assume that system  (\ref{eq:sys-res-0-full-xy-1}--\ref{eq:sys-res-0-full-center}) governs the dynamics, and we will use
$\varphi_T$ to map the neighborhood of the entry section $y_-=\sigma$ to the neighborhood of the exit section $x_+=\sigma$. Next we compose this local map with the global map from the domain of the $j$-th chart to the domain of the $(j+1)$-th chart.

 Below we postulate the precise form of this ``jump between charts'' (compare with subsection~\ref{subsubsec:tran-chart}),
 which should correspond to the movement along the heteroclinic connection between $\mathbb{T}_j$ and
$\mathbb{T}_{j+1}$, followed by the change of coordinates to the $(j+1)$-th chart.

In the formulas given below
the variables without tildes are the ones referring to the $j$-th node chart,
while those  with tildes  denote the variables with respect to the $(j+1)$-th chart.
\begin{eqnarray*}
  \tilde{c}_{k \leq j-2}&=& c_{k \leq j-2} \\
  \tilde{c}_{j-1}&=& g_1(x_-,y_-) \\
  \tilde{x}_- &=& y_+ \\
  \tilde{y}_- &=& x_+  \\
  (\tilde{x}_+,\tilde{y}_+) &=& g_2(c_{j+2}) \\
  \tilde{c}_{k \geq j+3} &=& c_{k \geq j+3},
\end{eqnarray*}
where
$g_1: \mathbb{R}^2 \to \mathbb{C}$ and $g_2: \mathbb{C} \to \mathbb{R}^2$ are defined  by
\begin{equation}
  g_1(c^-,c^+)=\omega c^- + \omega^2 c^+, \qquad g_2(c)=g_1^{-1}(c).
  \label{eq:gtwo}
\end{equation}

We will denote by $J$ the map defined by the above equations. In principle the map $J$ depends on $j$, but it will be always clear from the context
what is the $j$ used.

When comparing with the toy model in \cite{CK} we basically do the following simplification:
\begin{itemize}
\item we assume that entry section $y_-=\sigma$ in the $j+1$-th chart and  the exit section $x_+=\sigma$ in the $j$-th chart coincide. Therefore there is no need
 to study the evolution along the heteroclinic connection between these sections, which the authors in \cite{CK} were forced to do.
\item this simplification  makes perfect sense as in \cite{CK} this transition map was also studied using the system derived in Lemma~\ref{lem:eq-CK-Tj}, which the local system in our approach. Therefore in our model the transition along heteroclinic is achieved by the computation of the local map 
    with the entry and exit sections located at the macroscopic distance from the $\mathbb{T}_j$'s. 
\end{itemize}

\subsection{Fixed points and heteroclinic connections}

Observe that in our system we have fixed points $\mathbb{T}_j$ parameterized by $j \in \mathbb{Z}$
given in the coordinates centered in the $j$-th node by
\begin{equation}
  c_{k\leq j-2}=0, \quad c_{k\geq j+2}=0, \quad x_-=y_-=x_+=y_+=0.
\end{equation}
For the fixed point $\mathbb{T}_{j}$ the directions $ c_{k\leq j-2}, c_{k\geq j+2}$ are the center directions, $x_-,x_+$ are the unstable directions and
$y_-,y_+$ are the stable directions.

Consecutive fixed points $\mathbb{T}_{j}$ and $\mathbb{T}_{j+1}$
are connected by a heteroclinic connection escaping the
neighborhood of $\mathbb{T}_j$ along the solution $c_{k}(t)=0$ for
$k\leq j-2$ or $k\geq j+2$, $x_-(t)=y_-(t)=y_+(t)=0$,
$x_+(t)=\textrm{e}^{t}$  and continued later in the
coordinates centered on $\mathbb{T}_{j+1}$ as $c_{j}(t)=0$ for
$k\leq j-1$ or $k\geq j+3$, $x_-(t)=x_+(t)=y_+(t)=0$,
$y_-(t)=\textrm{e}^{-t}$.

\begin{figure}[H]
\label{fig:tms-heteroclinic}
 \centering
    \includegraphics
    {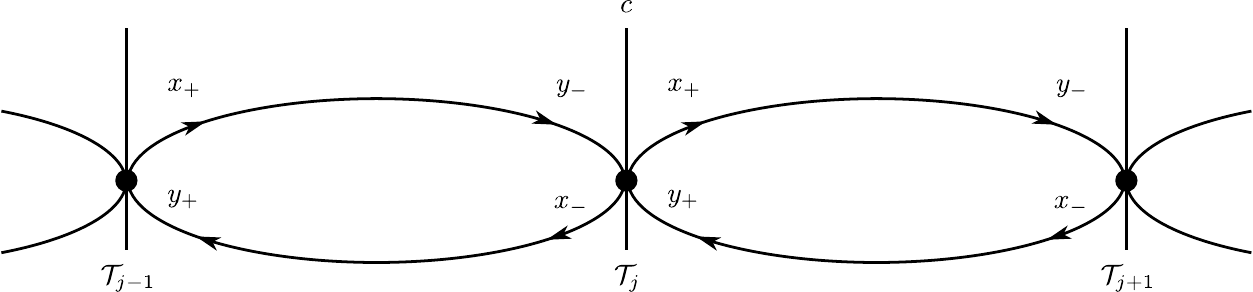}%
\caption{The heteroclinic connections in the toy model}
\end{figure}

There are more heteroclinic connections in the toy model
\eqref{eq:toy-model} see Figure~\ref{fig:tms-heteroclinic},
however we will only use the connection $\mathbb{T}_j \to \mathbb{T}_{j+1}$ described above.

\subsection{Some heuristic discussion}
\label{subsec:tms-heuristic}

Following the heteroclinic chain $\mathbb{T}_0 \to \mathbb{T}_1 \to \cdots \to \mathbb{T}_N$ requires the analysis of consecutive compositions of maps $\varphi_T$ (the  shift along the trajectory by $T$ in the $j$-th chart) with $J$, which maps from the $j$-th chart to the $(j+1)$-th chart.

After passing by the $j$-th torus $\mathbb{T}_j$ we have to drop the $x_+$-coordinate in the $j$-th chart
and then after the application of the $J$ map, this coordinate
becomes $y_-$ and it is related  to $c_j$. Therefore we decide  that in the vicinity of $\mathbb{T}_j$ we want $c_k$ with $k < j$ (the past modes) to be entry variables and $c_k$ with $k>j$ (the future modes) to be the exit directions.

The realization of the above idea must  take into account the following issues.
\begin{itemize}
\item we declare $x_-$ to be an entry direction, despite being  unstable.
\item in the center directions we can have also some growth or decay which might be non-desirable for the past modes and for the future modes, respectively. This growth or decay might turn out to be unbounded for $T \to \infty$.
\end{itemize}

To overcome these difficulties we have to chose carefully the relative sizes of the variables.
For heuristic reasons we will deal with three different sizes: macro, micro and nano.
 The macro size will be $O(1)$, the micro size $w(T)e^{-T}$ and the nano size $w(t) e^{-2T}$, where $w$ is some polynomial. By increasing $T$ we will be able
 to follow the heteroclinic chain as close as we desire, keeping the relative sizes of the variables.
 This idea appears also in \cite{CK}.

In the evolution related to passing by $\mathbb{T}_j$ in the $j$-th chart  on the exit section $x_+=\sigma$
(which is of macro size), we want all the other variables of micro size, i.e.
$w(T)e^{-T}$, where $w(T)$ is some polynomial.  To get this sizes on the exit section $x_+=\sigma$,
we will need to impose on the entry section $y_-=\sigma$ (the macro size) to the variable $x_-$ to be of nano size $w(T)e^{-2T}$,
whereas all the other variables should be of micro size.

 Regarding the behavior of the center directions, it turns out that it is possible to maintain the micro size throughout this transition for the variables $c_k$ for $k\leq j-2$ or $k \geq j+2$, because  the possible decay or growth is bounded by a constant, which does not depend on $T$. This is the content of Theorem~\ref{thm:center-estm}.

\subsection{Evolution estimates in the hyperbolic directions}
\label{subsec:scatt-estm} We consider the system
(\ref{eq:sys-res-0-full-xy-1}--\ref{eq:sys-res-0-full-center}).
Let
\begin{equation}
  \sigma'=1.01 \cdot \sigma\ . \label{eq:sigma'}
\end{equation}

\begin{theorem}
\label{thm:est-res-full-O}
Consider (\ref{eq:sys-res-0-full-xy-1}--\ref{eq:sys-res-0-full-center}) satisfying (\ref{eq:O-full-K-bound}) with initial conditions
\begin{equation}
  y_-(0)=\eta, \quad x_-(0)=e^{-2T} a_0, \quad y_+(0)=e^{-T} b_0, \quad x_+(0)=e^{-T}d_0  \label{eq:sys-res-full-ic}
\end{equation}
where  $\eta \in (0,\sigma')$ and $a_0 \in a$, $b_0\in b$, $d_0 \in d$ where  $a,b,d$ are real intervals such that
\begin{equation}
  a,b,d \subset [-T^k,T^k], \quad k \geq 1,  \label{eq:sys-res-full-ic-bounds}
\end{equation}
and
\begin{equation}
\sum_{l}|c_l(0)|^2 = 1.  \label{eq:scat-estm-c0}
\end{equation}

Assume  that $T \geq T_0>1$ is large enough, so that  the following inequalities are satisfied
\begin{eqnarray}
   e^T &\geq& (K (T^k + T   2K \sigma'  \left(T^k + \sigma' \right)^2 )  \left(T^k+ T 4 K \sigma'^2 (T^k + \sigma' )  \right)^2 \\
   e^T &\geq& K \left(T^k + T 4 K \sigma'^2 (T^k + \sigma' )  \right) (T^k + T   2K \sigma'  \left(T^k + \sigma' \right)^2  )^2 \\
    T e^{-3T} &<& \sigma'.  \label{eq:Te-tsigma}
\end{eqnarray}
Then for $t \in [0,T]$
\begin{eqnarray*}
   y_-(t)&\in & \eta e^{-t} + \alpha e^{-t} t e^{-3T} \\
  x_-(t)&\in& e^{-2T} a_0 e^t + \alpha t e^{-2T} e^t \left(  2K \sigma'  \left(|d| + \sigma' \right)^2  \right) \\
  y_+(t)&\in&e^{-T}b_0 e^{-t}+ \alpha t e^{-T}e^{-t} \left(4 K \sigma'^2 (|d|+ \sigma' ) \right) \\
  x_+(t)&\in &e^{-T} d_0 e^t + \alpha t e^{-4T}e^{t},
\end{eqnarray*}
where
\begin{equation*}
  \alpha=[-1,1].
\end{equation*}
\end{theorem}

In Theorem~\ref{thm:est-res-full-O} it is implicitly assumed that we are working in the $j_0$-th chart.
This index does not appear neither in the statement of the result nor in the proof, because
Theorem~\ref{thm:est-res-full-O} is concerned only about the hyperbolic directions.
Nevertheless, we will complete this result with the estimates in the center directions
in Theorem~\ref{thm:center-estm} of the next subsection,
and there the use of $j_0$ will be required.

Assumption (\ref{eq:scat-estm-c0}) implies an a-priori bound for all variables so we have a uniform estimate on all $O()$ terms in equations
(\ref{eq:sys-res-0-full-xy-1}--\ref{eq:sys-res-0-full-center}).

The proof can be obtained by a direct verification based on differential inequalities  and the continuation argument (as in \cite{CK}). However we prefer to do it in a more constructive  way
to show how these bounds have been derived. The proof uses iterations, which produce the stabilizing estimates after the third iterate.
The  result of the $l$-th iteration  step will be called the $l$-th approximation.

In the estimates we will repeatedly  use the following simple lemma.
\begin{lemma}
\label{lem:simple-estm}
 Consider the following one-dimensional ODE
 \begin{equation*}
   x'=\lambda x + e^{\lambda t} D(t), \quad x(0)=x_0
 \end{equation*}
 where the continuous function $D: \mathbb{R} \to \mathbb{R}$ satisfies
 \begin{equation*}
    |D(t)| \leq D, \quad 0 \leq t \leq T
 \end{equation*}
 Then
 \begin{equation*}
    x(t) \in e^{\lambda t} x_0 + [-1,1] \cdot D t e^{\lambda t}, \quad 0 \leq t \leq T.
 \end{equation*}
\end{lemma}

\subsubsection{First approximation}

We begin our process with
\begin{eqnarray*}
  y_-(t)&=&\eta e^{-t} \\
  x_-(t)&=&a_0 e^{-2T}  e^t \\
  y_+(t)&=&b_0 e^{-T} e^{-t} \\
  x_+(t)&=&d_0 e^{-T}  e^t.
\end{eqnarray*}

The estimates for the nonlinear terms (we skip the absolute value sign for $a,b,d$ in the estimates for the nonlinear terms) are
\begin{eqnarray*}
  K |x_-(t)(y_+(t))^2|  \leq  K (e^{-2T} a e^t)(e^{-T}b e^{-t})^2 \leq e^{-4T} e^{-t} \left(K a b^2 \right) \leq e^{-3T} e^{-t}
\end{eqnarray*}
for $T$ large enough to satisfy
\begin{equation}
  e^{T} > Ka b^2. \label{eq:res-full-cond1}
\end{equation}
\begin{eqnarray*}
  K  |y_-(t) (x_+(t))^2| &\leq& K (\eta e^{-t}) (e^{-T} d e^t)^2 \leq e^{-2T} e^t \left(K \sigma' d^2 \right), \\
  K  |x_+(t) (y_-(t))^2| &\leq& K \left( e^{-T} d e^t \right) \left( \eta e^{-t} \right)^2 \leq e^{-T}e^{-t} \left(Kd \sigma'^2 \right), \\
  K |y_+(t) (x_-(t))^2| &\leq& K \left(e^{-T}b e^{-t} \right)\left( e^{-2T} a e^t\right)^2=e^{-5T}e^{t} \left(K b a^2 \right) \leq e^{-4T}e^{t}
\end{eqnarray*}
for $T$ large enough to satisfy
\begin{equation}
  e^{T} > Kb a^2. \label{eq:res-full-cond2}
\end{equation}
We obtain the following estimates in $1$-st approximation
\begin{eqnarray*}
   y_-(t)&\in & \eta e^{-t} + \alpha e^{-t} t e^{-3T} \\
  x_-(t)&\in& e^{-2T} a_0 e^t + \alpha t e^{-2T} e^t \left(K \sigma' d^2 \right) \\
  y_+(t)&\in&e^{-T}b_0 e^{-t}+ \alpha t e^{-T}e^{-t} \left(Kd \sigma'^2 \right) \\
  x_+(t)&\in &e^{-T} d_0 e^t + \alpha t e^{-4T}e^{t}.
\end{eqnarray*}

\subsubsection{$2$-nd approximation}
Using the $1$-st approximation we obtain
\begin{eqnarray*}
  K |x_-(t)(y_+(t))^2|  &\leq&  K (e^{-2T}e^t)(a + t K \sigma' d^2)(e^{-T} e^{-t})^2 \left(b+ t Kd \sigma'^2 \right)^2 = \\
   & & e^{-4T} e^{-t} \left(K (a + t K \sigma' d^2)  \right)\left(b+ t K d \sigma'^2 \right)^2 \leq e^{-3T} e^{-t} \\
  K  |y_-(t) (x_+(t))^2| &\leq& K  e^{-t} \left(\sigma' + t e^{-3T}\right) (e^{-T}  e^t)^2 \left(d + t e^{-3T} \right)^2 = \\
     & & e^{-2T}e^t \left( K \left(\sigma' + t e^{-3T}\right) \left(d + t e^{-3T} \right)^2 \right) \\
     &\leq& e^{-2T}e^t   \left( 2K \sigma'  \left(d + \sigma' \right)^2 \right) \\
  K  |x_+(t) (y_-(t))^2| &\leq& K \left( e^{-T} e^t \right) \left(d + t e^{-3T} \right) e^{-2t}\left( \sigma' + t e^{-3T} \right)^2= \\
   & & e^{-t} e^{-T}\left(K \left(d + t e^{-3T} \right) \left( \sigma' + t e^{-3T} \right)^2  \right) \\
   &\leq& e^{-t} e^{-T} \left(4 K \sigma'^2 (d+ \sigma' ) \right) \\
  K |y_+(t) (x_-(t))^2| &\leq& K \left(e^{-T} e^{-t} \right)\left(b+ t K d \sigma'^2 \right)\left( e^{-2T} e^t\right)^2\left( a + t K \sigma'^2 d \right)^2 = \\
    & & e^{-5T} e^{t} \left( K \left(b+ t K d \sigma'^2 \right)\left( a + t K \sigma'^2 d \right)^2 \right) \leq e^{-4T} e^{t}
\end{eqnarray*}
provided $T$ is large enough for the following conditions to hold for $t \in [0,T]$
\begin{eqnarray*}
  e^{T} &>& \left(K (a + t K \sigma' d^2)  \right)\left(b+ t K d \sigma'^2 \right)^2, \\
  t e^{-3T} &<& \sigma',  \\
  e^T &>& \left(K \left(b+ t K d \sigma'^2 \right)\left( a + t K \sigma'^2 d \right)^2 \right).
\end{eqnarray*}

These bounds give us the following estimates in $2$-nd
approximation
\begin{eqnarray*}
   y_-(t)&\in & \eta e^{-t} + \alpha e^{-t} t e^{-3T} \\
  x_-(t)&\in& e^{-2T} a_0 e^t + \alpha t e^{-2T} e^t \left(  2K \sigma'  \left(d + \sigma' \right)^2  \right) \\
  y_+(t)&\in&e^{-T}b_0 e^{-t}+ \alpha t e^{-T}e^{-t} \left(4 K \sigma'^2 (d+ \sigma' ) \right) \\
  x_+(t)&\in &e^{-T} d_0 e^t + \alpha t e^{-4T}e^{t}.
\end{eqnarray*}
Observe that we have obtained the same formula for $y_-(t)$ and $x_+(t)$ in the second approximation as in the first approximation.

\subsubsection{Next iterate, the $3$-rd approximation}
Since the bounds  for the terms $y_- (x_+)^2$ and $x_+ (y_-)^2$ will be the same as in the previous subsubsection,   we just compute bounds for $x_- (y_+)^2$ and $y_+ (x_-)^2$.

We have
\begin{eqnarray*}
  K |x_-(t)(y_+(t))^2|  \leq \\
   K (e^{-2T}e^t)(a + t   2K \sigma'  \left(d + \sigma' \right)^2  )(e^{-T} e^{-t})^2 \left(b+ t 4 K \sigma'^2 (d+ \sigma' )  \right)^2 = \\
   e^{-4T} e^{-t} \left(K (a + t   2K \sigma'  \left(d + \sigma' \right)^2 )  \left(b+ t 4 K \sigma'^2 (d+ \sigma' )  \right)^2   \right) \leq e^{-3T} e^{-t} \\
  K |y_+(t) (x_-(t))^2| \leq \\
   K \left(e^{-T} e^{-t} \right)\left(b+ t 4 K \sigma'^2 (d+ \sigma' )  \right) \left( e^{-2T} e^t\right)^2 (a + t   2K \sigma'  \left(d + \sigma' \right)^2  )^2 = \\
    e^{-5T} e^{t} \left( K \left(b+ t 4 K \sigma'^2 (d+ \sigma' )  \right) (a + t   2K \sigma'  \left(d + \sigma' \right)^2  )^2 \right) \leq e^{-4T} e^{t}
\end{eqnarray*}
provided $T$ is large enough for the following conditions to hold for $t \in [0,T]$
\begin{eqnarray*}
  e^T &\geq& (K (a + t   2K \sigma'  \left(d + \sigma' \right)^2 )  \left(b+ t 4 K \sigma'^2 (d+ \sigma' )  \right)^2 \\
   e^T &\geq& K \left(b+ t 4 K \sigma'^2 (d+ \sigma' )  \right) (a + t   2K \sigma'  \left(d + \sigma' \right)^2  )^2
\end{eqnarray*}
This concludes the proof of Theorem~\ref{thm:est-res-full-O}.

\subsection{The estimates in the center direction}
\label{subsec:tms-center-estm}
The goal of this subsection to complete the estimates for the flow near $\mathbb{T}_j$ given in Theorem~\ref{thm:est-res-full-O} by providing bounds for $c_j$'s
for $j \leq j_0-2$ or $j \geq j_0+2$, where $j_0$ is the index of the torus on which our chart is centered.

Let $G=3/2$ be such that for $c=c_- \omega + c_+ \omega^2$ the following estimate holds
\begin{eqnarray*}
 |c|^2 \leq G(|c_-|^2 + |c_+|^2).
\end{eqnarray*}

\begin{theorem}
\label{thm:center-estm}
Under the same assumptions as in Theorem~\ref{thm:est-res-full-O}, take
\begin{equation}
  d = [-2.5 \sigma', 2.5 \sigma'] \label{def:d}
\end{equation}
and assume that $c_j=0$ for $j<0$ and $j>N$.

Let $j_0$ be the index of the chart used and assume that
\begin{equation}
  c_j(0)=u_j e^{-T} \in ue^{-T}, \quad u \subset [-T^k,T^k].  \label{eq:center-ic-size}
\end{equation}
for the indexes $j=0,\dots,N$ such that $j \leq j_0-2$ or $j \geq j_0+2$

Assume also that  $T$ is large enough to satisfy
\begin{eqnarray}
    G e^{-2T} \left(\left(T^k + T (2 K \sigma' (3.5 \cdot \sigma')^2)\right)^2  +
    \left(T^k + T (4K \sigma' (3.5 \cdot \sigma'))\right)^2 \right) + \nonumber \\
     2  N  \exp\left(21 \cdot  K G  \sigma'^2 \right)  T^{2k+1} e^{-2T} <  G \sigma'^2.\label{eq:BoundOnT}
\end{eqnarray}
Then
\begin{equation}
   |c_j(0)|^2 \exp(-21 \cdot K G \sigma^2)   < |c_j(t)|^2 < |c_j(0)|^2 \exp(21 \cdot K G \sigma^2). \label{eq:cj-lup-estm}
\end{equation}
\end{theorem}
\textbf{Proof:}
Since $\dot{c}_j=ic_j(1+ g(t))$, where $g(t)=O(c^2(t))$, we obtain
\begin{eqnarray*}
  \frac{d}{dt}|c_j(t)|^2=\dot{c}_j \overline{c}_j + c_j \dot{\overline{c}}_j= - 2 |c_j|^2 \im g(t).
\end{eqnarray*}
Therefore we have
\begin{equation*}
  -2|c_j(t)|^2 \cdot |g(t)| \leq \frac{d}{dt}|c_j(t)|^2 \leq 2|c_j(t)|^2 |g(t)|,
\end{equation*}
hence
\begin{equation}
|c_j(0)|^2 \exp\left(-2 \int_0^t |g(s)|ds\right)  \leq  |c_j(t)|^2 \leq |c_j(0)|^2 \exp\left(2 \int_0^t |g(s)|ds\right). \label{eq:cj-low-upp-estm}
\end{equation}

Therefore, it will important to show that under our assumptions on  the initial conditions  we  have a good bound
for $\int_0^T \sum_{j \neq j_0}|c_j|^2(s) ds$.
As the main step to achieve this goal we will show that there exists $B>1$ such that
\begin{equation}
   |c_j(t)| \in \alpha B |u_j| e^{-T}  \label{eq:cj-apriori}
\end{equation}
holds for all $t \in [0,T]$.
We use a continuation argument. Since $B>1$ then (\ref{eq:cj-apriori}) is satisfied for $t \in [0,T']$. We will show that $T' \geq T$.

We have the following estimate
\begin{eqnarray*}
  \sum_{j \neq j_0}|c_j(t)|^2 \leq    G \left(|x_-(t)|^2 + |y_-(t)|^2 + |x_+(t)|^2 + |y_+(t)|^2 \right) + \sum_{j \neq j_0,j_0\pm 1}|c_j(t)|^2
\end{eqnarray*}

From (\ref{eq:cj-apriori}) we have that, for $j \neq j_0,j_0\pm 1$ and $t \in [0,T']$,
\begin{eqnarray}
  \int_0^t |c_j(s)|^2 ds  \leq \int_0^t B^2 |u_j|^2 e^{-2T}ds \leq B^2 |u_j|^2 T e^{-2T}. \label{eq:int-cj}
\end{eqnarray}

To estimate $\int_0^t |y_-(s)|^2ds$ we use the bounds obtained in Theorem~\ref{thm:est-res-full-O}.
From assumption~(\ref{eq:Te-tsigma}) of Theorem~\ref{thm:est-res-full-O} we get
\begin{equation*}
  |y_-(t)| \leq 2 \sigma' e^{-t}
\end{equation*}
and therefore
\begin{eqnarray}
  \int_0^{T'} |y_-(t)|^2 ds \leq 4 \sigma'^2 \int_0^\infty e^{-2s}ds = 2 \sigma'^2. \label{eq:int-y-}
\end{eqnarray}

To estimate $\int_0^t |x_-(s)|^2 ds$  we use again Theorem~\ref{thm:est-res-full-O}. From
\begin{equation*}
  |x_-(t)| \leq e^t e^{-2T} \left(|a| + T (2 K \sigma' (|d| + \sigma')^2) \right),
\end{equation*}
we get
\begin{equation}
  \int_0^{T'} |x_-(s)|^2 ds
    < 0.5 \cdot e^{-2T} \left(|a| + T (2 K \sigma' (|d| + \sigma')^2)\right)^2. \label{eq:int-x-}
\end{equation}

We now estimate  $\int_0^t |y_+(s)|^2 ds$. From Theorem~\ref{thm:est-res-full-O} we have that
\begin{equation*}
  |y_+(t)| \leq e^{-T}e^{-t} \left(|b| + T (4K \sigma' (|d| + \sigma')) \right),
\end{equation*}
hence
\begin{equation}
  \int_0^{T'} |y_+(s)|^2 ds
<0.5 \cdot e^{-2T} \left(|b| + T (4K \sigma' (|d| + \sigma'))\right)^2. \label{eq:int-y+}
\end{equation}

It remains to estimate $\int_0^t |x_+(s)|^2 ds$. By assumption~(\ref{eq:Te-tsigma}) and the definition~(\ref{def:d}) of $d$
we get
\begin{equation*}
  |x_+(t)| \leq 3.5 \sigma' e^{-T} e^t,
\end{equation*}
so that
\begin{eqnarray}
  \int_0^{T'} |x_+(s)|^2 ds \leq 3.5^2 \sigma'^2 e^{-2T} \int_0^T e^{2s} ds <8 \cdot \sigma'^2. \label{eq:int-x+}
\end{eqnarray}

By combining (\ref{eq:int-cj},\ref{eq:int-y-},\ref{eq:int-x-},\ref{eq:int-y+},\ref{eq:int-x+}) we directly obtain
\begin{eqnarray*}
  \int_0^{T'} \sum_{j \neq j_0}|c_j(s)|^2 ds \leq N B^2  |u|^2 T e^{-2T} +  G \cdot 2 \cdot \sigma^2 + G \cdot 8 \cdot \sigma'^2 \\
   + G \cdot \left( 0.5 \cdot e^{-2T} \left(|a| + T (2 K \sigma' (|d| + \sigma')^2)\right)^2 \right) + \\
  G \cdot \left(  0.5 \cdot e^{-2T} \left(|b| + T (4K \sigma' (|d| + \sigma'))\right)^2 \right)  < \\
  N B^2 |u|^2 T e^{-2T} + 10 G \sigma'^2  + \\
  0.5 \cdot G e^{-2T} \left(\left(|a| + T (2 K \sigma' (|d| + \sigma')^2)\right)^2  +
    \left(|b| + T (4K \sigma' (|d| + \sigma'))\right)^2 \right)
\end{eqnarray*}

Introducing
\begin{eqnarray*}
  E&=& N B^2 |u|^2 T e^{-2T}  + \\
  & & 0.5 \cdot G e^{-2T} \left(\left(|a| + T (2 K \sigma' (|d| + \sigma')^2)\right)^2  +
    \left(|b| + T (4K \sigma' (|d| + \sigma'))\right)^2 \right),
\end{eqnarray*}
and taking into account~(\ref{eq:O-full-K-bound}) we get
\begin{equation*}
  \int_0^{T'} |g(s)|ds \leq K E + 10 K G \sigma'^2.
\end{equation*}
For the continuation argument~(\ref{eq:cj-apriori}), by~(\ref{eq:center-ic-size}) and~(\ref{eq:cj-low-upp-estm}), we need that
\begin{equation*}
  |u_j|^2 e^{-2T} \exp\left( 2K E + 20 K G \sigma'^2\right) < |u_j|^2 B^2 e^{-2T},
\end{equation*}
which is equivalent to
\begin{equation*}
   \exp\left( 2K E  + 20  K G \sigma'^2 \right) < B^2.
\end{equation*}
It is clear that this can be achieved if
we take $B$ such that
\begin{equation}
  \exp\left(21 \cdot  K G  \sigma'^2 \right) = B^2 \label{eq:central-B2}
\end{equation}
and then take $T$ large enough to have
\begin{eqnarray*}
 G e^{-2T} \left(\left(|a| + T (2 K \sigma' (|d| + \sigma')^2)\right)^2  +
    \left(|b| + T (4K \sigma' (|d| + \sigma'))\right)^2 \right) + \\
     2  N B^2 |u|^2 T e^{-2T} <  G \sigma'^2,
\end{eqnarray*}
which is guaranteed by hypothesis~(\ref{eq:BoundOnT}) on $T$.

We have also obtained estimate~(\ref{eq:cj-lup-estm}):
\begin{equation*}
   |c_j(0)|^2 \exp(-21 \cdot K G \sigma'^2)   < |c_j(t)|^2 < |c_j(0)|^2 \exp(21 \cdot K G \sigma'^2).
\end{equation*}
\qed

\subsection{Construction of the covering relations}
The goal of this section is to construct a sequence of covering relations for our model using the estimates
obtained in Theorems~\ref{thm:est-res-full-O} and~\ref{thm:center-estm}.

Let us set (compare the estimate for the expansion and contraction rates in the center direction (\ref{eq:cj-lup-estm}) in Theorem~\ref{thm:center-estm})
\begin{equation}
  A=\exp(21 \cdot K G \sigma'^2).  \label{eq:def-A}
\end{equation}

We will have two types of $h$-sets in the $j$-th chart, $N^j_{in}$ and $N^j_{out}$ (to be defined later in this section), such that
 the following covering relations are satisfied
\begin{eqnarray}
  N^j_{in} \cover{\varphi_T} N^j_{out},  \label{eq:tms-nin-nout}\\
  R_{<x_+,y_+>} N^{j}_{out}=\widetilde{N}^j_{out} \cover{J} N^{j+1}_{in}. \label{eq:tms-tran}
\end{eqnarray}
In  relation (\ref{eq:tms-nin-nout}) the map $\varphi_T$ is the shift along the trajectory by the time $T$. The map $J$ in relation (\ref{eq:tms-tran}) is the jump modeling
the transition along the heteroclinic and the change of coordinates. Recall (see Subsection~\ref{subsec:NM-def}) that $R_{<x_+,y_+>} N^{j}_{out}$ means that we drop
the directions $x_+$ and $y_+$, because we have  just `passed' by $T_j$ in our heteroclinic chain.

In the derivation we will use the following conventions
\begin{itemize}
\item $\gamma(variable)$ - will be used for the sizes in the entry directions,
\item $r(variable)$ - will be used for the sizes in the exit directions.
\end{itemize}
To be more precise, by the size will mean the radius or half-diameter of the balls or intervals used to define our h-sets.

Observe that when we are dropping some directions, the sizes in these directions become very close to zero (to set them to zero will not change anything, but it would require slight changes in Theorem~\ref{thm:top-gen}).

By $c_p$ (the past modes) we will denote the collection $\{c_k\}_{k \leq j-2}$ and by $c_f$ (the future modes) we will denote the collection $\{c_k\}_{k \geq j+2}$. On
$c_p$ and $c_f$ we use the sup norm, i.e. $\|c_p\|=\sup_{k \leq j-2} |c_k|$.

The structure of h-sets $N^j_{in}$ and $N^j_{out}$ is defined as follows (we are using the $j$-th chart):
\begin{itemize}
\item  the entry variables: $c_p$, $x_-$, $y_-$
\item  the exit variables: $x_+$, $y_+$, $c_f$
\item  parameters of $N^j_{in}$:

For the entry directions:
\begin{eqnarray*}
 |c_p| &\leq& \gamma^j_{in}(c_p) e^{-T}, \quad \mbox{(micro)} \\
  y_-&\in &\sigma + \alpha \gamma^j_{in}(y_-)e^{-T}, \qquad \gamma^j_{in}(y_-)\approx 0,  \quad \mbox{(macro)} \\
  |x_-|&\leq& \gamma^j_{in}(x_-)e^{-2T}, \qquad  \gamma^j_{in}(x_-)\approx 0.  \quad \mbox{(nano)}
\end{eqnarray*}

For the exit directions:
\begin{eqnarray}
  |x_+| &\leq& r^j_{in}(x_+) e^{-T}, \quad r^j_{in}(x_+)=2.1 \sigma,  \quad \mbox{(micro)} \label{eq:rjinplus}\\
  |y_+| &\leq& r^j_{in}(y_+) e^{-T},  \quad \mbox{(micro)} \nonumber \\
  |c_f| &\leq& r^j_{in}(c_f)e^{-T}.  \quad \mbox{(micro)} \nonumber
\end{eqnarray}

\item  parameters of $N^j_{out}$:

For the entry directions:
\begin{eqnarray*}
  |y_-|&\leq&\gamma^j_{out}(y_-)e^{-T},  \quad \mbox{(micro)} \\
  |x_-|&\leq&\gamma^j_{out}(x_-)e^{-T},  \quad \mbox{(micro)} \\
  |c_p| &\leq& \gamma^j_{out}(c_p) e^{-T}.  \quad \mbox{(micro)}
\end{eqnarray*}

For the exit directions:
\begin{eqnarray*}
  x_+ &\in& \sigma + \alpha r^j_{out}(x_+)e^{-T}, \quad r^j_{out}(x_+)\approx 0,  \quad \mbox{(macro)}\\
  |y_+| &\leq& r^j_{out}(y_+) e^{-2T}, \quad r^j_{out}(y_+)\approx 0,  \quad \mbox{(nano)} \\
  |c_f| &\leq& r^j_{out}(c_f)e^{-T}.  \quad \mbox{(micro)}
\end{eqnarray*}
\end{itemize}

\subsubsection{Covering $N_{in}^j \cover{\varphi_T} N_{out}^j$}

We use Theorems~\ref{thm:est-res-full-O} and~\ref{thm:center-estm}  for the shift along the trajectory by time $T$ with
\begin{eqnarray}
\eta &=& \sigma + \alpha \gamma^j_{in}(y_-)e^{-T}, \quad a=\alpha \gamma^j_{in}(x_-), \nonumber \\
  b&=&\alpha r^{j}_{in}(y_+), \quad d = \alpha r^j_{in}(x_+).  \label{eq:abc-res-full-O}
\end{eqnarray}

To satisfy the assumptions of Theorem~\ref{thm:est-res-full-O} about the size of $|y_-|$ we require that
\begin{equation}
  \gamma^j_{in}(y_-)e^{-T} < \sigma' - \sigma=0.01 \sigma.
\end{equation}

 The conditions for $N_{in}^j \cover{\varphi_T} N_{out}^j$ are
\begin{itemize}
\item entry conditions:
\begin{itemize}
\item
for $c_p$ variables:
\begin{equation}
  A \gamma^j_{in}(c_p) \leq \gamma_{out}^j(c_p).  \label{eq:tms-cp-flow}
\end{equation}
\item for $x_-$ from Theorem~\ref{thm:est-res-full-O} and (\ref{eq:abc-res-full-O}) we have the following condition
\begin{equation*}
   \gamma^j_{in}(x_-)e^{-T} + T e^{-T}(2K\sigma'(3.1 \sigma')^2) < \gamma_{out}^j(x_-) e^{-T}
\end{equation*}
therefore it is enough to take
\begin{equation}
  T (2.1 K\sigma'(3.1 \sigma')^2) \leq \gamma_{out}^j(x_-), \label{eq:tms-gamma-out-x-}
\end{equation}
if $ \gamma^j_{in}(x_-) \approx 0$, which will turn out to be compatible with other conditions.
In fact we had replaced $2$ by $2.1$ and $\sigma$ by $\sigma'$ ($\sigma$ appears in $d$) in order to make an explicit margin for  $\gamma^j_{in}(x_-)$ given by
\begin{equation}
  \gamma^j_{in}(x_-) <  T (0.1 K\sigma'(3.1 \sigma')^2).  \label{eq:tms-gamma-in-x-}
\end{equation}
\item for $y_-$ it is enough to have
   \begin{equation*}
   \sigma' e^{-T} + T e^{-4T} < \gamma_{out}^j(y_-) e^{-T},
   \end{equation*}
and from (\ref{eq:Te-tsigma}) it is enough to take
   \begin{equation}
     2 \sigma' \leq \gamma_{out}^j(y_-). \label{eq:tms-gamma-out-y-}
   \end{equation}
\end{itemize}
\item exit conditions:
\begin{itemize}
\item{$x_+$}
\begin{eqnarray*}
  r_{in}^j(x_+) - T e^{-3T} > \sigma + r^j_{out}(x_+)e^{-T},
\end{eqnarray*}
which in view of  assumption (\ref{eq:Te-tsigma}),(\ref{eq:sigma'}) and (\ref{eq:abc-res-full-O}) is satisfied , if
\begin{equation}
  r^j_{out}(x_+) e^{-T} < 0.09 \sigma.  \label{eq:tms-r-out-x+}
\end{equation}
Obviously  (\ref{eq:tms-r-out-x+}) is  compatible with
$r^j_{out}(x_+) \approx 0$, which is to be expected as this is the direction which will be dropped in the next covering relation.

\item{$y_+$},
 from Theorem~\ref{thm:est-res-full-O}  and (\ref{eq:abc-res-full-O}) it follows that the following estimate is sufficient
 \begin{equation}
  r_{in}^j(y_+) e^{-2T} - Te^{-2T}\left(4K \sigma'^2(r_{in}^j(x_+) + \sigma') \right) > r^j_{out}(y_+)e^{-2T}.
\end{equation}
Observe that this is the direction which is dropped in the next covering relation, hence any $r^j_{out}(y_+)  >0$, $r^j_{out}(y_+) \approx 0$ is good for our construction.
Therefore we can take (where we used also the known value~(\ref{eq:rjinplus}) of $r_{in}^j(x_+)$)
\begin{equation}
  r_{in}^j(y_+) \geq T\left(4.1 K \sigma'^2(3.1 \sigma') \right), \label{eq:tms-r-in-y+}
\end{equation}
which leaves some margin for $r^j_{out}(y_+)$, given by
\begin{equation}
  r^j_{out}(y_+) <  T\left(0.1 K \sigma'^2(3.1 \sigma') \right). \label{eq:tms-r-out-y+}
\end{equation}

\item{$c_f$}
\begin{equation}
  A^{-1} r^j_{in}(c_f) \geq r^j_{out}(c_f). \label{eq:tms-r-flow-cf}
\end{equation}
\end{itemize}
\end{itemize}

\subsubsection{Covering relation (\ref{eq:tms-tran})}

Let $L \geq 1$ be the Lipschitz constant which holds for both functions $g_1$ and $g_2$ introduced
in~(\ref{eq:gtwo}),
where on $\mathbb{R}^2$ the max-norm is used, whereas in $\mathbb{C}$ we use the euclidian norm.

The conditions are as follows.

In the entry directions
\begin{eqnarray}
  L(\gamma^j_{out}(y_-) + \gamma^j_{out}(x_-)) &\leq& \gamma_{in}^{i+1}(c_p), \label{eq:tms-tran-cov-cp-1} \\
  \gamma_{out}^{j}(c_p) &<& \gamma_{in}^{j+1}(c_p), \label{eq:cp-i-cond}\\
  r^j_{out}(x_+) &<& \gamma^{j+1}_{in}(y_-), \label{eq:tms-tran-cov-entry-x+} \\
  r^j_{out}(y_+) &<& \gamma^{j+1}_{in}(x_-). \label{eq:tms-tran-cov-entry-y+}
\end{eqnarray}

In the exit directions
\begin{eqnarray}
  r^j_{out}(c_f) &\geq& L r^{j+1}_{in}(x_+), \label{eq:tms-rcf-tran-x+} \\
  r^j_{out}(c_f) &\geq& L r^{j+1}_{in}(y_+), \label{eq:tms-rcf-tran-y+}\\
  r^j_{out}(c_f) &>&  r^{j+1}_{in}(c_{f}).  \label{eq:rcf-i}
\end{eqnarray}

\subsubsection{Solving the inequalities for coverings}

We have to find the following set of parameters
$\gamma_{in,out}^j(c_p,x_-,y_-)$ and \newline $r_{in,out}^j(x_+,y_+,c_f)$, such that
inequalities~(\ref{eq:tms-tran-cov-cp-1}-\ref{eq:rcf-i}) are satisfied.

We split these parameters into two groups: the ones related to the dropped directions
$\gamma^j_{in}(y_-)$, $\gamma^j_{in}(x_-)$, $r^j_{out}(x_+)$, $r^j(y_+)$,
 and the remaining ones.

In the first group we effectively should obtain
\begin{eqnarray*}
  \gamma^j_{in}(y_-)=0 , \quad  \gamma^j_{in}(x_-)=0,  \quad r^j_{out}(x_+)=0, \quad r^j_{out}(y_+)=0.
\end{eqnarray*}
The conditions involving these parameters are (\ref{eq:tms-gamma-in-x-}), (\ref{eq:tms-r-out-x+}), (\ref{eq:tms-r-out-y+}), (\ref{eq:tms-tran-cov-entry-x+}) and (\ref{eq:tms-tran-cov-entry-y+}). It is clear that these conditions can be easily satisfied with all these parameters being very close to zero.

Now we deal with the other directions. We already have set the value for $r^j_{in}(x_+)$ in~(\ref{eq:rjinplus})
and we now set the following parameters  (compare with (\ref{eq:tms-gamma-out-x-}),(\ref{eq:tms-gamma-out-y-}), (\ref{eq:tms-r-in-y+}))
\begin{eqnarray}
   r^j_{in}(x_+)&=&2.1 \sigma,\nonumber \\
   \gamma_{out}^j(x_-)&=& T Q_1, \quad Q_1=(2.1 K\sigma'(3.1 \sigma')^2) , \\
  \gamma^j_{out}(y_-)&=&2 \sigma', \\
   r_{in}^j(y_+) &=& T Q_2, \quad   Q_2=\left(4.1 K \sigma'^2(3.1 \sigma') \right),
\end{eqnarray}

The remaining inequalities involve only the sizes for the variables $c_p$ and $c_f$. These are as follows:
\begin{itemize}
\item for the entry directions (see (\ref{eq:tms-cp-flow}), (\ref{eq:tms-tran-cov-cp-1}), (\ref{eq:cp-i-cond}))
\begin{eqnarray*}
  A \gamma^j_{in}(c_p) &\leq& \gamma_{out}^j(c_p) \\
    L(2 \sigma' + T Q_1) &\leq& \gamma_{in}^{j+1}(c_p), \\
  \gamma_{out}^{j}(c_p) &<& \gamma_{in}^{j+1}(c_p),
\end{eqnarray*}
\item for the exit directions (see (\ref{eq:tms-r-flow-cf}), (\ref{eq:tms-rcf-tran-x+}), (\ref{eq:tms-rcf-tran-y+}) and (\ref{eq:rcf-i})
\begin{eqnarray*}
  A^{-1}  r^j_{in}(c_f) &\geq& r^j_{out}(c_f), \\
    r^j_{out}(c_f) &\geq& 2.1 L \sigma, \\
  r^j_{out}(c_f) &\geq& T L Q_2 , \\
  r^j_{out}(c_f) &>&  r^{j+1}_{in}(c_{f}).
\end{eqnarray*}
\end{itemize}

It is clear that there exist $\gamma^j_{in}(c_p)$, $\gamma^j_{out}(c_p)$ satisfying
\begin{equation}
   Q_3 T <  \gamma^j_{in}(c_p) \leq A^{-1} \gamma_{out}^j(c_p) < \gamma_{in}^{j+1}(c_p), \quad j=0,\dots,N
\end{equation}
where $Q_3=L\left( \frac{2 \sigma'}{T} + Q_1 \right)$.

For instance we take
\begin{equation*}
  \gamma^0_{in}(c_p)=\tilde{Q_3} T, \quad \gamma^{j+1}_{in}(c_p)=\tilde{A} \gamma^j_{in}(c_p), \quad \gamma^{j}_{out}(c_p)=A\gamma^j_{in}(c_p), \quad j=0,\dots, N-1
\end{equation*}
where $\tilde{Q}_3 > Q_3$ and $\tilde{A} > A$.
With such sequence we have solved the inequalities for the entry directions.
Observe that as a consequence we have
\begin{equation}
   Q_3 T <  \gamma^j_{in}(c_p) < \gamma_{out}^j(c_p) \leq \tilde{Q}_3 \tilde{A}^N T, \qquad j=0,\dots,N. \label{eq:tms-cv-sizes-cp}
\end{equation}

In the exit direction the situation is similar. We just take any sequence satisfying
\begin{equation}
  \max (TLQ_2, 2.1 L \sigma) \leq r^{j+1}_{out}(c_f) \leq A^{-1}r^{j+1}_{in}(c_f) < A^{-1} r^j_{out}(c_f).
\end{equation}
For example we can take for $j=N,\dots,1$
\begin{equation}
  r^N_{out}(c_f)=\max (TLQ_2, 2.1 L \sigma), \quad r^{j-1}_{out}(c_f)=\tilde{A}r^{j}_{out}(c_f), \quad r^j_{in}(c_f)=Ar^j_{out}(c_f).
\end{equation}
Hence we obtain
\begin{equation}
  r^j_{out}(c_f), r^j_{in}(c_f) \leq  \tilde{A}^N  \max (TLQ_2, 2.1 L \sigma), \quad j=0,\dots,N. \label{eq:tms-cv-sizes-cf}
\end{equation}

In (\ref{eq:tms-cv-sizes-cp}) and (\ref{eq:tms-cv-sizes-cf}) we introduced also an upper bound which is $O(T)$, so now all sizes are $O(T)$ times a suitable weight function ($e^{-T}$ or $e^{-2T}$). Observe that  this bound allows us to use
Theorems~\ref{thm:est-res-full-O} and~\ref{thm:center-estm} with $k=2$, for any $T \geq \tilde{A}^N\max(2LQ_2, 2Q_3,Q_2, Q_1,1)$, because then all sizes will be less than $T^2$.

\subsection{The conclusion}
From the chain of coverings constructed above  and Theorem~\ref{thm:top-gen}  since the  distance from the heteroclinics in the chain are $O(T)e^{-T}$ we obtain the following theorem.
\begin{theorem}
\label{thm:toytoymodel}
 For the system discussed in this section for any $N$, for all $\epsilon>0$ there exists a point $x_0$ close to $\mathbb{T}_0$ whose trajectory is $\epsilon$ close
 to the  chain of heteroclinic connections $\mathbb{T}_0 \to \mathbb{T}_1 \to \cdots \to \mathbb{T}_N$.
\end{theorem}


\begin{thebibliography}{KKY}

\bibitem[Ar]{Ar64} V. I. Arnold, \emph{Instability of dynamical systems with several degrees of freedom},
   Soviet Mathematics 5(1964),581--585


\bibitem[BM+]{BM+} R. Barrio, M.A. Martinez, S. Serrano, D. Wilczak,
\emph{When chaos meets hyperchaos: 4D R{\"o}ssler model},
Physics Letters A, Vol. 379, No. 38, 2300-2305 (2015).

\bibitem[CKS+]{CK} J. Colliander, M. Keel, G. Staffilani, H. Takaoka, T. Tao, \emph{Transfer of energy to high frequencies in the cubic defocusing nonlinear
  Schr\"odinger equation}, Invent. math (2010) 181:39-113

\bibitem[C]{C} C. C. Conley
{\em Isolated Invariant Sets and the Morse Index.} 1978. CBMS vol.
38, Amer. Math. Soc., Providence

\bibitem[GK]{GK} M. Guardia, V. Kaloshin, \emph{Growth of Sobolev norms in the cubic defocusing nonlinear Schr\"odinger
  equation}, Journal of the European Mathematical Society,  vol. 17, 1, 71--149, 2015




\bibitem[MM]{MM} K. Mischaikow, M. Mrozek,
The Conley Index, in: Handbook of Dynamical Systems II: Towards Applications,
(B. Fiedler, ed.) North-Holland, 2002.


\bibitem[S]{S} J. Smoller, Shock Waves and Reaction -Diffusion Equations. Springer, 1983.







\bibitem[W2]{W2} D. Wilczak,
\emph{The existence of Shilnikov homoclinic orbits in the Michelson system: a computer assisted proof.}
Foundations of Computational Mathematics, Vol.6, No.4, 495-535, (2006).


\bibitem[WBS]{WBS} D. Wilczak, S. Serrano, R. Barrio,
\emph{Coexistence and dynamical connections between hyperchaos and chaos in the 4D R{\"o}ssler system: a Computer-assisted proof},
SIAM Journal on Applied Dynamical Systems, 15 (2016), 356--390




\bibitem[ZGi]{ZGi}P. Zgliczy\'{n}ski and M. Gidea, \emph{Covering relations
for multidimensional dynamical systems}, Journal of Differential
Equations 202/1, 32--58 (2004)


\end{thebibliography}
\end{document}